\pdfoutput=1

\documentclass[mnsc,nonblindrev]{informs_modified}
\OneAndAHalfSpacedXII

\usepackage{natbib}
 \bibpunct[, ]{(}{)}{,}{a}{}{,}%
 \def\bibfont{\small}%
\usepackage{color}
\usepackage{algorithm}
\usepackage{algorithmic}
\usepackage{setspace}

\usepackage{array}
\usepackage{multirow}
\usepackage{comment}
\usepackage[toc,page]{appendix}
\usepackage{appendix}

\usepackage{setspace}
\usepackage{hyperref}
\usepackage[utf8]{inputenc}

\TheoremsNumberedThrough     
\ECRepeatTheorems

\EquationsNumberedThrough 

\usepackage{tabularx}
\usepackage{hyperref}
\usepackage{url}

\def\eventzero{\mathcal E_0}
\def\eventzeroc{\mathcal E_0^{\mathsf C}}
\def\eventone{\mathcal E_1}
\def\eventonec{\mathcal E_1^{\mathsf C}}
\def\eventtwo{\mathcal E_2(h,\ell)}

\def\eventthree{\mathcal E_3}

\def\eventfour{\mathcal E_4}
\def\eventfourc{\mathcal E_4^{\mathsf C}}
\def\eventfive{\mathcal E_5}
\def\eventfivec{\mathcal E_5^{\mathsf C}}
\def\rh{\mathcal R_3^2}
\def\rl{\mathcal R_3^1}
\def \mzero{\sqrt{m}}
\def\rewGood{\rew_{A}\left(\mathbf n, \boldsymbol \rho;  \mathcal G\right)}
\def\rewBad{\rew_{A}\left(\mathbf n, \boldsymbol \rho; \mathcal G^{\mathsf C}\right)}
\def\rewGen{\rew_{A}\left(\mathbf n, \boldsymbol \rho\right)}

\def\hzeroExp{\min\{y\ge 0:\frac{p(m-\mzero )}{1-p}\ge y  p+\sqrt{y }\}}
\def \hzero{h_0}
\def \mOneExp{\hzero (1-p)-\sqrt{p(1-p)\hzero }- \frac{\beta}{\sqrt{\hzero }}}
\def\mOne{h_1}
\def\mTwoExp{\min_{y\ge 1/p^4}\big\{ y: r_1- r_2 >\frac{2}{y^{1/8}\sqrt{y^{1/4} p-1}}\big\}}
\def \mTwo{m_1}
\def \CROneExp{\min\left\{\left(1-\frac{1-p}{p\mzero}\right)^{+},1-\frac{ r_1- r_2}{(r_1-r_2)+r_2\mzero}\right\}} 
\def \CROne {\text{CR}_1}
\def \CRTwoExp {\min\left\{\big(1-\frac{1}{\mzero}\big)\big(1-\frac{1}{m}\big), \frac{\mOne }{m} \right\}}
\def \CRTwo{\text{CR}_2}
\def \X{W}
\def \V{V}

\def\CRThreeEXPTwo{\V \cdot \min \Bigg\{(1-\frac{1}{m^2})\X, \frac{1}{2}(1-\frac{1}{m^2})\X+\frac{1}{2}\min\left\{(1-\frac{1}{\ell_0^2})\alpha, 1-\frac{1-p}{pm}\ell_1\right\} \Bigg\}}
\def\CRThreeExpOne{(1-\frac{1}{m})^2(1-\frac{1}{(\mzero p-m^{1/4})^2})^2\X}
\def\E{\mathbb{E}}
\def\rew{\textsc{rew}}
\def\ellzeroExp{\min\{y: (1-p)y- \sqrt{y} \geq m\}}
\usepackage{multirow}

 \usepackage{enumitem}

\begin{document}

\ARTICLEAUTHORS{
\AUTHOR{Negin Golrezaei}\AFF{Sloan School of Management, Massachusetts Institute of Technology,
\EMAIL{golrezaei@mit.edu}}
\AUTHOR{Patrick Jaillet}\AFF{Department of Electrical Engineering and Computer Science, Massachusetts Institute of Technology,
\EMAIL{jaillet@mit.edu}}
\AUTHOR{Zijie Zhou}\AFF{Operations Research Center,
Massachusetts Institute of Technology, Cambridge, MA, 02139,
\EMAIL{zhou98@mit.edu}}
}

\RUNAUTHOR{Golrezaei, Jaillet, and Zhou}
\RUNTITLE{Online Resource Allocation with Samples}

\TITLE{Online Resource Allocation with Samples}

\ABSTRACT{
 We study an online resource allocation problem under uncertainty about demand and about the reward of each type of demand (agents) for the resource. Even though  dealing with demand uncertainty in resource allocation problems has been the topic of many papers in the literature, the challenge of not knowing rewards has been barely explored. The lack of knowledge about agents' rewards is inspired by the problem of allocating units of a new resource (e.g., newly developed vaccines or drugs) with unknown effectiveness/value. For such settings, we assume that we can \emph{test} the market before the allocation period starts. During the test period, we sample each agent in the market with probability $p$. We study how to optimally exploit the \emph{sample information} in our online resource allocation problem under adversarial arrival processes.  We present an asymptotically optimal algorithm that achieves $1-\Theta(1/(p\sqrt{m}))$ competitive ratio, where $m$ is the number of available units of the resource. By characterizing  an upper bound on the competitive ratio of any randomized and deterministic algorithm, we show that our competitive ratio of $1-\Theta(1/(p\sqrt{m}))$ is tight for any $p =\omega(1/\sqrt{m})$. That  asymptotic optimality is possible with sample information highlights the significant advantage of running a test period for new resources.  
We  demonstrate the efficacy of our proposed algorithm using a dataset that contains the number of COVID-19 related hospitalized patients across different age groups.   
} 

\maketitle

\section{Introduction}

In online resource allocation problems, the goal is to allocate a limited number of a given resource to heterogeneous demand/agents that arrive over time. These problems are notoriously challenging mainly because of 
 the demand uncertainty  and scarcity of the  resource.  Such problems get even more challenging for newly developed resources (e.g., new drugs, products, and services). For such a resource, the effectiveness/value of the resource for different types of agents may not be fully known. To overcome this additional challenge,  businesses, for example,  offer free products in an exchange for honest feedback \citep{free}, and pharmaceutical companies test potential treatments/drugs in human volunteers  \citep{pfizertrial}.  These practices raise the following key question: \emph{can and to what extent such feedback improve the efficiency of online resource allocation?}

 To answer this question, we consider a decision-maker who aims to allocate her available resource to two types of unit-demand agents with unknown (expected) rewards, where type $1$ has a higher expected reward than type $2$.\footnote{In Section \ref{subsec:nagent}, we discuss the case of having more than two types of agents.} The total number of agents (i.e., the market size), as well as, the number of agents of type $1$ and $2$, denoted by $h$ and $\ell$ respectively, are chosen adversarially, 
 and hence are unknown to the decision-maker. 
 Before the allocation period starts, the decision-maker tests the market by, for example, making a public announcement and offering resources for free. We assume that with probability $p\in (0,1)$, each of the $h+\ell$ agents sees and reacts to the announcement,\footnote{In Section \ref{subsec:heterop}, we allow the probability $p$ to be different for type $1$ and $2$ agents.} and gets one unit of the resource, where we assume that $p$ is known to the decision-maker. (We will discuss this assumption later in this section.)
 These agents then provide feedback on their realized reward for the resource in return. That is, we assume that all the agents in the test period will get a resource. As we will show in Section \ref{sec:limitedsample}, this assumption can be relaxed by limiting the number of available resources during the test period. The test procedure supplies the decision-maker with some information about agents' expected rewards, as well as, the size of the market for each type of agent. We refer to this information as \emph{sample information}.

 After the test period ends, the remaining agents arrive over time according to an adversarially-chosen order. For each arriving agent, the decision-maker has to make an irrevocable decision about accepting him and allocating him one unit of the resource or rejecting him. The decision-maker who has $m$ units of the resource when the allocation period starts makes acceptance/rejection decisions while being uncertain about the number, type of future agents, and their expected rewards. The decision-maker is also uncertain about which types of agents earn higher rewards upon receiving the resource.  For such a demanding setting, our goal is to design efficient resource allocation algorithms that can optimally utilize the sample information under any possible arrival sequence.  In other words, we  measure the performance of the algorithm in terms of its competitive ratio, which is the expected ratio of the reward of the algorithm to the reward of the optimal clairvoyant solution that knows the arrival sequence and the expected reward of agents in advance; 
 see Section \ref{sec:model} for the model and the formal definition of the competitive ratio. 
 
 Before presenting our contributions, we make two remarks about our model. First, our model bears resemblance to 
the proposed models in \cite{correa2021secretary, kaplan2022online} for secretary and online bipartite matching problems, respectively.  In \cite{correa2021secretary}, each of the secretaries is placed in a sample set with probability $p$, where the value of the sampled secretaries will be disclosed to the decision-maker before the decision period starts. 
In \cite{kaplan2022online} that studies an online bipartite matching problem, each agent independently will be placed in a sample set with probability $p$. While at a high level, these works seek to design algorithms that can  take advantage of samples, the nature of their considered problems is different from ours, preventing us to use their designed algorithms for our setting. We discuss the details in Section \ref{sec:related}.

Second, our model is a special case of the single-leg revenue management problem, which has been widely studied in the literature; see, e.g., \cite{littlewood1972, amaruchkul2007single, ball2009toward, gallego2009choice,ferreira2018online, jasin2015performance, hwang2021online,  golrezaei2021online}.  In all of these aforementioned works, while the decision-maker may be uncertain about the demand (i.e., the number and the order of the arrivals), the obtained reward of different types of demand upon receiving the resource is fully known to the decision-maker. This is in sharp contrast with our model in which these rewards are not known to the decision-maker, adding extra complexity to our problem. (See also Section \ref{sec:related} for a discussion about previous works on revenue management problems with limited demand information, but full knowledge of rewards.)

 \subsection{Our Contributions}
In addition to our modeling contribution, our work makes the following contributions.

\textbf{Impossibility results.} To shed light on the value of sample information, in Section \ref{sec:negative}, we consider alternative scenarios under which either the sample information is not available or the sample information is not very informative due to the lack of some additional   knowledge (e.g., the sampling probability). In all of the  considered scenarios, we show that it is not possible to design asymptotically optimal algorithms whose competitive ratio goes to one as the number of resources $m$ increases.  
 While in the first scenario, no sample information is not available (i.e., the \emph{sampling probability} $p=0$),  in the second scenario, the sample information is available but the sampling  probability $p$ is not known to the decision-maker. 
 The result for the second scenario  justifies our assumption about knowing the sampling probability; see    \citep{correa2021secretary} for a similar assumption.  Nonetheless, this assumption  can be also justified  by the fact that the outreach program is designed by the decision-maker herself, and hence $p$ can be well estimated using historical data for similar outreach programs.  
See also our case study in Section \ref{sec:simulation} where we examine the robustness of our proposed algorithm (which we will present shortly) to the lack of knowledge of the sampling probability.\footnote{In Appendix \ref{subsec:partialknown}, using numerical studies and under adversarial arrival processes,  we further test the robustness of our proposed algorithm to the lack of knowledge of $p$. There we show that when the estimation error in $p$ is $10\%$, the competitive ratio  decreases by at most $6\%$. Similar results are obtained when the estimation error in $p$ is $20\%$. }

\textbf{Asymptotically optimal protection level algorithm.} In light of our impossibility results in Section \ref{sec:negative}, we design  a simple, yet effective protection level algorithm that optimally utilizes the sample information; see Algorithm \ref{alg:ptas} in Section \ref{sec:alg}. The algorithm uses the sample information to obtain an estimate of the expected reward of each type of agents. If  the estimated reward of type $1$ is greater than that of type $2$, the algorithm  \emph{protects} type $1$ agents, otherwise type $2$ agents will be protected. In each of these cases, the algorithm uses the sample information to estimate the protection level for the type that has a higher estimated average reward. 

In Theorem \ref{thm:2}, we present the competitive ratio of our proposed algorithm for any finite value of $m$. In addition, in Proposition \ref{thm:3}, we show that our algorithm is asymptotically optimal as $m$ goes to infinity. More precisely, we show that for any $p = \omega(1/\sqrt{m})$\footnote{We scale $p$ with $m$ because of the cost of testing the market. With a large value of $m$, if the market size turns out to be large too (e.g., proportional to $m$), the decision-maker will incur a large cost during the test period, justifying our choice of $p=\omega(1/\sqrt{m})$. Nonetheless, our results hold for constant $p$'s, and theoretically speaking, handling the case of $p=\omega(1/\sqrt{m})$ is more challenging than that of constant $p$'s. This is because, under constant $p$'s, the sample information is more informative, easing the analysis of the algorithm.} (which includes constant values for $p$ that does not scale with $m$), the asymptotic competitive ratio of our algorithm is on the order of $1-\Theta(1/(p\sqrt{m}))$.\footnote{ For the case of $p= O(\frac{1}{\sqrt{m}})$ where $p$ goes to zero very fast as $m$ grows, in Section \ref{sec:algpsmall}, we design another algorithm whose competitive ratio is $1/2$. We note that  our upper bound in Section \ref{sec:uppersmallp} confirms the challenges of such an extreme case for the sampling probability. As we show there, when $p$ goes to zero very fast (i.e., $p=O(1/m)$), obtaining asymptotic optimality is impossible.  } As we show in Section \ref{sec:limitedsample}, the same asymptotic competitive ratio continues to hold even when there is a capacity constraint during the test period. (There we show that when the number of resources during the test period is in the order of $\omega(\sqrt{m})$, the asymptotic competitive ratio  of (a slightly modified version of) Algorithm \ref{alg:ptas} is in the order of $1-\Theta(1/(p\sqrt{m}))$.)

This result shows that the sample information can be extremely useful to improve the performance of resource allocation algorithms. This is because in the absence of the sample information, and even when the expected rewards of agents are fully known, as shown in \cite{ball2009toward}, the competitive ratio of any algorithm cannot exceed $1/(2-r_2/r_1)$ no matter how large $m$ is, where $r_1 > r_2$ are the  expected rewards for the the two types of agents. Here, we show that we can break the barrier of  $1/(2-r_2/r_1)$ by taking advantage of the sample information in a very challenging setting, where the expected reward of agents is not known to the decision-maker. This is mainly  because the sample information can be used to infer some knowledge about the number of agents of different types in the online arrival sequence,  gaining some knowledge about the adversarially-chosen arrival sequence.  

That being said, due to the adversarial nature of the demand process, it is not obvious if the sample information can lead to asymptotic optimality when $p =\omega(1/\sqrt{m})$. Even when the sampling probability  $p$ is constant, with a large number of resources, the sample information may not be very illuminating. Consider  the case where either the market size $h$ or $\ell$ is small. In such cases, the number of agents in the sample set for at least one of the types will be small and hence the decision-maker cannot correctly estimate the expected rewards of the agents. This, in turn, can lead to the decision-maker protecting the wrong (low-reward) type. Considering this, it is quite remarkable that our algorithm obtains asymptotic optimality even when $p$ shrinks as a function of $m$. As we will discuss shortly, the competitive ratio of our algorithm is tight concerning both $m$ and $p$.

We now comment on our technical contributions in characterizing  the competitive ratio of our algorithm.  
The proof of Theorem \ref{thm:2} is quite involved and is divided into three main cases, where each case bounds the competitive ratio of the algorithm when the total number of agents of type $1$ and $2$ falls into a certain region. The most challenging region is the one in which the total number of type $1$ agents (i.e., $h$) is large. In this case, the algorithm may lose reward in three aspects: not protecting type $1$ agents, over-protecting type $1$ agents, and under-protecting type $1$ agents. Recall that in our setting, the decision-maker does not even know which type has the higher reward, and hence our algorithm may wrongfully protect type $2$ agents. In addition, for large values of $h$, even when the right type is protected, by over- and under-protecting the protected type, our algorithm can lose some rewards. At a high level, we overcome these challenges, by showing that either (i) the algorithm protects the right (high-reward) type with high probability, or  (ii) thanks to setting protection levels using the sample information, the loss of the algorithm is small if it protects the wrong (low-reward) type.

We note that while characterizing the bound for the region with large $h$ is the most challenging one, the proof shows that the competitive ratio is smallest when $h$ is small. This is because, with small $h$, the sample information may not even reveal which type has the highest reward, leading the algorithm to protect the wrong type.

\textbf{Upper bound on the competitive ratio of any deterministic and randomized algorithms.} Our proposed algorithm obtains the asymptotic competitive ratio of $1-\Theta(1/(p\sqrt{m}))$ when $p =\omega(1/\sqrt{m})$. Such a superb performance makes us wonder if we can design an algorithm with even a better asymptotic competitive ratio. In Section \ref{sec:upper}, we present Theorems \ref{thm:upperbound} and \ref{thm:upperboundrandom}, which show that with $p = \omega(1/\sqrt{m})$, even when the decision-maker is fully aware of the expected reward of agents, no deterministic and randomized algorithms can obtain an asymptotic competitive ratio better than $1-\Theta(1/(p\sqrt{m}))$. To show the upper bound of  $1-\Theta(1/(p\sqrt{m}))$ on the competitive ratio for any deterministic algorithm,  in the proof of Theorem \ref{thm:upperbound},  we consider a family of arrivals, wherein this family, a large number of type $2$ agents (i.e., the type with a lower reward) arrive first, followed by some number of type $1$ agents. Under this family of arrivals, any deterministic algorithm has to decide how many type $2$ agents to accept based on the number of type $1$ agents in the sample. 
We show that due to the lack of precise knowledge about the number of type $1$ agents, no deterministic algorithm can do better than $1-\Theta(1/(p\sqrt{m}))$ on the constructed family of arrivals. To show the same result for any randomized algorithms, we first derive a variation of Yao's Lemma that could be of independent interest; see Lemma \ref{lem:yao}. We then devise a distribution over the family of the arrival sequence that we considered in Theorem \ref{thm:upperbound} and show the upper bound of $1-\Theta(1/(p\sqrt{m}))$ using  Lemma \ref{lem:yao}.

For the case of $p =O(1/\sqrt{m})$, we present another upper bound in Section \ref{sec:uppersmallp}. This non-asymptotic upper bound, which is valid for any value of $m$, shows that when $p=O(1/m)$, it is not possible to exceed the  bound of $1/(2-r_2/r_1)$, where  $1/(2-r_2/r_1)$ is the upper bound in \cite{ball2009toward} for a setting with no sample information, but known rewards. This shows that when $p$ is very small, the sample information is not sufficient to overcome the challenge of not knowing the rewards. 

\textbf{Case study.} 
We perform a case study in Section \ref{sec:simulation} using the ``Laboratory-Confirmed COVID-19-Associated Hospitalizations" dataset, which contains the number of bi-weekly cases of COVID-19-associated hospitalizations in the US from March $7^{\text{th}}$, $2020$ to February  $5^{\text{th}}$, $2022$. The dataset is obtained from the following website: \url{gis.cdc.gov/grasp/covidnet/covid19_5.html}.  We study how to 
use our algorithm to allocate limited hospital resources (e.g., a certain medicine) to different types of COVID-19 patients while having access to some sample information. We show that the average  competitive ratio of  our algorithm in various  realistic settings exceed $0.88$ and our algorithm substantially outperforms  a naive algorithm that does not use the sample information. Further, we observe that  our algorithm maintains its performance when it only has access to an estimate of the sampling probability $p$.

\subsection{Other Related Works}\label{sec:related}
Our work is related to various streams of works in the literature.

\textbf{Online decision-making with samples.} As stated earlier, our model is related to some of the recent works on online decision-making with samples.  
 \cite{correa2021secretary} study the secretary problem under both adversarial and random order  arrival models with the independent sampling process, where in the secretary problem, the goal is to maximize the probability of selecting the best applicant. For the adversarial model,  with the knowledge of the true sample information, they design a simple threshold-based algorithm (with the threshold of $k=1/(1-p)$) that achieves a competitive ratio of $kp^k(1-p)$. They show that the bound of $kp^k(1-p)$ is  tight. (Similar results are obtained for the random order arrival models.)  \cite{kaplan2022online} study an online weighted bipartite matching problem under an adversarial arrival with a  similar independent sampling process. They analyze a simple greedy algorithm which does not depend on $p$ and show that it achieves a competitive ratio of $p(1-p)$ for $p\leq 1/2$ and $1/4$ for $p>1/2$. We note that without samples,  the best 
 competitive ratio for the adversarial secretary problem and the online matching problem in \cite{kaplan2022online} is 
 zero. Thus, similar to our work,  the aforementioned works  show that the sample information can be very helpful
 in enabling more efficient online decision-making. 
 
\textbf{Revenue management with limited demand information.} In practice, it is hard to fully predict demand, and hence revenue management with limited demand information becomes an essential problem. \cite{lan2008revenue} study a single-leg revenue management problem with predicted lower and upper bounds of the demand for each type of agent and present an optimal algorithm. \cite{perakis2010robust} study a network revenue management problem with lower and upper bounds on demand of each type of agent and design an approximation algorithm.\footnote{See also  \cite{hwang2021online, esfandiari2015online} for  works that consider a hybrid arrival model with both adversarial and stochastic components. For such models, the stochastic component can reduce the uncertainty in the demand process. }
Besides the quantity-based revenue management problem discussed above, some works study the price-based revenue management problem with limited demand information. For example, \cite{besbes2012blind} design an algorithm that dynamically adjusts prices to maximize  rewards under the model where the demand at each point in time is determined by the price. \cite{araman2011revenue} study a price-based revenue management problem with a parametric demand model (e.g., linear demand model with unknown coefficients). See also \cite{bu2020online} for a work that considers the problem of learning optimal price while having access to some offline demand data under certain prices.  Our work contributes to this literature by considering a novel and practical adversarial model with samples, where the samples provide limited demand information.      

  \textbf{Online algorithm design with machine learned advice.} The sample information in our setting provides some information to the decision-maker regarding the adversarially-chosen arrival sequence, allowing us to significantly improve the worst-case guarantee of our algorithm. Improving the worst-case guarantee of online algorithms with the help of extra information has been the topic of recent literature on algorithm design with machine-learned advice. See, for example,  \cite{antoniadis2020secretary} for using advice on the maximum value of secretaries in the online secretary problems, \cite{lattanzi2020online} for using advice   
  on the weights of jobs in online scheduling problems,  \cite{lykouris2018competitive} for using the advice in the online caching problem,  \cite{balseiro2022single} for using advice in single leg revenue management problems,\footnote{In \cite{balseiro2022single}, the advice is a predicted demand vector.} and \cite{jin2022online} for using advice in online bipartite matching problems. Our work contributes to this line of work by presenting the first algorithm that optimally exploits rather unstructured advice obtained through the sample information in an online resource allocation problem.   
  
  \textbf{Multi-armed bandits.} Our setting is also related to the vast literature on multi-armed bandits; see, for example, \cite{Thompson1933, Auer2002, Auer2002nonstochastic, balseiro2019contextual,  Bart, BW, chen2022bridging, chen2021dynamic}. In this literature, it is assumed that there are some arms/options with unknown expected rewards, and the goal is to identify the best arm by suffering from minimal regret, where the regret is computed against a clairvoyant policy that knows the expected  reward of arms in advance. 
  In our setting, similar to the bandit settings,  the rewards of agents are unknown. However, unlike the bandit setting,  the algorithm aims to partially learn these rewards via the sample information, rather than the feedback it receives throughout the allocation period. Our upper bound results, which are obtained under the setting when the rewards are fully known,  show that using the feedback during the allocation period does not improve the asymptotic competitive ratio of our algorithm. This is because, in our setting, the order of arrivals and the number of agents of different types are chosen adversarially. Hence, in the worst case,  the feedback throughout the allocation periods does not add any value, as the decision-maker is not able to acquire the right feedback at the right time.  What feedback the decision-maker can receive is mainly determined by the adversary. 
  
\textbf{Online resource allocation.} Our work is also related to the rich literature on online resource allocation.  \cite{devanur2009adwords, feldman2010online}, \cite{agrawal2014dynamic}, and \cite{chen2021approximation} study this problem in a stochastic setting with i.i.d. demand arrivals. In these works, the authors use the primal-dual technique to design algorithms with sub-linear regret, where the algorithms aim to learn the optimal dual variables associated with resource constraints.   The primal-dual technique is also effective for adversarial demand arrivals even though attaining sub-linear regret in adversarial settings is generally impossible
\citep{mehta2007adwords,buchbinder2007online, golrezaei2014real}. See also \cite{balseiro2020best} for a recent work that shows how the primal-dual technique results in well-performing algorithms for various demand processes. While in our work, we do not rely on the primal-dual technique,  our work contributes to this literature by presenting a model that---with the help of the sample information---bridges the gap between stochastic and adversarial arrivals, allowing us to bypass the impossibility results in the adversarial settings which is obtained by \cite{ball2009toward}.

\section{Model}\label{sec:model}
We consider a decision-maker who would like to allocate (identical) units of a resource/service to two types of unit-demand agents. ({In Section \ref{subsec:nagent}, we discuss how to extend our algorithm (Algorithm \ref{alg:ptas}) to a setting with more than two types.}) The expected reward of agents of type  $i\in \{1, 2\}$ upon receiving one unit of the resource  is $r_i \in (0,1)$, where without loss of generality,  we assume that $r_1> r_2$ and define $\alpha =\frac{r_2}{r_1}\in (0,1)$. The decision-maker, however, is not aware of the expected reward of agents; she does not even know which type of agents attains a higher expected reward. As stated in the introduction, this setting captures scenarios where we would like to allocate new drugs, services, and products to customers. Given the lack of knowledge about the expected rewards, before the allocation period starts, the decision-maker aims to collect some information about the unknown expected rewards of agents during a \emph{test period}.

\subsection{Test Period and Sample Information}
 During the test period, the decision-maker aims to outreach the market by, for example, making a public announcement.  
Let $h,\ell\ge 0$ be respectively the market size of agents of types $1$ and $2$, i.e., the total number of agents of types $1$ and $2$ that are interested in the resource.  The market size $h$ and  $\ell$, which are unknown to the decision-maker,   can take any arbitrary values; that is, they are chosen adversarially.  We assume that during the test period, with probability $p\in (0,1)$, each of $h+\ell$ agents sees and responds to the outreach program run by the decision-maker, where the \emph{sampling probability}  $p$ is known to the decision-maker. ({In Section \ref{subsec:heterop},   we show that our algorithm and results can be simply extended to a setting where each agent of type $i\in \{1,2\}$ reacts to the outreach program with probability $p_i$, where $p_1\ne p_2$.})

Note that  we do not require  $p$ to be a constant.  Our results  not only hold for any constant value of the {sampling probability} $p$, but also hold when  $p =\omega(\frac{1}{\sqrt{m}})$, where $m$ is the total number of resources when the test period ends.  Furthermore, 
 the  assumption of knowing the sampling probability  is motivated by the fact that the outreach program is designed by the decision-maker herself and hence an accurate estimate of the sampling probability $p$ is available to her. Nonetheless, in Section \ref{sec:negative}, we show that without knowing $p$, achieving asymptotic optimality is not possible.   See also   our case study in Section \ref{sec:simulation} and our numerical studies in Appendix \ref{subsec:partialknown} where  we investigate the robustness of our proposed algorithm to the lack of exact knowledge about $p$.

Any agent who responds to the outreach program gets one unit of the resource and reveals his realized reward to the decision-maker. \footnote{The assumption that any agent in the test period gets a resource can be relaxed as we show in Section \ref{sec:limitedsample}. More specifically, there we show that our asymptotic result still holds when the number of resources in the test period is in the order of $\omega(\sqrt{m})$, which includes the  case of $\Theta(m)$.}  Here, for simplicity, we assume that the realized reward of type $i\in \{1,2\}$ agents is either $0$ or $1$. This assumption can be  relaxed by letting the realized rewards  take more than two values. Let $s_i$, $i\in \{1,2\}$, be the (random) number of agents of type $i$ that are reached/\emph{sampled} during the test period. Note that $s_1$ is drawn from a binomial distribution with parameters $h$ and $p$ (i.e., $s_1\sim \text{Bin}(h, p)$) and $s_2\sim \text{Bin}(\ell, p)$. The set of realized rewards of $s_i$ sampled agents is denoted by $ \boldsymbol \rho_i \in \{0,1\}^{s_i}$, where $ \boldsymbol \rho_i= \{\rho_{i,1}, \rho_{i,2}, \ldots, \rho_{i, s_i}\}$ and $\rho_{i,j}$, $j\in [s_i]$, is drawn from a Bernoulli distribution with the success probability of $r_i$.    Throughout the paper, we denote $(s_1, s_2)$ and $(\boldsymbol \rho_1, \boldsymbol \rho_2)$ with $\mathbf s$ and  $\boldsymbol \rho$, respectively.  In addition, we refer to $\boldsymbol \psi :=(\mathbf s, \boldsymbol \rho)$ as the \emph{sample information} that decision-maker obtains during the test period. (Note that $\boldsymbol \rho$ depends on the random variables $\mathbf s$. But, to simplify the exposition, we do not show the dependency of $\boldsymbol \rho$ to  $\mathbf s$.)

\subsection{Online Allocation Period}
Having described the test period, we are now ready to explain the allocation period. Let $m\ge 2$  be the number of available resources  at the beginning of the allocation period. During this period, the rest of the market, i.e., $h-s_1$ type $1$ agents  and $\ell-s_2$ type $2$ agents arrive one by one over time in an arbitrary order. We denote the number of agents of type $i$ during the allocation period with $n_i$ and we further denote $(n_1, n_2)$ with $\mathbf n$.  Let $z_t\in \{1, 2\}$, $t=1, 2, \ldots, $ be the type of the agent in time period $t$ within the allocation period. The agent type $z_t$, which can represent any available information at the time of the decision,  is observable to the decision-maker at  time period $t$, but given the online nature of the problem during the allocation period,  $z_{t'}$ for any $t'> t$ is not observable at time period $t$. 
We define $I= (z_t)_{t\ge 1}$ as the \emph{online arrival sequence} and note that  $n_i =|\{t \le |I|: z_t= i\}|$, $i\in \{1, 2\}$. 

Upon the arrival of the agent of type $z_t$ in time period $t$, the decision-maker has to make an irrevocable acceptance/rejection decision regarding that agent. If the decision-maker accepts the agent, she allocates the agent one unit of the resource.
Otherwise, no resource will be allocated to the agent, and that agent will not come back. If the agent gets accepted, he  reveals his realized reward to the decision-maker.

In this work, to ease the analysis, we will focus on the continuous version of our problem, where the resources are divisible and  we are allowed to partially accept  agents.    We can still capture all the core ideas regarding sample information in continuous models. For the discrete model, our  algorithm can still be applied and are asymptotically optimal as $m \to \infty$. To deal with small values of $m$, one can apply a technique used in \cite{ball2009toward} to randomly accept or reject agents, obtaining the same performance guarantees in expectation.

\subsection{Performance Measure}
The goal of the decision-maker during the allocation period is to yield high total rewards while being uncertain about the number and the order of agents, as well as, their expected rewards. For an algorithm $A$, online arrival sequence  $I$, and the sample information,  $\boldsymbol \psi =(\mathbf s, \boldsymbol \rho)$,   let $\rew_{A}\left(I, \boldsymbol \psi \right)$ be the cumulative expected reward of algorithm $A$ across all the time periods.  
We measure the performance of an algorithm $A$ using the following competitive ratio (CR) definition, which compares our algorithm to the  optimal clairvoyant benchmark that knows the arrival sequence $I$ and the expected reward of agents ($r_i$'s) in advance.\footnote{Observe that the  optimal clairvoyant benchmark does not know the realized rewards in advance. This optimal clairvoyant benchmark is consistent with that used in \cite{ball2009toward}. Nevertheless, in Section \ref{sec:realcr}, we discuss a stronger benchmark that knows the realized rewards. We show that under a strong clairvoyant benchmark that knows the realized reward of agents, it is not possible to design an algorithm whose competitive ratio w.r.t. this strong  benchmark goes to one as $m$ goes to infinity. }
Let $\textsc{opt}(I)$ be the optimal clairvoyant  cumulative expected reward that can be obtained from $I$ using $m$ units of the resource.  The CR of an algorithm $A$ is then defined as 
\begin{equation} \label{eq:CR_early}
    \text{CR}_{A} =  \inf_{(h, \ell)}~\E_{\boldsymbol \psi}
    \left[\inf_{I}\frac{\E[\rew_{A}\left(I, \boldsymbol \psi\right)]}{\textsc{opt}(I)}\right]\,.
\end{equation}
Here, the inner expectation is with respect to  (w.r.t.) any randomness in algorithm $A$, and the outer expectation is w.r.t. the arrival sequence $I$ and the sample information $\boldsymbol \psi$. Note that given $(h, \ell)$, the arrival sequence $I$ and more precisely, the number of agents of type $i\in \{1, 2\}$ in $I$ (denoted by $n_i$), is random. Further, observe that in our definition of CR, we take infimum over (i) the size of the market, i.e., $h, \ell$, and (ii) the order of arrivals in $I$.

\section{Impossibility Results}\label{sec:negative}
In this section, to
shed light on the challenges in our setting and the necessity  of some of the assumptions we made, 
we present alternative scenarios for which  we show that it is impossible to obtain a CR that goes to one as the number of resources $m$ goes to infinity. 
These scenarios are listed below.

\begin{itemize}[leftmargin=*]
\item \textit{Scenario 1: No Sample Information.} Consider the setting in Section \ref{sec:model} where the rewards $r_i$ are unknown and the arrival sequence during the allocation period is chosen adversarially. But, unlike our original  setting in  Section \ref{sec:model}, the sampling probability $p=0$, and hence, no sample information is available. For this scenario, we show that  the CR is upper bounded by $1/(2-\alpha)$, where $\alpha = r_2/r_1$ and  $1/(2-\alpha)$ is the bound in \cite{ball2009toward}.
This impossibility result highlights the importance of the sample information in achieving the asymptotically optimal CR. See also Section \ref{sec:uppersmallp} where we show that even with $p= O(1/m)$,  the CR  cannot exceed $1/(2-\alpha)$.

\item \textit{Scenario 2: Unknown Sampling Probability.} Consider the setting in Section \ref{sec:model} with sample information, unknown rewards, and adversarial arrival sequences. But, assume that the sampling probability $p$ is completely unknown to the decision-maker.  For this scenario, we  show an upper bound of  $1/(2-\alpha)$ on the CR. This impossibility result emphasizes the importance of knowing   the sampling probability.  Note that as stated earlier, similar to our setting,  \cite{correa2021secretary} also assume that the sampling probability $p$ is known and they crucially  use the knowledge of $p$ to design optimal algorithms for the secretary problem. 

\end{itemize}

 \begin{proposition} [Impossibility  Results] \label{prop:negativeresults}
For the scenarios $1$ and $2$, defined above, the CR of any deterministic and randomized algorithm  is upper bounded by $1/(2-\alpha)$,  where $\alpha= r_2/r_1$. 
\end{proposition}

The proof of Proposition \ref{prop:negativeresults} is presented in Section \ref{sec:proof:negative} in the appendix.

\section{Online Resource Allocation Algorithm with Samples}\label{sec:alg}
In light of our impossibility results, here, under our original model presented in Section \ref{sec:model},  we present an online resource allocation algorithm that  is asymptotically optimal as the number of resources  $m$ goes to infinity; that is, the   CR of our algorithm converges to one as $m$ goes to infinity. 

In our setting, 
 the decision-maker receives some information about the number of agents of each type in the online arrival sequence.  The decision-maker further obtains some partial information about the expected reward of each type.   Recall that the decision-maker observes 
 $s_i$ samples, $i\in \{1, 2\}$, for each type $i$ agents where $s_1 \sim \text{Bin}(h,p)$ and  $s_2 \sim \text{Bin}(\ell,p)$. In addition, the decision-maker  observes the realized reward of $s_i$ agents of type $i\in \{1, 2\}$, denoted by $\boldsymbol \rho_i =(\rho_{i,1}, \rho_{i,2}, \ldots, \rho_{i, s_i})$.

Our proposed algorithm (Algorithm \ref{alg:ptas}) takes advantage of both $\boldsymbol s=(s_1, s_2)$ and $\boldsymbol{\rho} =(\boldsymbol{\rho}_1, \boldsymbol{\rho}_2)$. The algorithm uses these pieces of information to obtain an estimate of the expected reward of each type $i$, denoted by $\widehat r_i$. In particular, if $s_i> 0$, the algorithm simply uses the sample average of the realized reward observed in  $\boldsymbol{\rho}_i$ as an estimate for $r_i$; that is, 
\[\widehat r_i = \frac{\sum_{j=1}^{s_i} \rho_{i,j}}{s_i}\qquad i\in \{1, 2\}\,.\]
When $s_i =0$, $\widehat r_i$ is randomly drawn from a uniform distribution in $[0,1]$. 

Having access to these estimates, the algorithm follows  a protection level policy to protect the type that has the highest estimated expected reward. The definition of a protection level policy for type $i\in \{1, 2\}$ and protection level $x$ is stated below in Definition \ref{def:protect_policy}.

\begin{definition}[Protection Level Policy for Type $i$ and Protection Level $x$] \label{def:protect_policy}
 In this algorithm, agents of type $i$ are always accepted unless there is no resource left. Agents of type $-i$ will be accepted if the number of accepted type $-i$ agents is less than $m-x$ and there is  resource left.\footnote{In this algorithm, the last accepted type $j\in\{1, 2\}$ agent may get accepted partially.}
\end{definition}

When $\widehat r_1 > \widehat r_2$, then the algorithm assigns a protection level of  $x_1= \min\{m,s_1\frac{1-p}{p} \}$ for type $1$ agents, where we note that conditioned on $s_1$, 
$s_1\frac{1-p}{p}$ is equal to the expected number of type $1$ agents in the online arrival sequence, i.e., $n_1$. Recall that $n_1 = h- s_1$ and $s_1 \sim \text{Bin}(h, p)$. 
 On the other hand, when  $\widehat r_2 \ge \widehat r_1$,  the algorithm assigns a protection level of  $x_2= \min\{m,s_2\frac{1-p}{p} \}$ for type $2$ agents. The description of the algorithm can be found in Algorithm \ref{alg:ptas}.

\begin{algorithm}[htbp] 
\caption{Online Resource Allocation Algorithm with Samples}
\footnotesize
\begin{itemize}
\item[] \textbf{Input:} The number of resources $m$ and sample information $\boldsymbol \psi =(\mathbf s, \boldsymbol \rho)$, where $\mathbf s =(s_1, s_2)$, $\boldsymbol \rho = (\boldsymbol \rho_1, \boldsymbol \rho_2)$, and for any $i\in \{1, 2\}$,   $\boldsymbol \rho_i =\{\rho_{i,1}, \rho_{i,2}, \ldots, \rho_{i, s_i}\}$.

\end{itemize}
\begin{enumerate}
  \item If the number of samples for type $i\in \{1,2 \}$, i.e., $s_i$, is positive, define \[\widehat r_i = \frac{\sum_{j=1}^{s_i} \rho_{i,j}}{s_i}\] as an  estimate of the expected reward of type $i\in \{1, 2\}$. Otherwise,  $\widehat r_i\sim \text{uniform}(0,1)$, where $\text{uniform}(0,1)$ is the uniform distribution in $[0,1]$. 
  \item If $\widehat r_1 > \widehat r_2$, set the protection level $x_1 = \min\left\{m,s_1 \frac{1-p}{p}\right\}$ for type $1$, and run the protection level policy for type $1$ agents and protection level $x_1$ as per Definition \ref{def:protect_policy}.
  
\item If $\widehat r_1 \le \widehat r_2$,  set the protection level $x_2 = \min\left\{m,s_2 \frac{1-p}{p}\right\}$ for type $2$ and run the protection level policy for type $2$ agents and protection level $x_2$ as per Definition \ref{def:protect_policy}.
\end{enumerate}

\label{alg:ptas}
\end{algorithm}

We highlight that in the described algorithm, the adversary's  choice of $h$ and $\ell$  influences the algorithm's estimate for the expected rewards, as well as, the protection levels $x_1$ and $x_2$. Furthermore, the adversary's  choices  impact the algorithm's decision about protecting a certain type. Note that when the estimated expected rewards $\widehat r_i$, $i\in\{1, 2\}$, are noisy, the algorithm may end up protecting a wrong type (i.e., the type with a lower reward). Despite these challenges, as we show in  Theorem \ref{thm:2} and Proposition \ref{thm:3}, Algorithm \ref{alg:ptas} manages to perform very well. Theorem \ref{thm:2} presents a lower bound on the CR of Algorithm \ref{alg:ptas} for any value of $m$ and $p$, where we recall that $m$ is the number of resources and $p$ is the sampling probability.  Further,  Proposition \ref{thm:3} shows how this lower bound scales with $m$ and $p$ as $m$ goes to infinity. In particular,  Proposition \ref{thm:3} shows that 
the asymptotic CR of Algorithm \ref{alg:ptas} is   
$1-\Theta(1/(p\sqrt{m}))$.

Later in Section \ref{sec:upper}, we  show this bound is tight in the sense that  no randomized or deterministic algorithm can break the asymptotic bound of $1-\Theta(1/(p\sqrt{m}))$. The asymptotic optimality of Algorithm \ref{alg:ptas} is quite remarkable given the simplicity of the algorithm and our strong adversary. Note that the adversary through choosing the market size $h$ and $\ell$ can control to what extent  the sample information is informative.   
When $h$ and $\ell$ are both large,  the number of agents in the sample is also  large and hence the sample information can be very useful in reducing uncertainty in demand and expected rewards of agents. However, when either $h$ or $\ell$ is small---even when the number of resources $m$ goes to infinity and the sampling probability is constant---the sample information may not be  very informative, as it does not even allow the decision-maker to obtain an accurate estimate of rewards. As a result, our algorithm may end up protecting the wrong type. Thus, the fact  that our simple algorithm  obtains an asymptotically tight CR of $1-\Theta(1/(p\sqrt{m}))$ is outstanding.

In addition, Algorithm \ref{alg:ptas} does not use the realization of the reward of arriving agents  during the allocation period. One may wonder if using such information would help us design a better algorithm. The answer is no. Our upper bound result in Section \ref{sec:upper} shows that  adaptive algorithms that further use the realized rewards of agents during the allocation period  cannot obtain a better   asymptotic CR.  Intuitively speaking, this is because the arrival sequence is chosen by an adversary who can  manipulate   how fast the rewards are learned by  choosing the order of arrivals. That is, the adversary can prevent the decision-maker to have a good estimate about the rewards  when having such estimates is crucial. Therefore, using the realized rewards during the allocation period does not give algorithms an edge.

Nonetheless, 
one can consider an adaptive version of Algorithm \ref{alg:ptas}. In this algorithm,  estimates of rewards  $\widehat r_i$ are updated once an agent is accepted: we can again let $\widehat r_i$ be the empirical mean of realized rewards of all agents of type $i$ who receives the resource.  At any time $t$, if we observe $\widehat r_i > \widehat r_j$, we give type $i$ agent protection level $x_i=\min\{m-Z_{j,t}, s_i\frac{1-p}{p}-Y_{i,t}  \}$, where $Y_{i,t}$ is the number of type $i$ agents arrived before time $t$ and $Z_{j,t}$ is the number of type $j$ agents accepted before time $t$. In Section \ref{sec:simulation}, via  our case study, we evaluate this adaptive algorithm, and we show that even in realistic arrival scenarios, the CR of this adaptive algorithm is almost the same as that of Algorithm \ref{alg:ptas}, verifying our theoretical result.

\subsection{Competitive Ratio of Algorithm \ref{alg:ptas}}\label{sec:CR_proof}

Here, we present a lower bound on the CR of Algorithm \ref{alg:ptas}. To do so, we start with the following lemma. This lemma  shows that to characterize   a lower bound on the CR of the algorithm,   it suffices to only consider ordered arrival sequences in which any type $2$ agent arrives before any type $1$ agent. The proofs of all the lemmas in this section are presented in Section \ref{app:proof_thm_2} in the appendix.

\begin{lemma} [Worst Order] \label{lm:order}
For any realization $n_1$, $n_2$, where $n_i$ is the number of type $i$ agents in the arrival sequence, let $I_{\text{order}}$ be an \emph{ordered} arrival sequence under which $n_2$ type $2$ agents arrived first, followed by $n_1$ type $1$ agents. Let $I$ be any arrival sequence that contains $n_2$ type $2$ agents and $n_1$ type $1$ agents. Then,  under  Algorithm \ref{alg:ptas} (denoted by $A$), we have
\begin{equation}
\frac{\E[\rew_{A}\left(I_{\text{order}}, \boldsymbol \psi\right)]}{\textsc{opt}(I_{\text{order}})}\le  \frac{\E[\rew_{A}\left(I, \boldsymbol \psi\right)]}{\textsc{opt}(I)}
    \,.
\end{equation}
\end{lemma}

With a slight abuse of notation, let $\textsc{opt}(\mathbf n)$ be  the  optimal clairvoyant cumulative expected reward under an online  ordered arrival sequence $I$ with $n_i\in \{1,2\}$ type $i$ agents. Similarly, $\rew_A\left(\mathbf n, \boldsymbol \rho\right)$ is the (expected) reward of Algorithm  
\ref{alg:ptas} under $\boldsymbol \rho$ and  an online  ordered arrival sequence $I$ with $n_i$ type $i$ agents.
Then, 
 Lemma \ref{lm:order}  allows us to rewrite the CR of Algorithm \ref{alg:ptas}  as follows 
\begin{align}
    \text{CR}_{A} &= \inf_{(h, \ell)}~\E_{\boldsymbol \psi}
    \left[\inf_{I}\frac{\E[\rew_{A}\left(I, \boldsymbol \psi\right)]}{\textsc{opt}(I)}\right]\, = \inf_{(h, \ell)}~\E_{\boldsymbol \psi}
    \left[\frac{\rewGen}{\textsc{opt}(\mathbf n)}\right] \notag\\
     \label{eq:crdefinitionA}
     &=   \inf_{(h, \ell)}~\E_{\boldsymbol  {n}}
    \left[\frac{\E[\rewGen]}{\textsc{opt}(\mathbf n)}\right]
    \,. 
\end{align}
 Note that in the inner expectation of the last expression, we take expectation w.r.t. any randomness in the algorithm and  $\boldsymbol \rho$   given a realization of $\mathbf n$.

To bound $\text{CR}_A$, we define the following  \emph{good event}, denoted by,  $\mathcal G(\boldsymbol \psi; (h, \ell))$,  as the event under which Algorithm \ref{alg:ptas} correctly assumes that type $1$ agents have a higher expected reward than type $2$ agents when the sample information is $\boldsymbol{\psi}$ and the market size for type $1$ and $2$ agents is respectively $h$ and $\ell$:
\begin{align}\mathcal G(\boldsymbol \psi; (h, \ell)) =\{\widehat r_1 > \widehat r_2\} =  \left\{\frac{\sum_{j=1}^{s_1} \rho_{1,j}}{s_1}> \frac{\sum_{j=1}^{s_2} \rho_{2,j}}{s_2}  \right\}\,. \label{eq:good_event}\end{align} 
Then,  given a realization of $\mathbf s$ or equivalently $\mathbf n $, we have 
\begin{align}\notag
    \E[\rewGen] &= \E[\rew_A\left(\mathbf n, \boldsymbol \rho\right)\cdot\boldsymbol{1}({\mathcal G(\boldsymbol \psi;(h,\ell))})]+\E[\rew_A\left(\mathbf n, \boldsymbol \rho\right)\cdot \boldsymbol{1}({\mathcal G^{\mathsf{C}}(\boldsymbol \psi;(h,\ell))})]\\ \notag
    &=\E[\rew_A\left(\mathbf n, \boldsymbol \rho\right) \big\vert {\mathcal G(\boldsymbol \psi; (h, \ell))}] \cdot \Pr(\mathcal G(\boldsymbol \psi; (h, \ell)))\\
    &+ \E[\rew_A\left(\mathbf n, \boldsymbol \rho\right) \big\vert {\mathcal G^{\mathsf{C}}(\boldsymbol \psi; (h, \ell))}] \cdot \Pr(\mathcal G^{\mathsf{C}}(\boldsymbol \psi; (h, \ell)))\,, \label{eq:rew_good_event}
\end{align}
where $\boldsymbol{1}(\mathcal A)$ is $1$ if an event $\mathcal A$ occurs and zero otherwise. Here, given $\boldsymbol{n}=(n_1,n_2)$, the expectation is taken on $\rho_{1,j} \sim \text{Ber}(r_1)$ for $j \in [s_1]$ and $\rho_{2,j} \sim \text{Ber}(r_2)$ for $j \in [s_2]$, where $\text{Ber}(r)$ is a Bernoulli distribution with a success probability of $r\in [0,1]$. Further, $\mathcal G(\boldsymbol \psi;(h,\ell))$ is the good event defined in Equation \eqref{eq:good_event}.
In Equation \eqref{eq:rew_good_event}, the  expression  $\E[\rew_A\left(\mathbf n, \boldsymbol \rho\right) \big\vert {\mathcal G(\boldsymbol \psi; (h, \ell))}]$ presents the reward of Algorithm \ref{alg:ptas} under a specific realization of   $\mathbf n$ when the good event $\mathcal G(\boldsymbol \psi; (h, \ell))$ happens. Note that for a given realization of  $\mathbf n$ and under  the good event, the algorithm's action and hence its reward does not depend on $\boldsymbol \rho$. This is because when the good event happens, the algorithm assigns a protection level $x_1$ to type $1$, where $x_1$ only depends on $s_1$, not $\boldsymbol \rho$.  Finally, in $\Pr(\mathcal G(\boldsymbol \psi; (h, \ell)))$, we take expectation with respect to $\boldsymbol{\rho}$ and any randomness in the algorithm in defining $\widehat r_i$, $i\in [2]$, given a realization of $\mathbf n$. 

To bound the CR of Algorithm \ref{alg:ptas}, we consider the following three cases based on the number of type $1$ and $2$ agents chosen by the adversary (i.e., $h$ and $\ell$):
\begin{enumerate}[leftmargin =*]
    \item []\textit{Case 1.} In this case, the number of type $1$  agents is small (less than $\mzero $). That is, $(h,\ell) \in \mathcal R_1$, where $\mathcal R_1 =\{(h,\ell):  h<\mzero , \ell \ge 0\}$.

\item []\textit{Case 2.} 
In this case, while the number of type $1$ agents is large, the number of type $2$ agents is small.  That is, $(h,\ell) \in \mathcal R_2$, where $\mathcal R_2 =\{(h,\ell):  h\ge \mzero ,  \ell< \mzero \}$.
    \item []\textit{Case 3.} In this case, the numbers of both type $1$ and $2$ agents are large.
    That is, $(h,\ell) \in \mathcal R_3$, where $\mathcal R_3 =\{(h,\ell):  h\ge \mzero ,  \ell\ge  \mzero \}$. 
\end{enumerate}\medskip

Observe that the worst case CR of Algorithm \ref{alg:ptas} can be written as 
\begin{align}\min_{i\in [3]}\left\{\inf_{(h, \ell)\in \mathcal R_i}~\E_{\boldsymbol  {n}}
    \left[\frac{\E[\rewGen]}{\textsc{opt}(\mathbf n)}\right]\right\}\,,\end{align}
where by Equation \eqref{eq:rew_good_event}, for any $i\in [3]$, we have
\begin{align}\nonumber \inf_{(h, \ell)\in \mathcal R_i}~\E_{\boldsymbol  {n}}
    \left[\frac{\E[\rewGen]}{\textsc{opt}(\mathbf n)}\right] &=\inf_{(h, \ell)\in \mathcal R_i} \Bigg\{ \E_{\mathbf n}\left[\frac{\E[\rew_A\left(\mathbf n, \boldsymbol \rho\right) \big\vert {\mathcal G(\boldsymbol \psi; (h, \ell))}] \cdot \Pr(\mathcal G(\boldsymbol \psi; (h, \ell)))}{\textsc{opt}(\mathbf n)}\right]\\ \nonumber
    &+ \E_{\mathbf n}\left[\frac{\E[\rew_A\left(\mathbf n, \boldsymbol \rho\right) \big\vert {\mathcal G^{\mathsf{C}}(\boldsymbol \psi; (h, \ell))}] \cdot \Pr(\mathcal G^{\mathsf{C}}(\boldsymbol \psi; (h, \ell)))}{\textsc{opt}(\mathbf n)}\right]\Bigg\}\,.
    \end{align}
In the rest of the paper, we use the following shorthand notation to simplify the exposition: 
\begin{align*}&\E[\rew_A\left(\mathbf n, \boldsymbol \rho \right) \big\vert {\mathcal G(\boldsymbol \psi; (h, \ell))}] = \rewGood\\
&\E[\rew_A\left(\mathbf n, \boldsymbol \rho\right) \big\vert {\mathcal G^{\mathsf C}(\boldsymbol \psi; (h, \ell))}] = \rewBad\,.\end{align*}
With this notation, we have
\begin{align}\nonumber \inf_{(h, \ell)\in \mathcal R_i}~\E_{\boldsymbol  {n}}
    \left[\frac{\E[\rewGen]}{\textsc{opt}(\mathbf n)}\right] &=\inf_{(h, \ell)\in \mathcal R_i} \Bigg\{ \E_{\mathbf n}\left[\frac{\rewGood \cdot \Pr(\mathcal G(\boldsymbol \psi; (h, \ell)))}{\textsc{opt}(\mathbf n)}\right]\\ \label{eq:CR} 
    &+ \E_{\mathbf n}\left[\frac{\rewBad \cdot \Pr(\mathcal G^{\mathsf{C}}(\boldsymbol \psi; (h, \ell)))}{\textsc{opt}(\mathbf n)}\right]\Bigg\}\,.
    \end{align}

Having defined this shorthand notation, in the following, we present   three main lemmas: Lemmas \ref{lem:c1}, \ref{lem:c3}, and \ref{thm:c4}, where  each of these lemmas provides a lower bound on 
    $\inf_{(h, \ell)\in \mathcal R_i}~\E_{\mathbf{n}}
    \left[\frac{\E[\rew_A\left(\mathbf n, \boldsymbol \rho\right)]}{\textsc{opt}(\mathbf n)}\right]$ for one of the regions $\mathcal R_i$, $i\in [3]$.

\begin{lemma}[Region $\mathcal R_1$] \label{lem:c1} 

Let $\mathcal R_1 =\{(h,\ell):  h<\mzero , \ell \ge 0\}$.  Then, we have
\[
\inf_{(h, \ell)\in \mathcal R_1}~\E_{\mathbf{n}}
    \left[\frac{\E[\rew_A\left(\mathbf n, \boldsymbol \rho\right)]}{\textsc{opt}(\mathbf n)}\right]\ge \CROne \,,
\]
where 
\[\CROne   =\CROneExp\,,\]
and  $y^{+}=\max\{y,0\}$.
\end{lemma}
\medskip

Lemma \ref{lem:c1} shows when $(h, \ell)$ falls into region $\mathcal R_1$,   the CR of Algorithm \ref{alg:ptas} is lower bounded by $\CROne$. From the expression of $\CROne$, the asymptotic bound in region $\mathcal R_1$ is $1-\Theta(1/(p\sqrt{m}))$, and based on the result we will present in Section \ref{sec:upper}, this bound is tight with respect to both $p$ and $m$. The bound of $1-\Theta(1/(p\sqrt{m}))$ is achieved when the total number of arrivals  is less than $m$ (i.e., $n_1+n_2\le m$ ). 
When the total number of arrivals is greater than $m$, a better asymptotic bound of $1-\Theta(1/\sqrt{m})$---which does not depend on $p$---can be obtained.

To show Lemma \ref{lem:c1}, by Equation \eqref{eq:CR}, we need to bound $\E_{\boldsymbol n}
    \left[\frac{\rewGood \cdot \Pr(\mathcal G(\boldsymbol \psi; (h, \ell)))}{\textsc{opt}(\mathbf n)}  \right]$ and $\E_{\boldsymbol n}
    \left[\frac{\rewBad \cdot \Pr(\mathcal G^{\mathsf{C}}(\boldsymbol \psi; (h, \ell)))}{\textsc{opt}(\mathbf n)}  \right]$ respectively. As in this region, $h$ can be very small, we cannot use any concentration inequality to bound the expectation over $\boldsymbol{n}$. Instead, we  bound $\frac{\rewGood}{\textsc{opt}(\mathbf n)}$ and $\frac{\rewBad}{\textsc{opt}(\mathbf n)}$ for any realization of $\mathbf n$ by considering different ranges for the number of agents of type $i\in [2]$ in the online arrival sequence.

\begin{lemma} [Region $\mathcal R_2$]\label{lem:c3}
 Let $\mathcal R_2 =\{(h,\ell):  h\ge \mzero , \ell < \mzero  \}$.  Then, we have
\[
\inf_{(h, \ell)\in \mathcal R_2}~\E_{\mathbf{n}}
    \left[\frac{\E[\rew_A\left(\mathbf n, \boldsymbol \rho\right)]}{\textsc{opt}(\mathbf n)}\right]\ge \CRTwo\,,
\]
where 
\[\CRTwo   =
    \CRTwoExp\,.\]
Here,  $\mOne =\mOneExp$, $\hzero =\hzeroExp$,  and $\beta=0.4215 \cdot \frac{p^2+(1-p)^2}{p(1-p)}$.
\end{lemma}

Lemma \ref{lem:c3} shows when $(h, \ell)$ falls into region $\mathcal R_2$,   the CR of Algorithm \ref{alg:ptas} is lower bounded by $\CRTwo$.
In the expression of $\CRTwo$,  $\hzero={\frac{m-\mzero }{1-p}+\frac{1}{2p^2}+\sqrt{\frac{1}{4p^4}+\frac{(m-\mzero )}{p^2(1-p)}}}$, and as $p=\omega(1/\sqrt{m})$, we have $\hzero=\frac{m}{1-p}-\Theta(\frac{\mzero}{1-p})$, 
$h_1=m-\Theta(\mzero)$, and $\frac{h_1}{m}=1-\Theta(1/\sqrt{m})$. This implies that  the asymptotic bound in region $\mathcal R_2$ is $1-\Theta(1/\sqrt{m})$, which is tight with respect to $m$, but is not related to $p$. 

To show Lemma \ref{lem:c3}, we consider three factors which contribute to the loss of Algorithm \ref{alg:ptas}, where for a realization of $\mathbf n$ and $\boldsymbol \rho$, the (normalized) loss of Algorithm $A$ is $\frac{\textsc{opt}(\mathbf n)- \rewGen}{\textsc{opt}(\mathbf n)} = 1-\frac{\rewGen}{\textsc{opt}(\mathbf n)}$.  These factors are: (i) protecting the wrong (low-reward) type, (ii) over-protecting the right (high-reward) type, and (iii) under-protecting the right type. 
In case (i), as the number of type $2$ agents is very small ($\ell<\sqrt{m}$), the loss caused by wrongfully protecting type $2$ agents is bounded by $1/\sqrt{m}$, and hence,  the CR is lower bounded by $1-1/\sqrt{m}$. This bound is obtained  by considering the worst case among any realization of $n_1$ and $n_2$.  

To bound the expected loss caused by over- or under-protecting type $1$ agents for cases (ii) and (iii), we use concentration inequalities, as the number of type $1$ agents is large under region $\mathcal R_2$. In doing so, we face one main challenge: In computing the expected (normalized) loss (or equivalently the expected CR),  it is hard to compute the expectation w.r.t. $\mathbf n$ as both the reward/loss of the algorithm (the numerator in the CR, i.e.,  $\rewGen$)  and the optimal clairvoyant solution (the denominator in the CR, i.e.,  $\textsc{opt}(\mathbf n)$) are random variables. Recall that after choosing $h$ and $\ell$, the number of agents $n_i$'s in the arrival sequence are random.

To overcome this challenge, based on the value of $h$, we split the region $\mathcal R_2$ into two sub-regions  $\underline{\mathcal{R}_2}$ and $\overline{\mathcal{R}_2}$, where  $\underline{\mathcal{R}_2} =\{(h, \ell): (h, \ell)\in \mathcal R_2, h< \hzero  \}$ and $\overline{\mathcal{R}_2} =\{(h, \ell): (h, \ell)\in \mathcal R_2, h\ge \hzero  \}$. In each of the sub-regions, we can deal with the aforementioned  expectation challenge  by bounding  the  optimal clairvoyant solution in the denominator. The splitting  allows us to (i) provide a tight bound on $\textsc{opt}(\mathbf n)$ in each of the sub-regions. This then further enables us to only focus on computing the expectation for the numerator.  We show that for the sub-region $\underline{\mathcal{R}_2}$, the CR can be lower bounded by $(1-\frac{1}{\sqrt{m}})(1-\frac{1}{m})$, and for the sub-region $\overline{\mathcal{R}_2}$, the CR can be bounded by $\min\{1-\frac{1}{\sqrt{m}}, \frac{h_1}{m} \} $. As the former term (i.e., $(1-\frac{1}{\sqrt{m}})(1-\frac{1}{m})$) is dominated by the latter (i.e., $\min\{1-\frac{1}{\sqrt{m}}, \frac{h_1}{m} \} $),  we have the CR for $\mathcal{R}_2$ is at least $\CRTwoExp$.

 \begin{lemma}[Case 3: region $\mathcal R_3$] \label{thm:c4}
 Let $\mathcal R_3 =\{(h,\ell):  h\ge \mzero , \ell \ge \mzero  \}$.  Then, $\inf_{(h, \ell)\in \mathcal R_3}~\E_{\mathbf n}
    \big[\frac{\E[\rewGen]}{\textsc{opt}(\mathbf n)}\big]$ is at least 
\begin{align*}
   \left\{\begin{matrix}
    \; \overline{\text{CR}_3}, \quad & m \geq\mTwo , &\\  
    \;  \underline{\text{CR}_3},\quad  & m < \mTwo . & \\ 
    \end{matrix}\right.\,,
\end{align*}
where 
\begin{align*}\overline{\text{CR}_3}&=\CRThreeExpOne\, \\
\underline{\text{CR}_3}&=\CRThreeEXPTwo\,.
\end{align*}
Here, $\V=1-2(1-p)^{\sqrt{m}} $, $\X= \min\left\{1-\frac{1}{\sqrt{m}}, \frac{\mOne }{m} \right\}$, $\mOne =\mOneExp$,  $\beta=0.4215 \cdot \frac{p^2+(1-p)^2}{p(1-p)}$, $\hzero = \hzeroExp$,  
$\mTwo = \mTwoExp$, $\ell_1= \sqrt{p(1-p)\ell_0} +\frac{\beta}{\sqrt{m}} +\ell_0 p$, and $\ell_0=\ellzeroExp$.
\end{lemma}

\begin{remark} \label{remark:CR3}
We provide the exact formula for $\ell_0$ and $\mTwo$ here. Note that  $\ell_0= \ellzeroExp$ can be written as $\frac{m}{1-p}+\frac{1+\sqrt{4m(1-p)+1}}{2(1-p^2)}$. In addition, $\mTwo  = \mTwoExp$ can be written as 
$\mTwo =(\frac{b+\sqrt{b^2+4ac}}{2a})^4$, where  $a=(4-3p)p$, $b=(1-p)(4-3p)^{\frac{1}{2}}$ and $c=\frac{4(1-p)^2}{(r_1-r_2)^2}$. Further note that as $m \to \infty$, $\inf_{(h, \ell)\in \mathcal R_3}~\E_{\mathbf n}
    \big[\frac{\E[\rewGen]}{\textsc{opt}(\mathbf n)}\big]=\overline{\text{CR}_3}=1-\max\{\Theta(1/\sqrt{m}), \Theta(1/(p^2m))  \}$. 
\end{remark}

Lemma \ref{thm:c4} shows that when $(h, \ell)$ falls into region $\mathcal R_3$,   the CR of Algorithm \ref{alg:ptas} is lower bounded by $\text{CR}_3$. From the expression of $\text{CR}_3$ and Remark \ref{remark:CR3}, the asymptotic bound in region $\mathcal R_3$ is $1-\max\{\Theta(1/\sqrt{m}), \Theta(1/(p^2m))  \}$.  
This shows that our algorithm has a better asymptotic CR in  region $\mathcal R_3$, compared with region $\mathcal R_1$. (Recall that the asymptotic CR in region $\mathcal R_1$  is $1-\Theta(1/p\sqrt{m})$ and $p=\omega(1/\sqrt{m})$.) The fact that region $\mathcal R_1$ is more challenging is because  under region $\mathcal R_1$, the number of type $1$ agents, i.e., $h$, is small and hence the sample information may not be informative enough.

To show Lemma \ref{thm:c4}, the main idea is similar to the one in region $\mathcal R_2$. As the number of type $1$ agents is still large in $\mathcal R_3$, to bound the loss caused by over- or under-protecting type $1$ agents, we can use concentration inequalities to show that either the bad situation (i.e., over- or under-protecting type $1$ agents) does not happen or if it happens, it does not have  a significant amount of loss.

However, as the number of type $2$ agents is also large in $\mathcal R_3$, to bound the loss caused by wrongfully protecting type $2$ agents, we need to further split the analysis into two parts based on  the number of resources $m$.  When $m \geq \mTwo$, we show that the probability that Algorithm \ref{alg:ptas} wrongfully protects type $2$ agents is small because, in this case, we have enough samples to get an accurate estimation of the reward of each type.  When  $m < \mTwo$,  we directly compute this loss. The main challenge is still  dealing with the expectation of a fraction of two random variables. To handle this challenge, we again split region $\mathcal R_3$ into two sub-regions:  $\underline{\mathcal{R}_3}$ and $\overline{\mathcal{R}_3}$, where
$\overline{\mathcal{R}_3} =\{(h, \ell): (h, \ell)\in \mathcal R_3, \ell > \ell_0 \}$ and $\underline{\mathcal{R}_3} =\{(h, \ell): (h, \ell)\in \mathcal R_3, \ell \le \ell_0 \}$.
 By providing a proper bound for the  optimal clairvoyant solution  in the denominator, we then show that the CR under region $\underline{\mathcal{R}_3}$ is at least $\underline {\text{CR}_3}$ and the CR under region $\overline{\mathcal{R}_3}$ is at least $\overline {\text{CR}_3}$.

Next, by putting the results of Lemmas \ref{lem:c1}, \ref{lem:c3}, and \ref{thm:c4}, we get the following theorem.

\begin{theorem} [Competitive ratio of Algorithm \ref{alg:ptas}]\label{thm:2}
Consider the model presented in Section \ref{sec:model}, where the expected rewards of $r_i$, $i\in \{1, 2\}$ are unknown to the decision-maker. Let $\hzero =\hzeroExp= \Theta(\frac{m}{1-p})$, $\mOne =\mOneExp= 1-\Theta(\frac{(1+\sqrt{p})}{\sqrt{m}})$, $\V=1-2(1-p)^{\sqrt{m}} =1-\Theta((1-p)^{\sqrt{m}})$, $\X= \min\left\{1-\frac{1}{\sqrt{m}}, \frac{\mOne }{m} \right\}=1-\Theta(\frac{1}{\sqrt{m}})$, $\ell_0=\ellzeroExp=\Theta(\frac{m}{1-p})$, $\ell_1= \sqrt{p(1-p)\ell_0} +\frac{\beta}{\sqrt{m}} +\ell_0 p=1-\Theta(\frac{(1+\sqrt{p})}{\sqrt{m}})$, $\alpha=\frac{r_2}{r_1}>0$, and $\beta=0.4215 \cdot \frac{p^2+(1-p)^2}{p(1-p)}$. Then, for any $m\ge 2$
and $p\in (0,1)$,  the CR of Algorithm \ref{alg:ptas}, denoted by $\text{CR}_{A}$, is lower bounded by:
\begin{align*}
    \text{CR}_{A} \geq 
   \left\{\begin{matrix}
    \; \min\left\{\CROne , \CRTwo , \overline{\text{CR}_3}\right\},\quad  & m \geq\mTwo , &\\  
    \;  \min\left\{\CROne , \CRTwo , \underline{\text{CR}_3}\right\}, \quad & m < \mTwo . & \\ 
    \end{matrix}\right.\,,
\end{align*}
where $\mTwo = \mTwoExp$, and 
\begin{align*}
\CROne &=\CROneExp\,,
\\
\CRTwo &=\CRTwoExp\,,\\
\overline{\text{CR}_3}&=\CRThreeExpOne\, \\
\underline{\text{CR}_3}&=\CRThreeEXPTwo\,.
\end{align*}
\end{theorem}

To provide insights into the lower bound in Theorem \ref{thm:2}, before we formally present the asymptotic CR of Algorithm \ref{alg:ptas}, we  provide an example to see how the CR of Algorithm \ref{alg:ptas} vary with $m$ and $p$. 

\subsection{ A Simple Example: Evaluating the Lower Bound in Theorem \ref{thm:2}}
\begin{example} \label{example:CR}
Consider a setting in which $r_1=0.9$, and  $r_2$ is uniformly drawn from the interval $(0.5,0.9)$; that is, $r_2 \sim \text{Uniform}(0.5,0.9)$. 
Here, we would like to study how the true CR (which will be defined below) of Algorithm \ref{alg:ptas} and its lower bound in Theorem \ref{thm:2} change as a function of $m$ and $p$. To do so, we consider two scenarios.  In first scenario, we fix the sampling probability $p$ at $0.3$, and we let the number of resources  $m\in\{10,35,\ldots,360,385\}$. For each value of $m$, we generate $300$ instances where an instance is determined by a realized value of  $r_2\sim \text{Uniform}(0.5, 0.9)$. 
In the second scenario, we fix the number of resources $m$ to be $200$
and we let the sampling probability $p\in \{0.1, 0.15, \ldots, 0.55,0.6 \}$. For each value of $p$,  we generate $300$ instances where again  an instance is determined by the realized value of $r_2$.

For each instance of both scenarios,  we then compute three quantities : (i) the true CR of Algorithm \ref{alg:ptas}, (ii) the lower bound on CR of Algorithm \ref{alg:ptas}, which is presented in Theorem \ref{thm:2}, and (iii) $1/(2-r_2/r_1) = 1/(2-\alpha)$. The true CR of the algorithm is obtained by considering all the possible values for $h$ and $\ell$, and only focusing  on ordered arrival sequences in which type $2$ agents arrive first.  Note that focusing  on ordered arrival sequences is without loss of generality as shown in Lemma  \ref{lm:order}.  The last quantity $1/(2-\alpha)$ is the optimal CR when $r_i$, $i\in[2]$, is known, but no sample information is  available. We use this quantity as a benchmark.    

Figure \ref{fig:sub-first} shows the expected value of these three quantities versus $m$ in the first scenario where we fix the sampling probability $p= 0.3$. Here, the expectation is with respect to $r_2$.  
The figure shows both the true CR and its lower bound in Theorem \ref{thm:2} improves as $m$ increases. Our lower bound, however,  gets tighter as $m$ increases. Nonetheless, our lower bound is quite loose when $m$ is small. This is mainly because of our loose lower bound for $\text{CR}_1$, where we recall that $\text{CR}_1$ is the lower bound on the CR of the algorithm when $h\le \sqrt{m}$.  To characterize $\text{CR}_1$, we present a universal lower bound for any realization $\mathbf n$, which leads to a loose bound; see Lemma \ref{lem:c1} and its proof.

What is quite interesting is that the true CR of the algorithm  exceeds the benchmark of $1/(2-\alpha)$ when $m\ge 35$.  Recall that  the benchmark is the optimal CR when $r_i$'s are known, but sample information is not available. This highlights that for large enough $m$, the value of the sample information outweighs the drawbacks of not knowing the true rewards. 

Figure \ref{fig:sub-second} shows the average of these three quantities versus $p$ in the  second scenario where $m =200$. The figure shows that our lower bound gets tighter as $p$ increases. When $p \geq 0.4$, we can find that the lower bound is almost equal to the true CR of  Algorithm \ref{alg:ptas}. 
In addition, with $m=200$, even if the sampling probability is as small as $0.1$, the true CR of Algorithm \ref{alg:ptas} is still higher than the benchmark. 
In addition, we observe  that both the lower bound on the  CR (given by Theorem \ref{thm:2}) and the true CR of Algorithm \ref{alg:ptas} are larger in the high sampling probability case. 

\begin{figure}[ht] 
    \centering
  \includegraphics[width=0.47\linewidth]{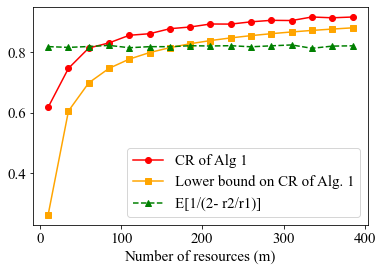}  
  \caption{ $p =0.3$}
  \label{fig:sub-first}
  \centering
  \includegraphics[width=0.47\linewidth]{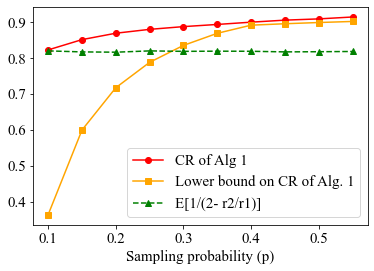}  
  \caption{ $m=200$}
  \label{fig:sub-second}

\caption{The CR of Algorithm \ref{alg:ptas}, the lower bound on the CR of Algorithm \ref{alg:ptas} presented in Theorem \ref{thm:2}, and the bound in \cite{ball2009toward} (i.e., $\E[1/(2-\alpha)]$). Here, $r_1=0.9$, and  $r_2$ is uniformly drawn from the interval $(0.5,0.9)$.}
\label{fig:fig}
\end{figure}

\end{example}

\subsection{Asymptotic Competitive Ratio of Algorithm \ref{alg:ptas}}
Here, we present an asymptotic CR of Algorithm \ref{alg:ptas} as a function of $m$ and $p$ when $m$ goes to infinity. 

\begin{proposition}[Asymptotic Competitive Ratio of Algorithm \ref{alg:ptas}] \label{thm:3}
Consider the model presented in Section \ref{sec:model}, where the expected rewards of $r_i$, $i\in \{1, 2\}$ is unknown to the decision-maker. As $m$ goes to infinity, the worst-case CR of Algorithm \ref{alg:ptas} scales with $1- \Theta(1/(p\sqrt{m}))$.  
\end{proposition}

Proposition \ref{thm:3} is shown in Section \ref{subsec:proofthm3} in the appendix. Proposition \ref{thm:3} shows that the CR of Algorithm \ref{alg:ptas} goes to one as $m$ goes to infinity if $p=\omega(1/\sqrt{m})$. The asymptotic CR of the algorithm (i.e., $1-\Theta(1/(p\sqrt{m}))$) is indeed tight for any $p=\omega(1/\sqrt{m})$; see Section \ref{sec:upper} for an upper bound on the CR of any deterministic  and randomized algorithms when $p=\omega(1/\sqrt{m})$. \footnote{For the case where  $p= O(1/\sqrt{m})$, in Section \ref{sec:upper}, we further present a non-asymptotic upper bound. 
This bound also shows that   when $p=O(1/{m})$, no randomized and deterministic algorithms can be asymptotically optimal.} 
Note that even  when $p$ is a constant, it is not obvious why achieving asymptotic optimality is possible using  Algorithm \ref{alg:ptas}. Suppose that $m$ is large and $p$ is constant. If the adversary chooses a large market size (i.e., when both $h$ and $\ell$ are  large), for a constant value of $p$, the number of agents in the sample set is large, providing significant information. However, the market size is chosen adversarially and hence it can be small either for type $1$ or $2$ agents, and in such a case, the sample information is not very informative. 
(Recall that the asymptotic CR of Algorithm \ref{alg:ptas} is governed by region $\mathcal R_1$, where the number of type $1$ agent is low.) Nevertheless, it turns out Algorithm \ref{alg:ptas} can very well handle the difficult cases where the market size is small, achieving an asymptotically optimal CR.

\section{Upper Bound on the Competitive Ratio of All Deterministic and Randomized Algorithms} \label{sec:upper}
In our problem setting, the average rewards of agents are unknown, but the decision-maker can learn about them  with the help of   the sample information. 
In this section, we relax the problem a little bit by
considering  setting in which the decision-maker knows $r_1$ and $r_2$ in advance.
We present an upper bound on the CR of any deterministic and randomized algorithm in the aforementioned relaxed setting when the number of resources $m$ goes to infinity. Clearly, any upper bound for this relaxed setting is  also a valid upper bound for our original setting in Section \ref{sec:model}.

First, in Section \ref{sec:upperlargep}, we consider the case where $p=\omega(1/\sqrt{m})$.
For the case of $p=\omega(1/\sqrt{m})$, we present an upper bound on the CR of any deterministic and randomized  algorithms under the relaxed setting. We show that the  CR of any randomized and deterministic algorithm is upper bounded by $1-\Theta({1}/{(p\sqrt{m})})$ when $p=\omega(1/\sqrt{m})$, implying the asymptotic optimality of  Algorithm \ref{alg:ptas}; see Proposition \ref{thm:3}. 

Put differently, the asymptotic CR of our algorithm---which does not use any feedback throughout the allocation period to tune its estimates for $r_i$'s---matches the upper bound on the CR of any deterministic and randomized algorithm in the relaxed setting where $r_i$'s are known in advance. This implies that asymptotically, feedback during the allocation period does not add any value in improving the CR of algorithms. The fact that feedback throughout the allocation period is not useful   is due to the adversarial nature of arrival sequences. Under adversarial arrivals, the decision-maker cannot control for what type of agents she receives feedback; this is governed by the adversarially-chosen order. 

Second, in Section \ref{sec:uppersmallp}, we consider the case where $p=O(1/\sqrt{m})$. For this case where the sampling probability is quite small even for large values of $m$, we present a non-asymptotic upper bound. Our bound  shows that when $p=O(1/m)$, 
 the CR of any randomized and deterministic algorithm is upper bounded by $1/(2-\alpha)$ as $m$ increases. This is mainly due to the lack of adequate sample information when $p =O(1/{m})$.

\subsection{Upper Bound for $p =\omega(1/\sqrt{m})$}\label{sec:upperlargep}
Consider the case   of $p =\omega(1/\sqrt{m})$. Under this case, in Theorem \ref{thm:upperbound}, we first present an upper bound of $1-\Theta({1}/{(p\sqrt{m})})$ on the CR of any deterministic algorithm under the relaxed setting. Second, in Theorem \ref{thm:upperboundrandom}, we show that the same upper bound holds for any randomized algorithm.

\begin{theorem}[Upper Bound of Any Deterministic Algorithm for $p = \omega(1/\sqrt{m})$] \label{thm:upperbound}
Consider a special case of our model, presented in Section \ref{sec:model}, where the expected rewards of $r_i$, $i\in \{1, 2\}$ is known to the decision-maker. When $p = \omega(1/\sqrt{m})$, as the number of resources $m$ goes to infinity,  any deterministic algorithm has the CR of at most $1-\Theta(1/({p\sqrt{m}}))$.
\end{theorem}

The proof of Theorem \ref{thm:upperbound} is presented  in Section \ref{sec:proof:upperbound}. 
 To show Theorem \ref{thm:upperbound}, we construct the following input family $\textbf{F}$: 
Let  $\underline{h}=0$ and  $\bar{h}=pm$. The input family $\textbf{F}$  contains all $(h, \ell)$ such that $h \in [\underline{h},\bar{h}]$ and $\ell = \frac{10000\cdot m}{p}$.  For any $h \in [\underline{h},\bar{h}]$, we then denote $I_h$ as a random arrival sequence under which $n_2$ type $2$ agents arrive followed by  $n_1$ type $1$ agents, where we recall that   $n_1 \sim \text{Bin}(h,1-p)$, $n_2 \sim \text{Bin}(\ell,1-p)$.
We characterize an upper bound on the CR of any deterministic algorithm  under the family $\textbf{F}$.

In this family, because $\ell\gg m$, we know  there will be more than $m$ type $2$ agents showing up. Therefore, the number of type $2$ agents in the sample does not impact the acceptance/rejection  decisions. Then, given that the online arrival sequences $I_h$ are all ordered,  any deterministic algorithm has to decide about how many type $2$ agents they accept provided that they observe $s_1$ samples from type $1$ agents. Put differently, any deterministic algorithm can be represented by a mapping that maps $s_1$ to the number of type $2$ agents it accepts. 
The proof of Theorem \ref{thm:upperbound} then shows that the best CR under any such mapping is $1-\Theta(1/(p\sqrt{m}))$. The main challenge in showing this result is  characterizing the optimal mapping. To overcome this, instead of characterizing the optimal mapping,  we first construct a  specific mapping  under which  upon observing $s_1$ type $1$ agents, we accept $m-s_1\frac{1-p}{p}$ type $2$ agents. Observe that $s_1\frac{1-p}{p}$ is the expectation of the number of type $1$ agents who will arrive given that there are $s_1$ type $1$ agents in the sample. 
Lemma \ref{lem:specific z} shows that under this mapping, the CR is at most $1-\Theta(1/(p\sqrt{m}))$.
To complete the proof of Theorem \ref{thm:upperbound}, we then compare any other mappings with this specific mapping; see Lemma \ref{lem:other_mappings}.

Next, we derive an upper bound on the CR of  all randomized algorithms.

\begin{theorem}[Upper Bound of Any Randomized Algorithms for $p = \omega(1/\sqrt{m})$] \label{thm:upperboundrandom}
Consider a special case of our model, presented in Section \ref{sec:model}, where the expected rewards of $r_i$, $i\in \{1, 2\}$ is known to the decision-maker. Then, when $p = \omega(1/\sqrt{m})$, as the number of resources $m$ goes to infinity,  any randomized algorithm has the CR of at most $1-\Theta(1/(p\sqrt{m}))$. 
\end{theorem}

To show Theorem \ref{thm:upperboundrandom}, unfortunately, we cannot use Von Neuman/Yao principle \cite{seiden2000guessing}. This is because 
in our setting, 
even when the input $(h,\ell)$ is realized, due to our sampling procedure, the online arrival sequence is still random. This is different from Von Neuman/Yao principle, because Von Neuman/Yao principle can only be applied to the model without any randomness. Nonetheless, in Lemma \ref{lem:yao}, we derive a result similar to the Von Neuman/Yao principle that can be applied to our setting.
We then  apply Lemma \ref{lem:yao} by  constructing  a distribution over the input family $\textbf{F}$ introduced above. This leads to the desired upper bound. See Section \ref{sec:proof:thm:upperboundrandom} for the proof of Theorem \ref{thm:upperboundrandom}.

\subsection{Non-asymptotic Upper Bound for any $p$} \label{sec:uppersmallp}

In the previous section, we give a tight asymptotic upper bound for $p=\omega(1/\sqrt{m})$. There are two remaining questions: (i) Does there exist any non-asymptotic upper bound (even if it is loose)?  (ii) What is the asymptotic upper bound when $p$ goes to zero very fast as $m$ grows? In this section, we try to answer these two questions.

The following theorem (Theorem \ref{thm:uppersmallp}) presents a non-asymptotic upper bound on the CR of any  deterministic or randomized algorithms as a   function $p$ and $m$. The theorem further  shows that when functions $p = O(1/m)$, achieving asymptotic optimality is not possible. With $p=O(1/m)$, the asymptotic upper bound is equal to $\frac{1}{2-\alpha}$, which is the bound in \cite{ball2009toward}. This shows that  when $p$ goes to $0$ very fast, the sample information is not adequate enough to make an impact. 
We  note that the upper bound in 
Theorem \ref{thm:uppersmallp} can be loose because  to compute the upper bound, we consider a relaxed setting where the rewards are known to the decision-maker. 
Even so, in light of the impossibility results in Theorem \ref{thm:uppersmallp}, in Section \ref{sec:algpsmall},  we present a simple algorithm whose CR is $0.5$ for any value of $p$ and $m$. While this CR does not match the upper bound in Theorem \ref{thm:uppersmallp}, the CR is tight when $r_2\rightarrow 0$ and $p=O(1/m)$. 

\begin{theorem}[Non-Asymptotic Upper Bound of Any Algorithms for any  $p$] \label{thm:uppersmallp}
Any deterministic or randomized algorithm cannot achieve a CR better than $\frac{m}{m+\min\{\tilde{h},m\}(1-\alpha)}+\frac{1}{\tilde{h}^2}$, where $\tilde{h}$ is the smallest integer such that $\tilde{h}p+\sqrt{\tilde{h} \log(\tilde{h})} \geq 1$, and 
 $\alpha = r_2/r_1$. 
In addition, when $p=O(1/m)$,  no algorithm can obtain a CR better than $\frac{1}{2-\alpha}$ as $m \to \infty$.
\end{theorem}

\begin{figure}\centering
      \includegraphics[width=0.4\textwidth]{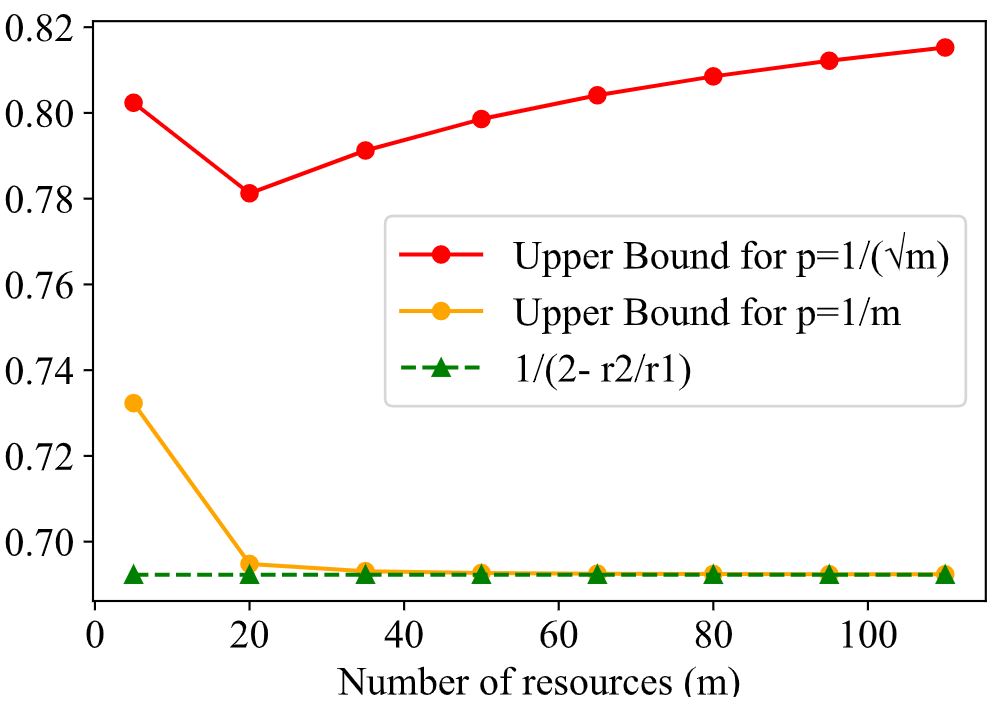}  
      \caption{}
      \label{fig:figupper}
       \includegraphics[width=0.4\textwidth]{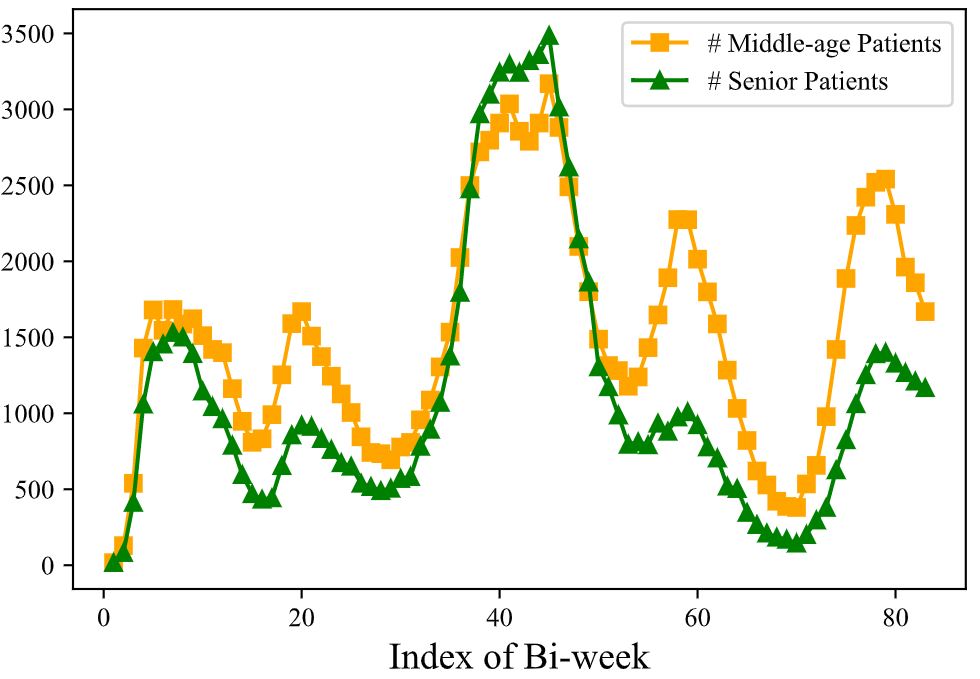}  
         \caption{}
       \label{fig:numberofpatients}
     \caption{(a) Upper bound in Theorem \ref{thm:uppersmallp} for $p=1/\sqrt{m}$ and $p=1/m$, and the bound  in \cite{ball2009toward}. Here, $r_1=0.9$, and $r_2=0.5$. (b) Number of middle-age patients and senior patients over weeks}
     \end{figure}

The proof is in Section \ref{subsec:proofoflowp}.  Figure \ref{fig:figupper}  depicts the upper bound in Theorem \ref{thm:uppersmallp} versus $m$ with
$p=1/\sqrt{m}$ and $p =1/m$, respectively. Here, $r_1 =0.9$ and $r_2 = 0.5$. As a benchmark, in both figures, we  also include  $1/(2-r_2/r_1)$, which is the bound in \cite{ball2009toward}. As expected, we observe that when $p=1/m$, the upper bound   converges to the benchmark. When $p= 1/\sqrt{m}$, however, the upper bound goes to $1$ because this bound is loose.  

To show Theorem \ref{thm:uppersmallp}, we construct an input family $\mathbf{F}_1$ as follows: The input family $\mathbf{F}_1$ contains all $(h,\ell)$ such that $h \in [0,\lfloor m/p \rfloor]$ and $\ell=\lceil \frac{10000\cdot m}{p} \rceil$. For any $h \in [0,\lfloor m/p \rfloor]$, we then denote $I_h$ as an ordered  random arrival sequence under which $n_2\sim \text{Bin}(\ell,1-p)$ type $2$ agents arrive followed by $n_1\sim \text{Bin}(h,1-p)$ type $1$ agents.  
We then show that no algorithm can do well simultaneously on arrival sequences  $I_{h}$ for any $h \in [0,\tilde{h}]$, where $\tilde{h}$ is the smallest integer such that $\tilde{h}p+\sqrt{\tilde{h} \log(\tilde{h})} \geq 1$. (Asymptotically, $\tilde{h}=\Theta(1/p)$) Under  these arrival sequences, with high probability, $s_1= 0$. But these sequences vary a lot in terms of the number of type $1$ agents. Yet it is not possible to distinguish them using the sample information.

\section{A Case Study on COVID-19 Dataset} \label{sec:simulation}

In this section, we do a case study in which we apply Algorithm \ref{alg:ptas} to the COVID-19-associated hospitalizations problem in the US. We use Algorithm \ref{alg:ptas} to allocate hospital resources (e.g., a medicine) to  arriving COVID-19 patients sequentially. This case study allows us to evaluate Algorithm \ref{alg:ptas} in real-world inspired settings. It further allows us to see how robust the algorithm is to the estimation errors in the sampling probability $p$. We observe that the average CR of the algorithm even for small values of $m$ ($m<50$) 
is at least $0.88$, and the algorithm maintains its performance when it does not know the true value of the sampling probability. Further, we  demonstrate the value of the sample information and running a test period by comparing our algorithm with a benchmark that does not use the sample information.

\textbf{Dataset.}  Here, we use the ``Laboratory-Confirmed COVID-19-Associated Hospitalizations" dataset.   This dataset contains the number of bi-weekly cases of COVID-19-associated hospitalizations in the US from March $7^{\text{th}}$, $2020$ to February  $5^{\text{th}}$, $2022$ across five age groups  of patients: 0-4 years, 5-17 years, 18-49 years, 50-64 years, and 65+ years.  Because the number of 0-17 years old  patients is small, we discard them in our analysis. We then consider two groups/types of patients: middle-age patients (18-64 years), and senior patients (65+ years).  In our studies,  the middle-age patients are considered to be  the low-reward agents (i.e., type $2$ is our setting), and senior  patients are considered to be high-reward agents (i.e., type $1$  in our setting). This assumption is partly motivated by high death-rate of senior patients when contracting COVID-19  \citep{covid}.

\textbf{Simulation setup.} In our studies, given the description of our dataset, we consider the allocation problem of a  resource (e.g., certain medicine) over the course of two weeks, where each resource can be assigned to at most one patient and each patient needs one unit of the resource.
For each such period,  we determine the number of hospitalized patients of the two aforementioned types (i.e., $h$ and $\ell$) using our dataset. (Note that we have $83$ periods in our dataset.) The value of $h$ and $\ell$ over $83$ periods can be found in Figure \ref{fig:numberofpatients}. At the beginning of each of the periods,  we observe the (realized) effectiveness of the resource for a sample of hospitalized patients, where each of the $(h+\ell)$  hospitalized patients falls into the sample with probability $p$. These sampled patients are, for example,  the patients that arrive at the beginning of the two-week period during which rationing the resource has not started yet. Let $n_1$ and $n_2$ be the number of hospitalized patients of types $1$ and $2$ who did not fall into the sample. We assume that the order over these $n_1+n_2$ patients are completely random. 
That is, we consider uniform permutations over these patients, modeling more realistic scenarios where the order over patients are not chosen adversarially.

We consider $9$ problem classes where each problem class is determined by two parameters $(p, \gamma)$. Here, $p\in \{0.1, 0.15, 0.2\}$ is the sampling probability, and $\gamma \in \{0.3, 0.5, 0.7\}$  determines the scarcity of the resource. In particular, we set $m = \gamma\cdot(\hat h+\hat \ell)$, where $\hat h$ and $\hat \ell$ are respectively  the number of type $1$ and type $2$ patients in the previous period (i.e., the last two weeks). By setting $m$ in this way, we would like to capture the hospital's inventory planning decisions based on  the most  recent observed demand. Further, in choosing the values for the number of resources $m$, we let $\gamma$ be less than one. This allows us to model scenarios where resources are scarce. When resources are not scarce,  the problem is not challenging.

For each problem class, we generate $1000$  instances, where each instance is determined by $(n_1, n_2)$  and an order over $n_1+n_2$ patients.  For each of the $9$ problem classes, we evaluate the performance of four algorithms: (i)
 Algorithm \ref{alg:ptas} that has access to the true value of $p$, (ii) an algorithm that does not use sample information (which we will  define later in this  section), (iii) an adaptive version of Algorithm \ref{alg:ptas} that has access to the true value of $p$ (see Section \ref{sec:alg} for the description of this algorithm),  and (iv) Algorithm \ref{alg:ptas} that has only access to the noisy estimate of $p$, denoted by $\hat p$, and uses $\hat p$ in place of $p$. Here, $\hat p\sim \text{Uniform}(\underline p, \bar p)$,  $\E[\hat p] =p$, and standard error of $(\hat p)$ is set to $0.3p$. For example, if $p=0.1$, we have  $\underline p \approx 0.05$ and $\bar p \approx 0.15$. 
 
The second algorithm is similar to Algorithm \ref{alg:ptas} in the sense that it also assigns protection levels to either type $1$ or type $2$ agents. But, unlike Algorithm \ref{alg:ptas}, it does not use the sample information to set the protection level and  decide what type of agents  to protect.  
This algorithm, which is parameterized by $x$, works as follows: with a probability $1/2$, it assigns a protection level of $x \in \{0,1,\ldots,m\}$ to type $1$ agents. Otherwise, it assigns a protection level of $x$ to type $2$ agents. 
 The optimal (CR-maximizing) value for  $x$ can be calculated as a function of $\alpha=r_2/r_1$ and $m$.\footnote{The optimal value for $x$ is $\frac{2-2\alpha}{2-\alpha}m$.} Although without the  sample information, rewards $r_1$ and $r_2$, and hence $\alpha$ are unknown, to give the algorithm an extra benefit, in our case study, we set $x$ optimally. We highlight that we use this algorithm in our case study to shed light on the value of the sample information. 

Finally, we note that in all of our evaluations, we assume that 
$r_1=0.6$ and  $r_2=0.2$. 
We also try to let $r_i$'s take other values: the  take-away messages of our case study   do not change.

\textbf{Performance evaluation.}
Table \ref{tab:multicol} presents the worst-case and average CR of the four algorithms across the 82 \footnote{Since in each two-week period, we set $m$ as a function of demand in the previous two-week period, we cannot have the results for the first two-week period.} allocation periods for each of our $9$ problem classes, where in each allocation period, the CR is the worse case (minimum) reward over optimal clairvoyant  reward among the $1000$ instances. Then, the worst-case and average CR reported in Table \ref{tab:multicol}
are computed by taking  the minimum and average  over the CR of  82 allocation periods. 

Interestingly, for all problem classes, the average CR of our algorithm and its noisy version is between  $0.88$ and $0.97$ while the average CR of the algorithm that does not use the sample information  can be as low as $0.66$. We also observe  that Algorithm \ref{alg:ptas} and the adaptive version have almost the same CR. This is because, under the adaptive algorithm, in average (in worst case respectively), we change our decision about whom to protect less than $1$ time (twice respectively). 

In addition, for all four algorithms, if we fix $\gamma$, we observe  that by increasing  the sampling probability $p$, both the average  and worst-case CR increase. This observation is consistent with the result in Proposition \ref{thm:3}.
However, if we fix the sampling probability $p$, there is no clear pattern of how the CR changes when $\gamma$ increases.

Overall, we observe that Algorithm \ref{alg:ptas} with true value of $p$ and  Algorithm \ref{alg:ptas} with noisy  $\hat p$ significantly outperform the naive algorithm that does not take advantage of the sample information.   

This suggests the importance using the sample information even when  we do not have access to the true value of  $p$.

One surprising result is that: among all 9 problem classes and 4 algorithms, the worst CR happens in the second allocation period among the $82$ periods. (Recall that we did not analyze the first allocation period because how we set $m$.) From Figure \ref{fig:numberofpatients}, we can find that the second allocation period has the smallest number of type 1 and 2 patients. Hence, this observation is consistent  with the  insights from our theoretical result, which  show that  when the market size is small, the CR is also small.

\begin{table}[htb]
\caption{The CR of four algorithms in Section \ref{sec:simulation}. The standard error of average CRs of the four algorithms is less than $0.001$. }
\begin{center}\footnotesize
\begin{tabular}{c|cc||c|c|c||c|c|c||c|c|c|}
    \cline{2-12}
    & Value of $p$ &   &  \multicolumn{3}{c}{p=0.1} &   \multicolumn{3}{|c}{p=0.15} &  \multicolumn{3}{|c|}{p=0.2}\\    \cline{2-12}
    & Value of $\gamma$ & & $\gamma=0.3$ & $\gamma=0.5$ & $\gamma=0.7$ &$\gamma=0.3$ & $\gamma=0.5$ & $\gamma=0.7$ &$\gamma=0.3$ & $\gamma=0.5$ & $\gamma=0.7$\\ \hline
        \multicolumn{1}{|c|}{Alg. 1}  & Avg. CR & & 0.977& 0.953 & 0.961 & 0.978 & 0.962 & 0.968 & 0.978 & 0.965 & 0.973\\
   \multicolumn{1}{|c|}{with true $p$}  &Worst CR& & 0.791 & 0.763 & 0.811 & 0.808 & 0.787 & 0.824 & 0.811 & 0.771 & 0.862\\ \hline
   \multicolumn{1}{|c|}{Alg that doesn't  } & Avg. CR & & 0.656 & 0.693 & 0.838 & 0.657 & 0.711 & 0.872 &0.671 &0.737 &0.912 \\
    \multicolumn{1}{|c|}{use samples} &Worst CR & & 0.567 & 0.573& 0.767 & 0.589 & 0.649 & 0.810 &0.593 &0.671 &0.860 \\\hline
   \multicolumn{1}{|c|}{Adaptive version of}  & Avg. CR & & 0.979& 0.955 & 0.962 & 0.980 & 0.962 & 0.968 & 0.979 & 0.967 & 0.973\\
   \multicolumn{1}{|c|}{Alg. 1 with true $p$}  &Worst CR& & 0.800 & 0.766 & 0.811 & 0.837 & 0.792 & 0.831 & 0.866 & 0.791 & 0.864\\ \hline
      \multicolumn{1}{|c|}{Alg. 1 with} & Avg. CR & & 0.953 & 0.892 & 0.888 & 0.944 & 0.889 & 0.910 & 0.931 & 0.881 & 0.917\\
    \multicolumn{1}{|c|}{noisy $p$} &Worst CR& & 0.790 & 0.723 & 0.775 & 0.817 & 0.800 & 0.811 & 0.799 & 0.781& 0.849 \\ \hline
\end{tabular}
\end{center}
\label{tab:multicol}
\end{table}

\section{Extensions and Discussions } \label{sec:exten}

In this section, we discuss some extensions of our model.  In Section \ref{sec:limitedsample}, we study a setting where the number of resources during the test period is limited, and in Section \ref{sec:realcr}, we investigate a strong clairvoyant benchmark that knows the realized reward of agents.
In Section \ref{subsec:heterop}, we discuss how to modify Algorithm \ref{alg:ptas} when the sampling probability is different across different types of agents. We further investigate how the modified version of  Algorithm \ref{alg:ptas} performs in the worst case. In Section \ref{subsec:nagent}, we present an extension of  Algorithm \ref{alg:ptas} when there are more than two types of agents.

\subsection{Capacity Constraints during the Test Period} \label{sec:limitedsample}
So far, we have assumed that there is no capacity constraint during the test period in the sense that any agent in the test period gets one unit of the resource. However, in practice, decision-makers may limit the number of available resources during the test period due to the cost of testing the market. Here, we show that when the number of available resources during the test period is $\omega(\sqrt{m})$ (e.g., $\Theta(m^{1/2+\epsilon})$ for any $\epsilon> 0$) the same asymptomatic CR, presented in Proposition \ref{thm:3}, continues to hold. 

Let $m_t$ and $m_a =m$ be the number of available resources during the test and allocation periods, respectively. 
Further, recall that $s_i$, $i\in \{1, 2\}$, is the number of agents of type $i$ in the sample, where we assume that  $s_i$'s are known to the decision-maker when enforcing the capacity constraint.\footnote{That is,  all agents who respond to the outreach program arrive first. After their arrivals, we decide about allocating resources to them.} To enforce the capacity constraint, when $s_1+s_2\le m_t$, clearly, we can accept all the (sampled) agents in the test period. On the other hand, when  $s_1+s_2> m_t$ with 
 $s_1, s_2> m_t/2$, we accept $m_t/2$ agents of type $i\in \{1, 2\}$. Finally, when  $s_1+s_2> m_t$ with either $s_1$ or $s_2$ (not both) is greater than $m_t/2$, we accept a total of $m_t$ agents, while making sure all the agent of the type with a lower number of samples is accepted. At the end of the test period, for simplicity, we assume that the leftover units do not carry over to the allocation period.

The following proposition presents the asymptotic CR of Algorithm \ref{alg:ptas} with a slight modification: in calculating  the empirical average awards $\widehat{r_i}$'s, instead of using $s_i$'s, we use the number of accepted agents of type $i$ at the end of the test period. 

\begin{proposition}[Capacity Constraint during the Test Period] \label{prop:limitedm}
Consider the model presented in Section \ref{sec:limitedsample} where the number of resources during the test period is at most $m_t$. Then,  when $m_t =\omega(\sqrt{m})$, the worst-case CR of a modified version of  Algorithm \ref{alg:ptas},  stated above,  scales with $1- \Theta(1/(p\sqrt{m}))$ as $m$ goes to infinity. Here, $m$ is the number of resources during the allocation period. 
\end{proposition}

The proof can be found in Section \ref{sec:prooflimitm}. The proof is similar to that of Theorem \ref{thm:3}. As the main difference, we need to re-evaluate the probability of good event which is  needed  in the proof of region $\mathcal R_3$.  Proposition \ref{prop:limitedm} shows that even when the number of resources during the test period is quite small, Algorithm \ref{alg:ptas} maintains its good performance.

\subsection{ A Different Notion of Competitive Ratio using Realized Rewards} \label{sec:realcr}
In our setting, we compared the performance of our algorithm w.r.t. an optimal clairvoyant benchmark  that knows the arrival sequence and the expected reward of the agents in advance, but  does not know the realized rewards of the agents. See Equation \eqref{eq:CR_early} for the definition of our CR notion. This benchmark is indeed consistent with the benchmark used in \cite{ball2009toward}.  Nevertheless, one may wonder how our algorithm performs against a stronger benchmark that also knows the realized reward of the agents, denoted by $\mathcal{R}=\{R_t\}_{t\in [|I|]}$. (Here, when  $z_t =i$, the realized reward $R_t\sim \text{Ber}(r_i)$, and the realized reward is not known at the time of the decision.)
Such a clairvoyant benchmark only accepts agents with realized reward of one, where we recall that the realized reward of agents is either zero or one. 

Let $\textsc{opt-re}(\mathcal{R}, I)$ be the maximum realized reward under $\mathcal{R}$ and the arrival sequence $I$, where ``\textsc{re}" stands for ``realized."  Further, let $\rew\textsc{-re}_{\mathcal A}(\mathcal{R}, I, \boldsymbol \psi)$ be the realized reward of an algorithm $\mathcal A$ under  $\mathcal{R}$, arrival sequence $I$, and  the sample information $\boldsymbol \psi$. Then,  the CR of an algorithm $\mathcal A$ w.r.t. the strong benchmark $\textsc{opt-re}(\mathcal{R}, I)$, denoted by $\text{CR}\textsc{-re}$, is defined as 
\begin{align}\label{eq:CRdefreal}
\text{CR}\textsc{-re}_{\mathcal A}  =  \inf_{I}~\E
    \left[\inf_{\mathcal R}\frac{\E[\rew\textsc{-re}_{\mathcal A}(\mathcal{R}, I, \boldsymbol \psi)]}{\textsc{opt-re}(\mathcal{R}, I)}\right]\,,
\end{align}
where the inner expectation is w.r.t. any randomness in algorithm $\mathcal A$, and the outer expectation is w.r.t. the sample information and realized rewards $\mathcal R$. We refer to $\text{CR}\textsc{-re}_{\mathcal A}$ as the \emph{realized CR} of algorithm $\mathcal A$.

In this section, we first show that for any randomized and deterministic algorithm $\mathcal A$, the realized CR cannot exceed $r_2$ even when expected rewards are known. See Proposition \ref{prop:negativereal}.
This highlights that is is very challenging to beat the realized benchmark   $\textsc{opt-re}(\mathcal{R}, I)$. We then revisit  Example \ref{example:CR} in which we evaluate Algorithm \ref{alg:ptas}  using the realized CR. The example shows that the realized CR of Algorithm \ref{alg:ptas} converges to $r_2$ as $m$ increases, implying Algorithm \ref{alg:ptas} might be also optimal in terms of   the realized CR.  

\begin{proposition} [Upper  Bound on the Realized CR]\label{prop:negativereal}
Consider a relaxed version of our setting where  the expected rewards are known, but we still have access to sample information $\boldsymbol \psi$ that reveals some information about the market size $h$ and $\ell$. Then, the realized CR of any deterministic or randomized algorithm (per Equation \eqref{eq:CRdefreal}) is at most $r_2$.
\end{proposition}
The proof can be found in Appendix \ref{sec:proof:prop:negativereal}. By Proposition \ref{prop:negativereal}, $r_2$ is  an upper bound among all deterministic and randomized algorithms for our setting where the expected rewards are unknown.  
In light of the lower bound in Proposition \ref{prop:negativereal}, next we revisit Example \ref{example:CR} where we evaluate Algorithm \ref{alg:ptas} and the algorithm of \cite{ball2009toward} (called the BQ algorithm) using the notion of realized CR.

\begin{figure}[ht]
   \centering
  \includegraphics[width=0.7\linewidth]{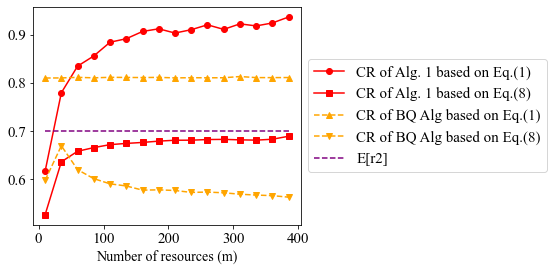}  
  \caption{The  CR and the realized CR of Algorithm \ref{alg:ptas} and the BQ algorithm, where the CR is computed using Equation \eqref{eq:CR_early} and the realized CR is computed using Equation \eqref{eq:CRdefreal}. Here, $p =0.3$, $r_1=0.9$, and  $r_2$ is uniformly drawn from the interval $(0.5,0.9)$.} 
  \label{fig:newnotion}
\end{figure}

\begin{example}[Revisiting Example \ref{example:CR}]
    Consider the same setting in Example \ref{example:CR}  with the sampling probability $p =0.3$,  $m\in\{10,35,\ldots,360,385\}$,  $r_1 =0.9$, and   $r_2\sim \text{Uniform}(0.5, 0.9)$. For this setting, we compute four quantities: (i) the  CR of Algorithm \ref{alg:ptas} per Equation \ref{eq:CR_early}, (ii)  
    the realized CR of Algorithm \ref{alg:ptas} per Equation \eqref{eq:CRdefreal}, (iii) the CR of  the BQ algorithm (i.e., $1/(2-r_2/r_1) = 1/(2-\alpha)$), (iv) the realized CR of the BQ algorithm. (Recall that (i) and (iii) are computed in Example \ref{example:CR} and here  they are just used for comparison.) Note that
    the BQ algorithm  is evaluated in a setting where the expected rewards are known. 
   To compute these quantities (i.e., quantities  (i), (ii), and (iv)), we considered all the possible values for $h$ and $\ell$, any arrival order, and the realization of reward. 

Figure \ref{fig:newnotion} shows the expected value of these four quantities versus $m$ with $p= 0.3$. Here, the expectation is with respect to $r_2$.  
The figure shows as expected, the realized CR of Algorithm \ref{alg:ptas} and the BQ algorithm  is much smaller than their CR per Equation \eqref{eq:CR_early}. As another important observation, 
 the realized CR of Algorithm \ref{alg:ptas} goes to $\E[r_2]$, which is the upper bound shown in Proposition \ref{prop:negativereal}. In addition, in terms of realized CR,  Algorithm \ref{alg:ptas} outperforms the BQ algorithm  for $m \geq 60$, which is consistent with what we observed in Example \ref{example:CR}. 
\end{example}

\subsection{Heterogeneous Sampling Probability} \label{subsec:heterop}
In our setting, we assumed that any of the $h+\ell$ agents reacts to the outreach program with probability $p$. In practice, however, different types of agents may react to the outreach program differently. To capture that, here we study a setting where type $i \in \{1,2\}$ agents react to the outreach program (i.e., falls into the sample set) with probability $p_i$.

For this setting, let us consider a slightly modified version of Algorithm \ref{alg:ptas}. 
In this modified algorithm, once $\widehat r_1> \widehat r_2$ ($\widehat r_1\le \widehat r_2$ respectively), the algorithm assigns a   protection level of  $\min\{m,s_1\frac{1-p_1}{p_1} \}$ ($\min\{m,s_2\frac{1-p_2}{p_2} \}$ respectively) to type $1$ (type $2$ respectively). The following theorem presents a lower bound on the CR of this modified algorithm.

\begin{proposition}[Heterogeneous $p$] \label{thm:heterop}
Consider the model presented in Section \ref{sec:model}, where the expected rewards of $r_i$, $i\in \{1, 2\}$ is unknown to the decision-maker. Let the sampling probability of type $1$, $2$ agents be $p_1$, $p_2$ respectively. Then, the asymptotic CR of the modified algorithm, defined above, equals to 
\[
\text{CR}=1-\Theta\left(\max\{\frac{1}{p_1\sqrt{m}}, \frac{1}{\min\{p_1,p_2\}^2m} \}\right).
\]
\end{proposition}

 The proof can be found in Section \ref{app:proofheterop}.
 From the expression of the CR, we can find that if $p_1 < p_2$, then $\Theta(\frac{1}{p_1\sqrt{m}})$ is the dominant term if and only if $p_1=\omega(1/\sqrt{m})$. That is, in this case, the asymptotic CR is $1-\Theta(\frac{1}{p_1\sqrt{m}})$. If $p_1 \geq p_2$, then $\Theta(\frac{1}{p_1\sqrt{m}})$ is the dominant term if and only if $\frac{p_1}{p_2^2}=o(\sqrt{m})$. 

\subsection{Beyond Two Types of Agents} \label{subsec:nagent}
So far, we have assumed that there are only two types of agents,  and under this assumption, we designed an asymptotically optimal algorithm (Algorithm \ref{alg:ptas}). Here, we show how to generalize Algorithm \ref{alg:ptas} to a setting with $k\ge 3$ type. 
In the generalized algorithm, similar to Algorithm \ref{alg:ptas}, we first estimate the average reward of each type of agents using the sample information. We then decide about an order over the agents using the estimated rewards, and follow 
a nesting protection level policy given the order; see Algorithm \ref{alg:ptasn}. 
{
{Consider an order $(j_1, j_2, \ldots, j_k)$  over types and protection levels $x_1, x_2, \ldots, x_k $. Given this order and  protection levels $x_1, x_2, \ldots, x_k $, the nested protection policy works as follows. It accepts the arriving agent of type $j_i$ if and only if (i) there is  resource left (ii) the total number of accepted agents of type $j_q$ for $q = i, i+1, \ldots, k$ is less than $m-x_{i-1}$. Here, we set $x_0 =0$.} }

\begin{algorithm}[htbp] 
\footnotesize
\caption{Online Resource Allocation Algorithm for $k$ types of Agents}
\begin{itemize}
\item[] \textbf{Input:} The number of resources $m$ and sample information $\boldsymbol \psi =(\mathbf s, \boldsymbol \rho)$, where $\mathbf s =(s_1, s_2,\ldots,s_k)$, $\boldsymbol \rho = (\boldsymbol \rho_1, \boldsymbol \rho_2, \ldots, \boldsymbol \rho_k)$, and for any $i\in [k]$,   $\boldsymbol \rho_i =\{\rho_{i,1}, \rho_{i,2}, \ldots, \rho_{i, s_i}\}$.

\end{itemize}
\begin{enumerate}
  \item If the number of samples for type $i\in [k]$, i.e., $s_i$, is positive, define $\widehat r_i = \frac{\sum_{j=1}^{s_i} \rho_{i,j}}{s_i}$ as an  estimate of the expected reward of type $i\in [k]$. Otherwise,  $\widehat r_i\sim \text{uniform}(0,1)$. 
  \item Sort $\widehat r_i$ for $i\in [k]$ from the largest to the smallest, i.e., $\widehat r_{j_1} \geq \widehat r_{j_2} \ldots \geq \widehat{r_{j_k}}$. 
  \item For $i \in [k-1]$, we give $x_i=\min\{m,\frac{1-p}{p}\sum_{q=1}^{i}s_{j_q} \}$ protection level to type $j_1, j_2, \ldots, j_i$ agents and run a nested protection level policy with the order  $j_1, j_2, \ldots, j_k$ and  protection levels $x_i$, $i\in [k]$.  
  
\end{enumerate}

\label{alg:ptasn}
\end{algorithm}

We believe that Algorithm \ref{alg:ptasn} is also asymptotically optimal although we find it hard to characterize a lower bound on its CR due to a large number of cases/regions that one needs to consider in its analysis.  To demonstrate its good performance though, in the following we present an example in which we evaluate Algorithm \ref{alg:ptasn} numerically.

\begin{example} \label{example:CRn}
Consider a setting in which $r_1=0.9$, $r_2\sim \text{Uniform}(0.5, 0.9)$, and $r_3\sim \text{Uniform}(0.1, 0.5)$. 
Figure \ref{fig:subN-first} shows the average of (i) the true CR of Algorithm \ref{alg:ptasn}, and (ii) $1/(3-r_3/r_2-r_2/r_1)$ versus $m$ when $p =0.3$. Here, the average is taken w.r.t. the randomness in $r_2$ and $r_3$ in our $300$ randomly generated instances. Note that $1/(3-r_3/r_2-r_2/r_1)$ is the optimal CR in \cite{ball2009toward}
 for a setting with three types and known rewards. We use $1/(3-r_3/r_2-r_2/r_1)$ as a benchmark to measure the value of the sample information. To compute the true CR, we take all $h, \ell \geq 0$, and for each $(h,\ell)$, we simulate $1000$ instances with random order (the first instance is fixed to be the ordered instance). By taking the minimum CR among all sequences described above, we compute the true CR. 
 
 Figure \ref{fig:subN-first} shows that  by increasing $m$, the CR of Algorithm \ref{alg:ptasn} improves. In addition, the true CR of the algorithm  exceeds the benchmark when $m\ge 60$.  
This highlights that even in the $3$-type case, for large enough $m$, the value of the sample information outweighs the drawbacks of not knowing the true rewards. 

Figure \ref{fig:subN-second} shows the average of the true CR of Algorithm \ref{alg:ptasn} and the benchmark  versus $p$ when $m =200$. Interestingly, we observe that  when $m=200$, even if the sampling probability is as low as $0.1$, the true CR is much larger  than the benchmark.

\begin{figure}[ht] 
  \centering
  \includegraphics[width=0.47\linewidth]{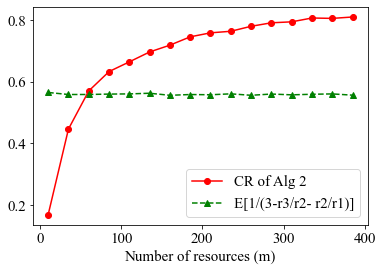}  
  \caption{ $p =0.3$}
  \label{fig:subN-first}
  \centering
  \includegraphics[width=0.47\linewidth]{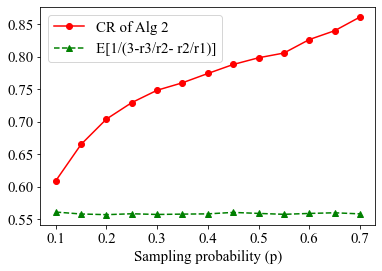}  
  \caption{ $m=200$}
  \label{fig:subN-second}
\caption{The CR of Algorithm \ref{alg:ptasn}, and the benchmark in \cite{ball2009toward}. Here, $r_1=0.9$, $r_2$ is uniformly drawn from the interval $(0.5,0.9)$, and $r_3$ is uniformly drawn from the interval $(0.1,0.3)$.}
\label{fig:figN}
\end{figure}

\end{example}

\section{Concluding Remarks  and Future Directions}
In this paper, we consider an online resource allocation problem in which the decision-maker has uncertainty about the arrival process, as well as, the obtained rewards upon the allocation. The decision-maker has access to the sample information that is often acquired through an initial test period. We study how to optimally exploit the sample information that provides partial knowledge about the arrival process and rewards. We propose a protection-level algorithm that achieves the competitive ratio of $1-\Theta(1/(p\sqrt{m}))$  and show that the obtained competitive ratio is asymptotically tight in terms of both the initial number of resources $m$ and the sampling probability $p$.

Our result shows that the sample information can be significantly beneficial, and hence opens up several new directions for future research.   One  natural  direction is to study how to optimize the number of resources $m$ upon receiving the sample information. While in many settings, due to a long lead production time, it is not possible to react to the sample information to optimize the inventory decisions, for some other settings when products/resources are produced domestically, the decision-makers can better plan their inventory decisions using the sample information.   
Another research question of interest is to
study how to generalize our setting 
to the network revenue management problem where there are multiple types of agents and multiple types of resources. The network revenue management problem under adversarial arrival settings does not admit a constant competitive ratio. Thus, it is interesting to explore if one can even obtain a  constant competitive ratio under adversarial arrival setting in the presence of sample information.

\bibliographystyle{ACM-Reference-Format}
\bibliography{sample-bibliography}

\newpage

\begin{APPENDICES}
\section{Proof of Proposition \ref{prop:negativeresults}} \label{sec:proof:negative}
We show the results for the three scenarios respectively. Recall that in scenario 1, $p=0$; that is, no sample information is available.  For this scenario,    even if the decision-maker knows the rewards $r_i$ ($i\in [2]$), without any sample information, as shown in  \cite{ball2009toward}, 
$1/(2-\alpha)$ is an upper bound for any algorithm.

Now consider  scenario 2,  where $p$ is completely unknown.  Suppose that the decision-maker has access to some sample information $(s_1,s_2)$, which means that there are $s_1$ type $1$ agents and $s_2$ type $2$ agents in the sample. As $p$ is completely unknown, $p$  can be any real number between $0$ and $1$. If $p \to 0$, the number of arriving agents of type $1$ $n_1 \to \infty$ and $n_2 \to \infty$. If $p =1$, we have $n_1=n_2=0$. Therefore, the decision-maker has no demand information because $n_i$ can be any non-negative real number for $i \in \{1,2\}$. Thus, again by \cite{ball2009toward}, the upper bound of CR among all deterministic  and randomized  algorithms is $1/(2-\alpha)$.

    \section{Proof of Statements in Section \ref{sec:CR_proof}}
   \label{app:proof_thm_2} 
   
 \subsection{Proof of Lemma \ref{lm:order}} 

    We will show that for any realization of $(n_1, n_2)$ and sample information $\boldsymbol \psi$, we have ${\rew_{A}\left(I_{\text{order}}, \boldsymbol \psi\right)}\le  {\rew_{A}\left(I, \boldsymbol \psi\right)}$. This implies that $\frac{\rew_{A}\left(I_{\text{order}}, \boldsymbol \psi\right)}{\textsc{opt}(I_{\text{order}})}\le  \frac{\rew_{A}\left(I, \boldsymbol \psi\right)}{\textsc{opt}(I)}$, as $\textsc{opt}(I) = \textsc{opt}(I_{\text{order}})$.  The last inequality is the desired result.    
    
    It remains to show that ${\rew_{A}\left(I_{\text{order}}, \boldsymbol \psi\right)}\le  {\rew_{A}\left(I, \boldsymbol \psi\right)}$. First consider a case where under  $\boldsymbol \psi$, $\widehat r_1> \widehat r_2$. In this case, Algorithm  \ref{alg:ptas} assigns a protection level of $x_1 = s_1(1-p)/p$ to type $1$ agents. Then, it is clear that (i) the number of accepted  type $2$  agents under $I_{\text{order}}$ is  larger  than or equal to that under $I$, and (ii) the number of accepted type $1$  agents under $I_{\text{order}}$ is  smaller   than or equal to that under $I$. Hence in this case,  Algorithm  \ref{alg:ptas} under $I_{\text{order}}$ obtains a lower reward than $I$.  On the other hand, when  $\widehat r_1\le  \widehat r_2$, Algorithm  \ref{alg:ptas} assigns a protection level of $x_2 = s_2(1-p)/p$ to type $2$ agents. Again, in this case, we have (i) the number of accepted type $2$  agents under $I_{\text{order}}$ is  larger  than or equal to that under $I$, and (ii) the number of accepted type $1$  agents under $I_{\text{order}}$ is  smaller   than or equal to that under $I$. Hence,  we have ${\rew_{A}\left(I_{\text{order}}, \boldsymbol \psi\right)}\le  {\rew_{A}\left(I, \boldsymbol \psi\right)}$.

\subsection{Proof of Lemma \ref{lem:c1}}
By Equation \eqref{eq:CR},  $\inf_{(h, \ell)\in \mathcal R_1}~\E_{\boldsymbol \psi}
    \left[\frac{\rewGood}{\textsc{opt}(\mathbf n)}\right]$ can be lower bounded as follows:
    \begin{align}\nonumber 
       \inf_{(h, \ell)\in \mathcal R_1}~\E_{\boldsymbol \psi}
    \left[\frac{\E[\rew_A\left(\mathbf n, \boldsymbol \rho\right)]}{\textsc{opt}(\mathbf n)}\right] \ge \inf_{(h, \ell)\in \mathcal R_1} \Bigg\{ \min\Bigg\{ \E_{\mathbf n}\left[\frac{\rewGood }{\textsc{opt}(\mathbf n)}\right], \\ \label{eq:cr1cond}
     \E_{\mathbf n}\left[\frac{\rewBad }{\textsc{opt}(\mathbf n)}\right]
   \Bigg\}\Bigg\}\,.
    \end{align}
In light of Equation \eqref{eq:cr1cond}, we divide the rest of the proof into two parts where 
 in the first part, for any realization of $\mathbf s$, we bound $\frac{\rewGood }{\textsc{opt}(\mathbf n)}$ and in the second part, we bound $\frac{\rewBad }{\textsc{opt}(\mathbf n)}$.

\textit{Part 1: bounding $\frac{\rewGood }{\textsc{opt}(\mathbf n)}$.} Recall that for a fixed realization of $\mathbf s$ (or equivalently $\mathbf n$),  $\rewGood$ is the reward of our algorithm when there are $n_i$ type $i\in\{1,2\}$ agents in the online arrival sequence  given that  $\widehat r_1 > \widehat r_2$. (See the definition of the good event $G(\boldsymbol \psi; (h, \ell))$ in Equation \eqref{eq:good_event}). Conditioned on $\widehat r_1 > \widehat r_2$, 
 Algorithm \ref{alg:ptas} assigns a protection level of $x_1=\min\{m, \frac{1-p}{p}s_1\}$ to type $1$ agents, and hence we have 
\begin{align} \notag
\rewGood&=\min\left\{n_2,\left(m-\frac{1-p}{p}s_1\right)^{+}\right\}\cdot r_2\\
&+\min\left\{m-\min\left\{n_2,\left(m-\frac{1-p}{p}s_1\right)^{+}\right\},n_1 \right\}\cdot r_1\,,\label{eq:good_reward}
\end{align}
where $\min\left\{n_2,\left(m-\frac{1-p}{p}s_1\right)^{+}\right\}$ is the number of accepted agents of type $2$ and $\min\left\{m-\min\left\{n_2,\left(m-\frac{1-p}{p}s_1\right)^{+}\right\},n_1 \right\}$ is the number of accepted agents of type $1$.  In what follows, we bound 
$\frac{\rewGood }{\textsc{opt}(\mathbf n)}$ by considering the following cases:

\begin{itemize}
        \item \textit{Case 1: $n_2 < m-n_1$.} In this case, under the  optimal clairvoyant algorithm, all agents are accepted as  we have $\min\{(m-n_1)^+,n_2\}=n_2$, and hence
        $
        \textsc{opt}(\mathbf{n})=n_2r_2+n_1r_1.
        $
        In addition, 
          the number of accepted type $1$ agents under Algorithm \ref{alg:ptas}, i.e., $\min\left\{m-\min\left\{n_2,\left(m-\frac{1-p}{p}s_1\right)^{+}\right\},n_1 \right\}$, is $n_1$; see Equation \eqref{eq:good_reward}. To see why note that  $\min\left\{n_2,\left(m-\frac{1-p}{p}s_1\right)^{+}\right\} \leq n_2$, and we have $n_1<m-n_2 \leq m-\min\left\{n_2,\left(m-\frac{1-p}{p}s_1\right)^{+}\right\}$, which implies that $\min\left\{m-\min\left\{n_2,\left(m-\frac{1-p}{p}s_1\right)^{+}\right\},n_1 \right\}=n_1$.
 Then, by Equation \eqref{eq:good_reward}, we have
              \begin{align*}
            \frac{\rewGood}{\textsc{opt}(\mathbf{n})}&=\frac{\min\left\{n_2,\left(m-\frac{1-p}{p}s_1\right)^{+}\right\}\cdot r_2+n_1 r_1}{n_2r_2+n_1r_1}\\&\geq \frac{\min\left\{n_2,\left(m-\frac{1-p}{p}s_1\right)^{+}\right\}}{n_2} = \min\left\{1,\frac{\left(m-\frac{1-p}{p}s_1\right)^{+}}{n_2} \right\} \\& = \frac{\left(m-\frac{1-p}{p}s_1\right)^{+}}{n_2}
            > \frac{\left(m-\frac{1-p}{p}s_1\right)^{+}}{m-n_1} \geq \frac{\left(m-\frac{1-p}{p}h\right)^{+}}{m} \geq \left( 1-\frac{1-p}{p}\frac{1}{\mzero} \right)^{+}\,,
        \end{align*}
        where in the last  inequality holds because in region $\mathcal R_1$, $h< \sqrt{m}$.
        
        \item \textit{Case 2: $n_2 \geq m-n_1$.} In this case,  the  optimal  clairvoyant  algorithm, cannot accept all the agents in the online arrival sequence as we have $\min\{(m-n_1)^+,n_2\}=m-n_1$. In addition, since  $n_1\leq h <\mzero < m$, we have
      $
        \textsc{opt}(\mathbf{n})=n_1r_1+(m-n_1)r_2.
     $
        Next, to bound $\frac{\rewGood }{\textsc{opt}(\mathbf n)}$, we consider the following two subcases:
        \begin{itemize}
            \item \textit{Case 2.1: $\min\{m-\min\{n_2,(m-\frac{1-p}{p}s_1)^{+}\},n_1 \}=m-\min\{n_2,(m-\frac{1-p}{p}s_1)^{+}\}$.} 
            In this case, Algorithm \ref{alg:ptas} does not  accept all type $1$ agents (see Equation \eqref{eq:good_reward}): the number of accepted  type $1$ agents is $m-\min\{n_2,(m-\frac{1-p}{p}s_1)^{+}\}$. Then, by Equation \eqref{eq:good_reward}, we have 
            \begin{align*}
             \frac{\rewGood}{\textsc{opt}(\mathbf{n})}&=\frac{\min\{n_2,(m-\frac{1-p}{p}s_1)^{+}\}r_2+(m-\min\{n_2,(m-\frac{1-p}{p}s_1)^{+}\})r_1}{(m-n_1)r_2+n_1r_1} \\&\geq \frac{mr_2}{(m-h)r_2+hr_1} \geq \frac{mr_2}{(m-\mzero )r_2+\mzero r_1} = 1-\frac{\mzero r_1-\mzero r_2}{\mzero r_1+(m-\mzero )r_2}.
            \end{align*}
            The first inequality is due to $\min\{n_2,(m-\frac{1-p}{p}s_1)^{+}\}r_2+(m-\min\{n_2,(m-\frac{1-p}{p}s_1)^{+}\})r_1 \geq \min\{n_2,(m-\frac{1-p}{p}s_1)^{+}\}r_2+(m-\min\{n_2,(m-\frac{1-p}{p}s_1)^{+}\})r_2 = mr_2$.

            \item \textit{Case 2.2:  $\min\{\frac{1-p}{p}s_1,n_1 \}=n_1$.} In this case,  Algorithm \ref{alg:ptas}  accepts all type $1$ agents. We further   have $n_2 \geq m-n_1 \geq (m-\frac{1-p}{p}s_1)^{+}$. Then, by Equation \eqref{eq:good_reward}, we have
            \begin{align*}
             \frac{\rewGood}{\textsc{opt}(\mathbf{n})}&=\frac{\min\{n_2,(m-\frac{1-p}{p}s_1)^{+}\}r_2+n_1r_1}{(m-n_1)r_2+n_1r_1} \\&=\frac{(m-\frac{1-p}{p}s_1)^{+}r_2+n_1r_1}{(m-n_1)r_2+n_1r_1} \geq \frac{(m-\frac{1-p}{p}s_1)^{+}}{m-n_1} \geq \frac{(m-\frac{1-p}{p}h)^{+}}{m} \geq \left(1-\frac{1-p}{p}\frac{1}{\mzero}\right)^{+}.
            \end{align*}
          \end{itemize}
    \end{itemize}
    Putting these cases together, we have
    \begin{align}\frac{\rewGood}{\textsc{opt}(\mathbf{n})}\ge \min\left\{\left(1-\frac{1-p}{p}\frac{1}{\mzero}\right)^{+},1-\frac{\mzero r_1-\mzero r_2}{\mzero r_1+(m-\mzero )r_2}\right\}\,. \label{eq:cr1_part1}\end{align}
    
    \textit{Part 2: bounding $\frac{\rewBad }{\textsc{opt}(\mathbf n)}$.} Recall that for fixed realization of $\mathbf s$,  $\rewBad$ is the reward of our algorithm when there are $n_i$ type $i\in\{1,2\}$ agents in the online arrival sequence  and   $\widehat r_1 \le  \widehat r_2$. In this case, the algorithm assigns a protection level of $x_2 = \frac{1-p}{p} s_2$, and hence we have
\begin{align}
\rewBad\geq \min\{m,n_2\}\cdot r_2+\min\{n_1,(m-n_2)^+\}\cdot r_1. \label{eq:bad_event_reward}
\end{align}

   To bound our desired ratio, we consider the following cases:
   
    \begin{itemize}
        \item \textit{Case 1: $n_2 \leq m-n_1$.} In this case, we have $\min\{m,n_2\}=n_2$, and $\min\{n_1,(m-n_2)^+\}=n_1$, and hence,  by Equation \eqref{eq:bad_event_reward}, $\rewBad\ge n_1 \cdot r_1 + n_2 \cdot r_2$. 
    Thus, we have
        \[
        \frac{\rewBad}{\textsc{opt}(\mathbf{n})} \ge\frac{n_1r_1+n_2r_2}{n_1r_1+n_2r_2}=1.
        \]
        
        \item \textit{Case 2: $n_2>m-n_1$.} In this case, we have $\min\{n_1,(m-n_2)^+\}=(m-n_2)^+$, $\min\{(m-n_1)^+,n_2\}=m-n_1$. Then, by Equation \eqref{eq:bad_event_reward}, we have 
        \begin{align*}
        \frac{\rewBad}{\textsc{opt}(\mathbf{n})} &\ge \frac{\min\{m,n_2\}r_2+(m-n_2)^{+}r_1}{n_1r_1+(m-n_1)r_2} 
        \\
        &\geq \frac{mr_2}{\mzero r_1+(m-\mzero )r_2}= 1-\frac{\mzero r_1-\mzero r_2}{\mzero r_1+(m-\mzero )r_2}.
        \end{align*}

    \end{itemize}
    Considering the following two cases, we have 
      \begin{align}\frac{\rewBad}{\textsc{opt}(\mathbf{n})}\ge 1-\frac{\mzero r_1-\mzero r_2}{\mzero r_1+(m-\mzero )r_2}\,. \label{eq:cr1_part2}\end{align}
    By plugging in  Equations \eqref{eq:cr1_part1} and \eqref{eq:cr1_part2} into 
     Equation \eqref{eq:cr1cond}, we get the desired result.

\subsection{Proof of Lemma \ref{lem:c3}}
We follow similar steps in the proof of Lemma \ref{lem:c1}. That is, we bound 
    \begin{align}\nonumber 
       \inf_{(h, \ell)\in \mathcal R_2}~\E_{\mathbf{s}}
    \left[\frac{\E[\rew_A\left(\mathbf n, \boldsymbol \rho\right)]}{\textsc{opt}(\mathbf n)}\right] \ge \inf_{(h, \ell)\in \mathcal R_2} \Bigg\{ \min\Bigg\{ \E_{\mathbf n}\left[\frac{\rewGood }{\textsc{opt}(\mathbf n)}\right], \\ \label{eq:cr2cond}
     \E_{\mathbf n}\left[\frac{\rewBad }{\textsc{opt}(\mathbf n)}\right]
   \Bigg\}\Bigg\}\,.
    \end{align}
In light of Equation \eqref{eq:cr2cond}, we divide the rest of the proof into two parts where 
 in the first part, for any realization of $\mathbf s$, we bound $\frac{\rewBad}{\textsc{opt}(\mathbf n)}$ and in the second part,  we bound $\E_{\mathbf n}\big[\frac{\rewGood }{\textsc{opt}(\mathbf n)}\big]$. While bounding the former ratio is straightforward, bounding the latter one is quite involved. This is because when the good event $\mathcal G(\boldsymbol \psi; (h, \ell))$ happens, Algorithm \ref{alg:ptas} may still over-protect or under-protect type $1$ agents given the randomness in the sample information. As $h$ is large in this  case, the number of type $1$ agents that the algorithm  over-protects or under-protects can be  large too, leading to a low reward for the algorithm  under some realization of $\mathbf s$. Nonetheless, we show that the expected reward of the algorithm, in this case, is still large once we take expectation with respect to $\mathbf s$.

\textit{Part 1: bounding $\frac{\rewBad }{\textsc{opt}(\mathbf n)}$. } Recall that for fixed realization of $\mathbf s$,  $\rewBad$ is the reward of our algorithm when there are $n_i$ type $i\in\{1,2\}$ agents in the online arrival sequence  and   $\widehat r_1 \le  \widehat r_2$. In this case, the algorithm assigns a protection level of $x_2 = \min\{m,\frac{1-p}{p} s_2\}$, and hence we have
     \[
    \rewBad \geq \min\{m,n_2\}r_2+\min\{n_1,(m-n_2)^+\}r_1.
    \]
    As $\ell <\mzero $, we have $n_2 \leq \ell <\mzero  <m$. Hence, 
    \begin{align*}
              \frac{\rewBad}{\textsc{opt}(\mathbf n)} &\ge \frac{\rewBad}{n_1r_1+\min\{m-n_1,n_2\}r_2}\\
              &=\frac{n_2r_2+\min\{m-n_2,n_1\}r_1}{n_1r_1+\min\{m-n_1,n_2\}r_2} \\&\ge  \min\left\{\frac{n_1r_1+n_2r_2}{n_1r_1+n_2r_2},\frac{(m-n_2)r_1+n_2r_2}{n_1r_1+(m-n_1)r_2} \right\} \\&\overset{a}{=} \frac{(m-n_2)r_1+n_2r_2}{n_1r_1+(m-n_1)r_2} \geq \frac{(m-\mzero )r_1}{mr_1} =1-\frac{1}{\mzero}.
            \end{align*}
            The second inequality is because when $\min\{m-n_2,n_1\}=n_1$, we have $n_2<m-n_1$ and $\min\{m-n_1,n_2\}=n_2$, and when $\min\{m-n_2,n_1\}=m-n_2$, we have $m-n_1<n_2$ $\min\{m-n_1,n_2\}=m-n_1$.

\textit{Part 2: bounding $\inf_{(h,\ell) \in \mathcal R_2}\E_{\mathbf n}\Big[\frac{\rewGood }{\textsc{opt}(\mathbf n)}\Big]$. } We begin by partitioning region $\mathcal R_2$ based on the number of type $1$ agents, (i.e. $h$):
\[\overline{\mathcal{R}_2} =\{(h, \ell): (h, \ell)\in \mathcal R_2, h\ge \hzero  \} \qquad \underline{\mathcal{R}_2} =\{(h, \ell): (h, \ell)\in \mathcal R_2, h< \hzero  \}\,, \]
where  $\hzero =\hzeroExp$. Then, $\inf_{(h,\ell) \in \mathcal R_2}\E_{\mathbf n}\big[\frac{\rewGood }{\textsc{opt}(\mathbf n)}\big]$ is equal to the minimum of $\inf_{(h,\ell) \in \underline{\mathcal{R}_2}}\E_{\mathbf n}\big[\frac{\rewGood }{\textsc{opt}(\mathbf n)}\big]$ and $\inf_{(h,\ell) \in \overline{\mathcal{R}_2}}\E_{\mathbf n}\big[\frac{\rewGood }{\textsc{opt}(\mathbf n)}\big]$. We bound the first ratio (i.e., the one concerns $\underline{\mathcal{R}_2}$) in Lemma \ref{lem:cr2case2} and the second ratio  (i.e., the one concerns $\overline{\mathcal{R}_2}$) will be bounded in Lemma \ref{lem:cr2case3}. In addition, $\hzero>\sqrt{m}$ is guaranteed because $\hzero={\frac{m-\mzero }{1-p}+\frac{1}{2p^2}+\sqrt{\frac{1}{4p^4}+\frac{(m-\mzero )}{p^2(1-p)}}}>m-\sqrt{m}+\sqrt{2(m-\sqrt{m})}>\sqrt{m}$ for all $m>0$.

    \begin{lemma} \label{lem:cr2case2}
Let $\underline{\mathcal{R}_2} =\{(h, \ell): (h, \ell)\in \mathcal R_2, h< \hzero  \}$, where $\mathcal R_2 =\{(h,\ell):  h\ge \mzero , \ell < \mzero  \}$, $\hzero =\hzeroExp$. 
Then,
\[
\inf_{(h, \ell)\in \underline{\mathcal{R}_2}}\E_{\mathbf n}\left[\frac{\rewGood }{\textsc{opt}(\mathbf n)}\right] \geq (1-\frac{1}{\sqrt{m}})(1-\frac{1}{m})\,.
\]
\end{lemma}

Lemma \ref{lem:cr2case2} gives the lower bound of the ratio conditional on the good event happens. In this case, the algorithm gives a certain protection level to type $1$ agents. Then, we show that when $(h, \ell)\in \underline{\mathcal{R}_2}$, with high probability, the number of type $2$ agents is less than the total number of resources minus the protection level. Thus, in this regime, we do not reject any type $2$ agents with high probability. Conditional on this event, we can  obtain a lower bound of the ratio for any realization $\mathbf{s}$.

Next, we lower bound $\inf_{(h,\ell) \in \overline{\mathcal{R}_2}}\E_{\mathbf n}\Big[\frac{\rewGood }{\textsc{opt}(\mathbf n)}\Big]$ in Lemma \ref{lem:cr2case3}.

\begin{lemma} \label{lem:cr2case3}
Let $\overline{\mathcal{R}_2} =\{(h, \ell): (h, \ell)\in \mathcal R_2, h\ge \hzero  \}$, where $\mathcal R_2 =\{(h,\ell):  h\ge \mzero , \ell < \mzero  \}$,  where $\mOne =\mOneExp$,  $\beta=0.4215 \cdot \frac{p^2+(1-p)^2}{p(1-p)}$, $\hzero =\hzeroExp$. Then,
\[
\inf_{(h, \ell)\in \overline{\mathcal{R}_2}}
\E_{\mathbf n}\left[\frac{\rewGood }{\textsc{opt}(\mathbf n)}\right] \geq \min\left\{1-\frac{1}{\sqrt{m}}, \frac{\mOne }{m}\right\}.
\]

\end{lemma}

To show Lemma \ref{lem:cr2case3}, notice that given the good event happens, the algorithm may over- or under-protect type $1$ agents. If it under-protects type $1$ agents, we further split the analysis into two cases: the number of type $2$ agents is less than, or is larger than the total number of resources minus the protection level. In both cases, we can easily find a lower bound for any realization $\mathbf{s}$. If the algorithm over-protects type $1$ agents, then it accepts all arriving type $1$ agents, but some resources are wasted due to over-protecting. In this case, we lower bound the ratio by a linear function with respect to the number of arriving type $1$ agents, $n_1$. Then, we used the results  in \cite{nagaev2011bound} and \cite{berend2013sharp} to estimate the conditional expectation of $n_1$ given that the algorithm over-protects type $1$ agents.

Finally, by Lemmas \ref{lem:cr2case2} and  \ref{lem:cr2case3}, we have the following inequality, which is the desired result.
\[
\inf_{(h, \ell)\in \mathcal R_2}~\E_{\mathbf{s}}
    \left[\frac{\rewGood}{\textsc{opt}(\mathbf n)}\right]\ge \CRTwo   =
    \CRTwoExp\,.
\]

\subsubsection{Proof of Lemma \ref{lem:cr2case2}}
Define event $\eventzero$ as $\{n_2 \leq m-\frac{1-p}{p}s_1\}$. We provide a lower bound on $\E_{\mathbf n}\left[\frac{\rewGood }{\textsc{opt}(\mathbf n)}\right]$ as follows:
\begin{align}\nonumber
    &\inf_{(h,\ell) \in \underline{\mathcal{R}_2}}\E_{\mathbf n}\left[\frac{\rewGood }{\textsc{opt}(\mathbf n)}\right] \\&\nonumber= \inf_{(h,\ell) \in \underline{\mathcal{R}_2}}\E_{\mathbf n}\left[\frac{\rewGood }{\textsc{opt}(\mathbf n)}\big \vert{\eventzeroc}\right]\Pr(\eventzeroc)+\E_{\mathbf n}\left[\frac{\rewGood }{\textsc{opt}(\mathbf n)}\big \vert{\eventzero}\right]\Pr(\eventzero) \\&\geq \label{eq:condcr2two} \inf_{(h,\ell) \in \underline{\mathcal{R}_2}}\E_{\mathbf n}\left[\frac{\rewGood }{\textsc{opt}(\mathbf n)}\big \vert{\eventzero}\right]\Pr(\eventzero)
\end{align}
In light of Equation \eqref{eq:condcr2two}, we will provide a lower bound for $\E_{\mathbf n}\left[\frac{\rewGood }{\textsc{opt}(\mathbf n)}\big \vert{\eventzero}\right]$ and $\Pr(\eventzero)$.

\textit{Part 1: bounding $\E_{\mathbf n}\left[\frac{\rewGood }{\textsc{opt}(\mathbf n)}\big \vert{\eventzero}\right]$.} For a fixed realization of $\mathbf s$ (or equivalently $\mathbf n$), we have
\begin{align} \notag
\rewGood&=\min\left\{n_2,\left(m-\frac{1-p}{p}s_1\right)^{+}\right\}\cdot r_2\\
&+\min\left\{m-\min\left\{n_2,\left(m-\frac{1-p}{p}s_1\right)^{+}\right\},n_1 \right\}\cdot r_1\,\label{eq:good_reward_2}
\end{align}
Conditional on $\eventzero=\{n_2 \leq m-\frac{1-p}{p}s_1\}$, we have $\min\{n_2,(m-\frac{1-p}{p}s_1)^{+}\}=n_2$, and $\min\{m-\min\{n_2,(m-\frac{1-p}{p}s_1)^{+}\},n_1 \}=\min\{m-n_2,n_1 \}$. Then, we have
\begin{align}\nonumber
              \frac{\rewBad}{\textsc{opt}(\mathbf n)} &\ge \frac{\rewBad}{n_1r_1+\min\{m-n_1,n_2\}r_2}\\
              &\nonumber=\frac{n_2r_2+\min\{m-n_2,n_1\}r_1}{n_1r_1+\min\{m-n_1,n_2\}r_2} \\&\nonumber \ge  \min\left\{\frac{n_1r_1+n_2r_2}{n_1r_1+n_2r_2},\frac{(m-n_2)r_1+n_2r_2}{n_1r_1+(m-n_1)r_2} \right\} \\&= \label{eq:cr2part1} \frac{(m-n_2)r_1+n_2r_2}{n_1r_1+(m-n_1)r_2} \geq \frac{(m-\mzero )r_1}{mr_1} =1-\frac{1}{\mzero}.
            \end{align}
            The second inequality is because when $\min\{m-n_2,n_1\}=n_1$, we have $\min\{m-n_1,n_2\}=n_2$, and when $\min\{m-n_2,n_1\}=m-n_2$, we have $\min\{m-n_1,n_2\}=m-n_1$.
            
            \textit{Part 2: bounding $\Pr(\eventzero)$.} We want to show that $\Pr(\eventzero)$ is a high probability event for $(h,
            \ell) \in \underline{\mathcal{R}_2}$. As $\ell\in \mzero $ for any $(h,
            \ell) \in \underline{\mathcal{R}_2}$, we have
            $n_2 \leq \mzero $.  Then, because $\eventzero = \{n_2 \le m-\frac{1-p}{p}s_1\}$, we have
            \[
            \Pr(\eventzero) = \Pr\left(\frac{1-p}{p}s_1 \leq m-n_2\right) \geq \Pr\left(\frac{1-p}{p}s_1 \leq m-\mzero \right) = \Pr\left(s_1 \leq \frac{p(m-\mzero )}{1-p}\right).
            \]
            Next, we show that with high probability, the upper bound of $s_1$ is no more than $\frac{p(m-\mzero )}{1-p}$.
            Because $s_1 \sim \text{Bin}(h,p)$, by Hoeffding's inequality, we have
                \[
                \Pr(|s_1-hp|\leq \sqrt{h}) \geq 1-\frac{1}{h^2} \geq 1-\frac{1}{m}.
                \]
                This implies that $s_1 \leq hp+ \sqrt{h}$ happens with high probability.
                  We then have \begin{align}
               1- \frac{1}{m}\le \Pr(s_1 \le hp+\sqrt{h}) \le  \Pr(s_1 \le \hzero p+\sqrt{\hzero }) \leq \Pr\left(s_1 \le \frac{p(m-\mzero )}{1-p}\right)\,, \label{eq:cr2part2}
                \end{align}
                where the second  inequality holds because $hp+ \sqrt{h}$ is an increasing function, and $h<\hzero $ in region $\underline{\mathcal R_2}$. These imply that $hp+\sqrt{p}\le \hzero  p+\sqrt{p_0}$.  The last inequality holds because $\frac{p(m-\mzero )}{1-p}\ge \hzero  p+\sqrt{\hzero }$. Recall that $\hzero =\hzeroExp$.

 By Equations \eqref{eq:cr2part1} and \eqref{eq:cr2part2}, we have 
        \begin{equation}
            \inf_{(h,\ell) \in \underline{\mathcal{R}_2}}\E_{\mathbf n}\left[\frac{\rewGood }{\textsc{opt}(\mathbf n)}\big \vert{\eventzero}\right]\Pr(\eventzero) \geq (1-\frac{1}{\sqrt{m}})(1-\frac{1}{m})\,,
        \end{equation}
        which is the desired result.

\subsubsection{Proof of Lemma \ref{lem:cr2case3}}
We define event $\eventone=\{n_1 > \frac{1-p}{p}s_1\}$. We then have
\begin{align} \nonumber
    \inf_{(h,\ell) \in \overline{\mathcal{R}_2} }\E_{\mathbf n}\left[\frac{\rewGood }{\textsc{opt}(\mathbf n)}\right]   
    \nonumber &\geq \inf_{(h,\ell) \in \overline{\mathcal{R}_2} } \Bigg\{ \min\Bigg\{ \E_{\mathbf n}\left[\frac{\rewGood }{\textsc{opt}(\mathbf n)} \big\vert{\eventone}\right], \\ \label{eq:condcr2second}
     &\E_{\mathbf n}\left[\frac{\rewGood }{\textsc{opt}(\mathbf n)}\big\vert{\eventone^{\mathsf{C}}}\right]
   \Bigg\}\Bigg\}\,
\end{align}

We will provide the lower bound of $\inf_{(h,\ell) \in \overline{\mathcal{R}_2} } \E_{\mathbf n}\big[\frac{\rewGood }{\textsc{opt}(\mathbf n)}\big\vert{\eventone} \big]$ in the first part  of the proof, and the lower bound of $\inf_{(h,\ell) \in \overline{\mathcal{R}_2} } \E_{\mathbf n}\big[\frac{\rewGood }{\textsc{opt}(\mathbf n)}\big\vert{\eventone^{\mathsf{C}}} \big]$ in the second part  of the proof.

\textit{Part 1: bounding $\inf_{(h,\ell) \in \overline{\mathcal{R}_2} } \E_{\mathbf n}\big[\frac{\rewGood }{\textsc{opt}(\mathbf n)}\big\vert{\eventone} \big]$.} We then  have
\begin{align*}
    &\inf_{(h,\ell) \in \overline{\mathcal{R}_2} } \E_{\mathbf n}\big[\frac{\rewGood }{\textsc{opt}(\mathbf n)}\big\vert{\eventone} \big] \\&\geq \inf_{(h,\ell) \in \overline{\mathcal{R}_2} } \min \Bigg\{ \E_{\mathbf n}\big[\frac{\rewGood }{\textsc{opt}(\mathbf n)}\big\vert{\eventone, \eventzeroc} \big], \E_{\mathbf n}\big[\frac{\rewGood }{\textsc{opt}(\mathbf n)}\big\vert{\eventone, \eventzero} \big] \Bigg\}\,,
\end{align*}
where we recall $\eventzeroc=\{n_2>m-\frac{1-p}{p}s_1\}$.

To provide a lower bound for $\E_{\mathbf n}\big[\frac{\rewGood }{\textsc{opt}(\mathbf n)}\big\vert{\eventone, \eventzeroc} \big]$, we note that under event 
$\eventzeroc$, we have  $n_2>m-\frac{1-p}{p}s_1$, and hence  
            \begin{equation} 
            \E_{\mathbf n}\left[\frac{\rewGood }{\textsc{opt}(\mathbf n)}\big\vert{\eventone, \eventzeroc} \right]\geq \frac{(m-\frac{1-p}{p}s_1)^{+}r_2+\min\{m,\frac{1-p}{p}s_1\}r_1}{mr_1} \geq \frac{(m-n_2) r_1}{mr_1} \geq 1-\frac{1}{\mzero}.
            \end{equation}
            
Then, we give a lower bound to $\E_{\mathbf n}\big[\frac{\rewGood }{\textsc{opt}(\mathbf n)}\big\vert{\eventone, \eventzero} \big]$. Conditional on $\eventone$ and  $\eventzero$, we have $\frac{1-p}{p}s_1 \geq m-n_2$ and $\frac{1-p}{p}s_1 \leq n_1$. That is, the number of type $1$ agents in the online arrival sequence is greater than or equal to the protection level of Algorithm \ref{alg:ptas} for type $1$ agents (i.e., $\frac{1-p}{p}s_1 \leq n_1$), and in addition, the algorithm ends up rejecting some of type $2$ agents as $\frac{1-p}{p}s_1 \geq m-n_2$. Therefore, 
we have
\[
\E_{\mathbf n}\left[\frac{\rewGood }{\textsc{opt}(\mathbf n)}\big\vert{\eventone, \eventzero} \right] \geq \frac{(m-\frac{1-p}{p}s_1)^{+}r_2+\min\{m,\frac{1-p}{p}s_1\}r_1}{mr_1} \geq \frac{(m-n_2)r_1}{mr_1} \geq 1-\frac{1}{\mzero}\,,
\]
where the last inequality holds because $n_2\le \ell \le \mzero $ for any $(h,\ell)\in \underline{\mathcal R_2}$.

\textit{Part 2: bounding $\inf_{(h,\ell) \in \overline{\mathcal{R}_2} } \E_{\mathbf n}\big[\frac{\rewGood }{\textsc{opt}(\mathbf n)}\big\vert{\eventone^{\mathsf{C}}} \big]$. } Recall that for a fixed realization of $\mathbf s$ (or equivalently $\mathbf n$), 
\begin{align*} \notag
\rewGood&=\min\left\{n_2,\left(m-\frac{1-p}{p}s_1\right)^{+}\right\}\cdot r_2\\
&+\min\left\{m-\min\left\{n_2,\left(m-\frac{1-p}{p}s_1 \right)^{+}\right\},n_1 \right\}\cdot r_1\,.
\end{align*}

Conditional on $\eventone^{\mathsf{C}} = \{n_1 \leq \frac{1-p}{p}s_1 \}$, we do not receive enough number of type $1$ agents and thus, we have
\[
\min\left\{m-\min\left\{n_2,\left(m-\frac{1-p}{p}s_1 \right)^{+}\right\},n_1 \right\}=\min\left\{\max\left\{m-n_2,\min\{m,\frac{1-p}{p}s_1\}\right\},n_1 \right\}=n_1\,.
\]
That is, Algorithm \ref{alg:ptas}  accepts all type $1$ agents in the online arrival sequence. Then, we have
\[
\E_{\mathbf n}\Big[\frac{\rewGood }{\textsc{opt}(\mathbf n)}\big\vert{\eventone^{\mathsf{C}}} \Big] \geq \E_{\mathbf n}\Big[\frac{\min\left\{n_2,(m-\frac{1-p}{p}s_1)^{+}\right\}\cdot r_2+n_1r_1 }{mr_1}\big\vert{\eventone^{\mathsf{C}}} \Big] \geq \frac{1}{m}\E_{\mathbf n}\big[ n_1 \big\vert{\eventone^{\mathsf{C}}} \big].
\]
We will show that $\E_{\mathbf n}\big[ n_1 \big\vert{\eventone^{\mathsf{C}}} \big] \ge \mOne $, where $\mOne  = \mOneExp$ and $\beta=0.4215 \cdot \frac{p^2+(1-p)^2}{p(1-p)}$.

 Since $\eventone^{\mathsf{C}} = \{n_1 \leq \frac{1-p}{p}s_1 \}= \{h(1-p)-n_1 \geq 0 \}$, and $h(1-p)=\E[h(1-p)]$, we have
\begin{align} \notag
    \E_{\mathbf n}\big[ n_1 \big\vert{\eventone^{\mathsf{C}}} \big] &= h(1-p)-\E_{\mathbf n}\big[ h(1-p)-n_1 \big\vert{h(1-p)-n_1 \geq 0} \big]\\ \notag
       &\ge h(1-p)-\E_{\boldsymbol {s}}
    \left[|h(1-p)-n_1|\right] -  \frac{\beta}{\sqrt{h}}\\ \notag
    &\ge h(1-p)-\sqrt{p(1-p)h}-  \frac{\beta}{\sqrt{h}}\\ \notag
      &\ge  \hzero (1-p)-\sqrt{p(1-p)\hzero }-  \frac{\beta}{\sqrt{\hzero }}
  \\ &= \notag
   \mOne \,,
\end{align} 
where the first inequality holds because by 
 \cite{nagaev2011bound}, we have
$
    \Big|\E_{\boldsymbol {s}}
    \Big[(h(1-p)-n_1) \big \vert \{h(1-p)-n_1 \geq 0\}\Big]\,-\E_{\boldsymbol {s}}
    \left[|h(1-p)-n_1|\right]\Big| \le \frac{\beta}{\sqrt{h}}
$. The second inequality follows from  \cite{berend2013sharp} that shows  
$
\E_{\boldsymbol s}\left[|h(1-p)-n_1|\right]\, \leq \sqrt{p(1-p)h}.
$ The third inequality holds because when $h>\hzero $, $h(1-p)-\sqrt{p(1-p)h}- \frac{\beta}{\sqrt{h}}$ is a non-decreasing function for any $p$. To see why note that  \begin{equation} \label{eq:finalcr2}
    \frac{\partial (h(1-p)-\sqrt{p(1-p)h}- \frac{\beta}{\sqrt{h}})}{\partial h} \geq (1-p)-\frac{1}{2}\sqrt{p(1-p)}\frac{1}{\sqrt{h}}\ge 0\,.
\end{equation}
where the last inequality holds for any  $h \geq \frac{1}{4}\frac{p}{1-p}$, and we have $\hzero  > \frac{m-\sqrt{m}}{1-p} > \frac{1}{4}\frac{p}{1-p}$ for all $m \geq 2$.

\subsection{Proof of Lemma \ref{thm:c4}}
By Equation \eqref{eq:CR}, we have
\begin{align}\nonumber \inf_{(h, \ell)\in \mathcal R_3}~\E_{\boldsymbol \psi}
    \left[\frac{\rewGood}{\textsc{opt}(\mathbf n)}\right] &=\inf_{(h, \ell)\in \mathcal R_3} \Bigg\{ \E_{\mathbf n}\left[\frac{\rewGood \cdot \Pr(\mathcal G(\boldsymbol \psi; (h, \ell)))}{\textsc{opt}(\mathbf n)}\right]\\ \label{eq:CR3} 
    &+\E_{\mathbf n}\left[\frac{\rewBad \cdot \Pr(\mathcal G^{\mathsf{C}}(\boldsymbol \psi; (h, \ell)))}{\textsc{opt}(\mathbf n)}\right]\Bigg\}
    \end{align}
We split the analysis into two cases based on the initial number of resources $m$. In the first case, we assume that $m \geq \mTwo $, and in the second case, we assume that $m<\mTwo $. Here, $\mTwo =\mTwoExp$  is  a constant that only depends on $p$, $r_1$, and $r_2$. 
Case 1 is considered in Lemma \ref{lem:cr3largem} while case 2 is  studied in Lemma \ref{lem:cr3smallm}.  
\begin{lemma}\label{lem:cr3largem}
Let $\mathcal R_3 =\{(h,\ell):  h\ge \mzero , \ell \ge \mzero  \}$. Let $\mTwo =\mTwoExp$.

When the initial number of resources  $m \geq \mTwo $, we have
\[
\inf_{(h, \ell)\in \mathcal R_3}~\E_{\mathbf n}
    \left[\frac{\rewGen}{\textsc{opt}(\mathbf n)}\right] \geq \CRThreeExpOne\,.
\]
where $\X= \min\left\{1-\frac{1}{\sqrt{m}}, \frac{\mOne }{m} \right\}$, $\mOne =\mOneExp$,  $\beta=0.4215 \cdot \frac{p^2+(1-p)^2}{p(1-p)}$,  and $\hzero =\hzeroExp$.
\end{lemma}

To show Lemma \ref{lem:cr3largem}, we observe that when the initial number of resources $m \geq m_1$, the good event happens with high probability. Then, conditional on the algorithm protects type $1$ agents, we give the lower bound of the ratio by the similar methodology introduced in Lemma \ref{lem:cr2case2} and Lemma \ref{lem:cr2case3}: we discuss whether the algorithm over- or under-protects type $1$ agents, and whether the number of type $2$ agents in the arriving sequence is less than, or is larger than the total number of resources minus the protection level. 

\medskip
\begin{lemma}\label{lem:cr3smallm}
Let $\mathcal R_3 =\{(h,\ell):  h\ge \mzero , \ell \ge \mzero  \}$, $\mTwo =\mTwoExp$, and $\alpha=\frac{r_2}{r_1}$.
When the initial number of resources  $m < \mTwo $, we have
\begin{align*}
    &\inf_{(h, \ell)\in \mathcal R_3}~\E_{\mathbf n}
    \left[\frac{\E[\rew_A\left(\mathbf n, \boldsymbol \rho\right)]}{\textsc{opt}(\mathbf n)}\right] \geq \V \cdot \min \Bigg\{(1-\frac{1}{m^2})\X, \frac{1}{2}(1-\frac{1}{m^2})\X+\frac{1}{2}\min\left\{(1-\frac{1}{\ell_0^2})\alpha, 1-\frac{1-p}{pm}\ell_1\right\} \Bigg\}\,.
\end{align*}
where $\V=1-2(1-p)^{\sqrt{m}} $, $\X= \min\left\{1-\frac{1}{\sqrt{m}}, \frac{\mOne }{m} \right\}$, $\mOne =\mOneExp$,  $\beta=0.4215 \cdot \frac{p^2+(1-p)^2}{p(1-p)}$, $\hzero =\hzeroExp$, $\ell_0=\ellzeroExp$, and $\ell_1= \sqrt{p(1-p)\ell_0} +\frac{\beta}{\sqrt{m}} +\ell_0 p$.
\end{lemma}

To show Lemma \ref{lem:cr3smallm}, we should use a different technique rather than the one we use in Lemma \ref{lem:cr3largem}. This is because if the initial number of resources $m < \mTwo$, the bad event is no longer a rare event. Then, if the good event happens, the ratio is derived in Lemma \ref{lem:cr3largem}. If the bad event happens, the algorithm gives protection level to type $2$ agents. Then, we first consider that if the algorithm under-protects tyepe $2$ agents, the algorithm accepts all type $2$ agents, and give the remaining resources to type $1$ agents. In this case, the ratio is bounded by $\alpha$. Otherwise, if the algorithm over-protects type $2$ agents, the algorithm rejects some type $1$ agents and wastes some resources. Then, we use concentration inequality to bound the loss.

\subsubsection{Proof of Lemma \ref{lem:cr3largem}}
We start by defining the  event
\[\eventtwo=\left\{s_1 \in [hp-\sqrt{h},hp+\sqrt{h}]\} \cap \{s_2 \in [\ell p-\sqrt{\ell},\ell p+\sqrt{\ell}]\right\}.
\]
Under event $\eventtwo$, both $s_1$ and $s_2$ are concentrated around their average. Recall that $s_1 \sim \text{Bin}(h, p)$ and  $s_2 \sim \text{Bin}(\ell, p)$.  By Equation \eqref{eq:CR}, we have
\begin{align}
    &\label{eq:event_s} \inf_{(h, \ell)\in \mathcal R_3}~\E_{\boldsymbol \psi}
    \left[\frac{\rewGen}{\textsc{opt}(\mathbf n)}\right] \\
    &\geq \label{eq:CR3EQ2} \inf_{(h, \ell)\in \mathcal R_3}  \E_{\mathbf n}\left[\frac{\rewGood \cdot \Pr(\mathcal G(\boldsymbol \psi; (h, \ell)))}{\textsc{opt}(\mathbf n)} \big\vert{\eventtwo}\right]\Pr(\eventtwo)\,.
\end{align}
We lower bound $\Pr(\eventtwo)$ in {Part 1} of the proof. Next, we explain how to handle the conditional expectation in Equation \eqref{eq:CR3EQ2}. 
We start by  defining the following event
\[\eventthree=\left\{\widehat{r_1} \in [r_1-\frac{1}{\sqrt{s_1}},r_1+\frac{1}{\sqrt{s_1}}  ]\} \cap \{\widehat{r_2} \in [r_2-\frac{1}{\sqrt{s_2}},r_2+\frac{1}{\sqrt{s_2}}  ]\right\}.
\]
Under 
event $\eventthree$, our estimates for $r_1$ and $r_2$ are concentrated around their true values. 
The conditional expectation in Equation \eqref{eq:CR3EQ2} is lower bounded as follows:
\begin{align}\nonumber
    &\inf_{(h, \ell)\in \mathcal R_3}  \E_{\mathbf n}\left[\frac{\rewGood \cdot \Pr(\mathcal G(\boldsymbol \psi; (h, \ell)))}{\textsc{opt}(\mathbf n)} \big\vert{\eventtwo}\right] 
    \\&\geq \label{eq:CR3EQ3} \inf_{(h, \ell)\in \mathcal R_3} \E_{\mathbf n}\left[\frac{\rewGood \cdot \Pr(\mathcal G(\boldsymbol \psi; (h, \ell)))}{\textsc{opt}(\mathbf n)} \big\vert{\eventtwo, \eventthree}\right]\Pr(\eventthree \vert \eventtwo)
\end{align}
Given \eqref{eq:CR3EQ3}, we lower bound $\Pr(\eventthree \vert \eventtwo)$ in part 2 of the proof. In part 3,  conditioned on $\eventtwo, \eventthree$, we show that  $\Pr(\mathcal G(\boldsymbol \psi; (h, \ell)))$, for any realization $\mathbf{s}$, equals to $1$ in {part 3}. Finally, in part 4, we we lower bound 
$\inf_{(h, \ell)\in \mathcal R_3} \E_{\mathbf n}\left[\frac{\rewGood}{\textsc{opt}(\mathbf n)} \big\vert{\eventtwo, \eventthree}\right]$.

\textit{Part 1: bounding $\Pr(\eventtwo)$. } By Hoeffding's inequality, we have
\[
\Pr(s_1 \in [hp-\sqrt{h},hp+\sqrt{h}]) \geq 1-\frac{1}{h^2} \geq 1-\frac{1}{m},
\]
\[
\Pr(s_2 \in [\ell p-\sqrt{\ell},\ell p+\sqrt{\ell}]) \geq 1-\frac{1}{\ell^2} \geq 1-\frac{1}{m}.
\]
The last steps of both equations are because in region $\mathcal R_3$, we have $h \geq \mzero $, $\ell \geq \mzero $. Then, because $s_1$ and $s_2$ are independent of each other, we have 
\[
\Pr(\eventtwo) = \Pr\left(s_1 \in [hp-\sqrt{h},hp+\sqrt{h}]\right) \cdot \Pr\left(s_2 \in [\ell p-\sqrt{\ell},\ell p+\sqrt{\ell}]\right) \geq (1-\frac{1}{m})^2.
\]

\textit{Part 2: bounding $\Pr(\eventthree \vert \eventtwo)$. } 
Conditioned on $s_1$, we know 
$\widehat{r_1}=\frac{\text{Bin}(s_1,r_1)}{s_1}$, and similarly conditioned on $s_2$,  $\widehat{r_2}=\frac{\text{Bin}(s_2,r_2)}{s_2}$.
Then, Hoeffding's inequality implies that 
\begin{align}
\Pr\left(\widehat r_1 \in [r_1-\frac{1}{\sqrt{ s_1}},r_1 +\frac{1}{\sqrt{s_1}}]~\Big \vert ~ s_1 = \tilde s_1\right) \geq 1-\frac{1}{\tilde s_1^2}\,. \label{eq:concentrate_r}
\end{align}
Now consider event  $\eventtwo$ that implies that  $s_1 \geq hp-\sqrt{h}$. As shown in in the proof of Lemma \ref{lem:cr2case2},  $hp-\sqrt{h}$ is an increasing function for any  $h \geq \mzero $, and hence under event  $\eventtwo$,  $s_1 \geq hp-\sqrt{h} \geq \mzero p-m^{1/4}$. Applying this lower bound on $s_1$ to Equation \eqref{eq:concentrate_r} leads to
\[
\Pr\left(\widehat{r_1} \in [r_1-\frac{1}{\sqrt{s_1}},r_1+\frac{1}{\sqrt{s_1}} ] ~\big\vert ~ \eventtwo \right) \geq 1-\frac{1}{(\mzero p-m^{1/4})^2}.
\]
Note that the same inequality holds $r_2$ and its estimate. As a result,  we  have
\[
\Pr(\eventthree \vert \eventtwo) \geq (1-\frac{1}{(\mzero p-m^{1/4})^2})^2.
\]

\textit{Part 3: bounding $\Pr(\mathcal G(\boldsymbol \psi; (h, \ell)))$ conditioned on events $\eventtwo$ and  $\eventthree$. } We show that conditioned on $\eventtwo$ and $\eventthree$, for any realization $\mathbf{s}$, $\Pr(\mathcal G(\boldsymbol \psi; (h, \ell)))=1$. We show this by proving $\widehat{r_1} \geq \widehat{r_2}$ always happens under events  $\eventtwo$ and  $\eventthree$. Under events $\eventtwo$ and $\eventthree$, we know that (i) $s_1, s_2 \geq \mzero p-m^{1/4}$ and (ii) $\widehat{r_1} \in [r_1-\frac{1}{\sqrt{s_1}},r_1+\frac{1}{\sqrt{s_1}}  ]$ and $\widehat{r_2} \in [r_2-\frac{1}{\sqrt{s_2}},r_2+\frac{1}{\sqrt{s_2}}]$. To show the result, we confirm that the lower bound on $\widehat r_1$ is smaller than the upper bound on $\widehat r_2$. That is, 
\[r_2+\frac{1}{\sqrt{s_2}} \le r_1-\frac{1}{\sqrt{s_1}} \qquad \text{for any $s_1, s_2  \geq \mzero p-m^{1/4}$}\]
Given the lower bound on $s_1$ and $s_2$, it suffices to show that 
\[r_2+\frac{1}{\sqrt{\mzero p-m^{1/4}}} \le r_1-\frac{1}{\sqrt{\mzero p-m^{1/4}}} ~~\Longleftrightarrow ~~ r_1-r_2\ge \frac{2}{m^{1/8}\sqrt{m^{1/4}p-1}} \]
Be definition of $\mTwo$, the above inequality holds for any $m\ge \mTwo$. Recall that $\mTwo = \mTwoExp$.

\textit{Part 4: bounding $\inf_{(h, \ell)\in \mathcal R_3} \E_{\mathbf n}\big[\frac{\rewGood}{\textsc{opt}(\mathbf n)} \big\vert{\eventtwo, \eventthree}\big]$. } First, observe that the ratio $\frac{\rewGood}{\textsc{opt}(\mathbf n)}$ is independent with $\eventthree$ because both the reward under the good event $\rewGood$,  and the optimal reward, i.e., $\textsc{opt}(\mathbf n)$, do not depend on $\widehat{r_1}$ and $\widehat{r_2}$. We can simply delete $\eventthree$, and bound $\inf_{(h, \ell)\in \mathcal R_3} \E_{\mathbf n}\big[\frac{\rewGood}{\textsc{opt}(\mathbf n)} \big\vert{\eventtwo}\big]$.

We begin by partitioning region $\mathcal R_3$ based the number of type $1$ agents, (i.e.. $h$): 
\[\rl =\{(h, \ell): (h, \ell)\in \mathcal R_3, h< \hzero  \}\, \qquad  \text{and} \qquad \rh =\{(h, \ell): (h, \ell)\in \mathcal R_3, h\ge \hzero  \}, \]
where  $\hzero =\hzeroExp$. Then, $\inf_{(h,\ell) \in \mathcal R_3}\E_{\mathbf n}\big[\frac{\rewGood }{\textsc{opt}(\mathbf n)} \big\vert \eventtwo \big]$ is equal to the minimum of $\inf_{(h,\ell) \in \rl}\E_{\mathbf n}\big[\frac{\rewGood }{\textsc{opt}(\mathbf n)}\big\vert \eventtwo \big]$ and 
$\inf_{(h,\ell) \in \rh}\E_{\mathbf n}\big[\frac{\rewGood }{\textsc{opt}(\mathbf n)}\vert \eventtwo \big]$.

\begin{lemma}\label{lem:copylow} Let $\rl =\{(h, \ell): (h, \ell)\in \mathcal R_3, h< \hzero  \}$. We then have 
\[
\inf_{(h,\ell) \in \rl}\E_{\mathbf n}\big[\frac{\rewGood }{\textsc{opt}(\mathbf n)}\big\vert \eventtwo \big] \geq 1-\frac{1}{\sqrt{m}}\,,
\]
where $\hzero =\hzeroExp$ and event $\eventtwo$ is defined in Equation \eqref{eq:event_s}. 
\end{lemma}
The proof of Lemma \ref{lem:copylow} is similar to the proof of Lemma  \ref{lem:cr2case2}, and is presented in Section \ref{sec:proof_lem:copylow}.

\begin{lemma} \label{lem:copyhigh}
Let $\rh =\{(h, \ell): (h, \ell)\in \mathcal R_3, h \geq \hzero  \}$. We then  have
\[
\inf_{(h,\ell) \in \rh}\E_{\mathbf n}\big[\frac{\rewGood }{\textsc{opt}(\mathbf n)}\big\vert \eventtwo \big] \geq \frac{h_1}{m},
\]
where $\mOne =\mOneExp$,  $\beta=0.4215 \cdot \frac{p^2+(1-p)^2}{p(1-p)}$,  and $\hzero =\hzeroExp$. Further, event  $\eventtwo$ is defined in Equation \eqref{eq:event_s}. 
\end{lemma}

The proof of Lemma \ref{lem:copyhigh} is similar to the proof of Lemma \ref{lem:cr2case3} and for completeness, is presented in Section \ref{sec:proof_lem:copyhigh}.  Lemmas \ref{lem:copylow} and \ref{lem:copyhigh} imply that    \[
    \inf_{(h, \ell)\in \mathcal R_3}~\E_{\mathbf n}
    \left[\frac{\rewGood}{\textsc{opt}(\mathbf n)} \right] \geq \min\left\{1-\frac{1}{\sqrt{m}}, \frac{\mOne }{m} \right\}=\X\,.
\]

Finally, by the results in parts 1-4, we have
\[
\inf_{(h, \ell)\in \mathcal R_3}~\E_{\mathbf n}
    \left[\frac{\rewGood}{\textsc{opt}(\mathbf n)}\right] \geq \CRThreeExpOne\,,
\]  which is the desired result.

\subsubsection{Proof of Lemma \ref{lem:cr3smallm}}
In Equation \eqref{eq:CR}, we have 
\begin{align}\nonumber \inf_{(h, \ell)\in \mathcal R_3}~\E_{\mathbf{s}}
    \left[\frac{\rewGen}{\textsc{opt}(\mathbf n)}\right] &=\inf_{(h, \ell)\in \mathcal R_3} \Bigg\{ \E_{\mathbf n}\left[\frac{\rewGood \cdot \Pr(\mathcal G(\boldsymbol \psi; (h, \ell)))}{\textsc{opt}(\mathbf n)}\right]\\ \label{eq:CR3EQ4'} 
    &+\E_{\mathbf n}\left[\frac{\rewBad \cdot (1-\Pr(\mathcal G(\boldsymbol \psi; (h, \ell))))}{\textsc{opt}(\mathbf n)}\right]\Bigg\}\,.
    \end{align}

Define event $\eventfive=\{s_1 \neq 0\} \cap \{s_2 \neq 0\}$. We can lower bound Equation \eqref{eq:CR3EQ4'} by 
\begin{align}\nonumber \inf_{(h, \ell)\in \mathcal R_3}~\E_{\mathbf{s}}
    \left[\frac{\rewGen}{\textsc{opt}(\mathbf n)}\right] &\geq \inf_{(h, \ell)\in \mathcal R_3}~\E_{\boldsymbol \psi}
    \left[\frac{\rewGen}{\textsc{opt}(\mathbf n)}\big \vert{\eventfive}\right]\Pr(\eventfive) \\&\geq \nonumber \inf_{(h, \ell)\in \mathcal R_3} \Bigg\{ \E_{\mathbf n}\left[\frac{\rewGood \cdot \Pr(\mathcal G(\boldsymbol \psi; (h, \ell)))}{\textsc{opt}(\mathbf n)} \big \vert{\eventfive}\right]\\ \label{eq:CR3EQ4} 
    &+\E_{\mathbf n}\left[\frac{\rewBad \cdot (1-\Pr(\mathcal G(\boldsymbol \psi; (h, \ell))))}{\textsc{opt}(\mathbf n)}\big\vert{\eventfive}\right]\Bigg\}\,\Pr(\eventfive).
    \end{align}
Note that
\begin{align}\nonumber
\Pr(\eventfive) &= 1-\Pr(\eventfivec) \geq 1-\Pr(s_1=0)-\Pr(s_2=0) \\& \label{eq:pevent5} = 1-{h \choose 0}(1-p)^h-{\ell \choose 0}(1-p)^{\ell} \geq 1-2{\sqrt{m} \choose 0}(1-p)^{\sqrt{m}}=1-2(1-p)^{\sqrt{m}} .
\end{align}
We denote $\V$ as $1-2(1-p)^{\sqrt{m}}$.
Furthermore, conditioned on event $\eventfive$, it is obvious that the probability of the good event $\mathcal G(\boldsymbol \psi; (h,\ell)) =\{\widehat r_1 \geq \widehat r_2\}$ is greater than or equal $1/2$. This is because when event $\eventfive$ holds, (i) our estimate for $r_i$, $i\in \{1, 2\}$ is simply the sample average of the $s_i$ realized rewards, and (ii) we have $r_1> r_2$. In light of Equation \eqref{eq:CR3EQ4}, if $ \E_{\mathbf n}\Big[\frac{\rewGood}{\textsc{opt}(\mathbf n)}\big \vert{\eventfive} \Big] <  \E_{\mathbf n}\Big[\frac{\rewBad}{\textsc{opt}(\mathbf n)} \big \vert{\eventfive} \Big]$, then Equation \eqref{eq:CR3EQ4} is lower bounded by setting the probability of the good event to be $1$. 
\begin{align} \inf_{(h, \ell)\in \mathcal R_3}~\E_{\mathbf{s}}
    \left[\frac{\rewGen}{\textsc{opt}(\mathbf n)}\big \vert{\eventfive} \right] &\geq 
    \label{eq:CR3c1}\inf_{(h, \ell)\in \mathcal R_3} \E_{\mathbf n}\left[\frac{\rewGood}{\textsc{opt}(\mathbf n)}\big \vert{\eventfive}\right].
    \end{align}
Otherwise, Equation \eqref{eq:CR3EQ4} is lower bounded by setting the probability of good event to be $\frac{1}{2}$.
\begin{align} \inf_{(h, \ell)\in \mathcal R_3}~\E_{\mathbf{s}}
    \left[\frac{\rewGen}{\textsc{opt}(\mathbf n)}\big \vert{\eventfive}\right] &\geq
     \inf_{(h, \ell)\in \mathcal R_3} \frac{1}{2}\E_{\mathbf n}\left[\frac{\rewGood}{\textsc{opt}(\mathbf n)}\big \vert{\eventfive}\right]+ \label{eq:CR3c2}\frac{1}{2}\E_{\mathbf n}\left[\frac{\rewBad }{\textsc{opt}(\mathbf n)}\big \vert{\eventfive}\right]
    \end{align}

Given Equations \eqref{eq:CR3c1} and \eqref{eq:CR3c2}, we need to bound $\E_{\mathbf n}\Big[\frac{\rewGood}{\textsc{opt}(\mathbf n)}\big \vert{\eventfive}\Big]$ and 
$\E_{\mathbf n}\Big[\frac{\rewBad}{\textsc{opt}(\mathbf n)}\big \vert{\eventfive}\Big]$. Notice that the former ratio can be lower bounded by
\begin{align*}
    \E_{\mathbf n}\Big[\frac{\rewGood}{\textsc{opt}(\mathbf n)}\big \vert{\eventfive}\Big] &\geq \E_{\mathbf n}\Big[\frac{\rewGood}{\textsc{opt}(\mathbf n)}\big \vert{\eventfive, \eventtwo}\Big]\Pr(\eventtwo) \\&= \E_{\mathbf n}\Big[\frac{\rewGood}{\textsc{opt}(\mathbf n)}\big \vert{\eventtwo}\Big]\Pr(\eventtwo)\,,
\end{align*}
where  event $\eventtwo=\left\{s_1 \in [hp-\sqrt{h},hp+\sqrt{h}]\} \cap \{s_2 \in [\ell p-\sqrt{\ell},\ell p+\sqrt{\ell}]\right\}$. 
The equality holds because $\eventtwo \subset \eventfive$. Then, by {part 1} and {4} of the proof of Lemma \ref{lem:cr3largem}, we have the former ratio is bounded by 
\begin{equation} \label{eq:cr3finalgood}
\E_{\mathbf n}\Big[\frac{\rewGood}{\textsc{opt}(\mathbf n)}\big \vert{\eventfive}\Big] \geq (1-\frac{1}{m})^2\X,
\end{equation}
where
$\X=\min\left\{1-\frac{1}{\sqrt{m}}, \frac{\mOne }{m} \right\}\,.$

Bounding the latter ratio, i.e., $\E_{\mathbf n}\Big[\frac{\rewBad}{\textsc{opt}(\mathbf n)}\big \vert{\eventfive}\Big]$,  is quite involved because when the bad event $\mathcal G^{\mathsf{C}}(\boldsymbol \psi; (h, \ell))$ happens, Algorithm \ref{alg:ptas} may over-protect or under-protect type $2$ agents given the randomness in the sample information. Over-protecting type $2$ agents leads to reject a certain number of type $1$ agents and waste some of the resources. As $\ell$ is large in this case, the number of wasted resources  can be large too. Under-protecting type $2$ agents leads to accept many type $1$ agents, which is a good thing for the algorithm.  Nonetheless, we provide the lower bound of the ratio below and show that it is still not too small in this case once we take expectation with respect to $\mathbf{s}$.   

We begin by partitioning region $\mathcal R_3$ based the number of type $2$ agents chosen by the adversary, (i.e. $\ell$):
\[\overline{\mathcal{R}_3} =\{(h, \ell): (h, \ell)\in \mathcal R_3, \ell > \ell_0 \}\, \qquad \underline{\mathcal{R}_3} =\{(h, \ell): (h, \ell)\in \mathcal R_3, \ell \le \ell_0 \}\,, \]
where $\ell_0=\ellzeroExp$. Then, $\inf_{(h,\ell) \in \mathcal R_3}\E_{\mathbf n}\big[\frac{\rewBad}{\textsc{opt}(\mathbf n)}\big\vert{\eventfive}\big]$ is equal to the minimum of $\inf_{(h,\ell) \in \overline{\mathcal{R}_3}}\E_{\mathbf n}\big[\frac{\rewBad}{\textsc{opt}(\mathbf n)}\big\vert{\eventfive}\big]$ and $\inf_{(h,\ell) \in \underline{\mathcal{R}_3}}\E_{\mathbf n}\big[\frac{\rewBad}{\textsc{opt}(\mathbf n)}\big\vert{\eventfive}\big]$. We bound the first ratio (i.e., the one concerns $\overline{\mathcal{R}_3}$) in Lemma \ref{lem:cr3low1} and the second ratio  (i.e., the one concerns $\underline{\mathcal{R}_3}$) will be bounded in Lemma \ref{lem:cr3low2}.

\begin{lemma} \label{lem:cr3low1}
Let $\overline{\mathcal{R}_3} =\{(h, \ell): (h, \ell)\in \mathcal R_3, \ell > \ell_0 \}$, where $\mathcal R_3 =\{(h,\ell):  h\ge \mzero , \ell \ge \mzero  \}$, and $\ell_0=\ellzeroExp$. Then,
\[
\inf_{(h, \ell)\in \overline{\mathcal{R}_3}}\E_{\mathbf n}\left[\frac{\rewBad }{\textsc{opt}(\mathbf n)}\big\vert{\eventfive}\right] \geq (1-\frac{1}{\ell_0^2})\alpha\,,
\]
where $\alpha = \frac{r_2}{r_1}$.
\end{lemma}

In this case, we show that with high probability, the number of type $2$ agents in the online arrival sequence $n_2 \geq m$. Then, the algorithm accepts $m$ type $2$ agents, and the ratio is bounded by $\alpha$.

\begin{lemma} \label{lem:cr3low2}
Let $\underline{\mathcal{R}_3} =\{(h, \ell): (h, \ell)\in \mathcal R_3, \ell \leq \ell_0 \}$, where $\mathcal R_3 =\{(h,\ell):  h\ge \mzero , \ell \ge \mzero  \}$, and  $\ell_0=\ellzeroExp$. Let $\beta=0.4215 \cdot \frac{p^2+(1-p)^2}{p(1-p)}$ and $\ell_1= \sqrt{p(1-p)\ell_0} +\frac{\beta}{\sqrt{m}} +\ell_0 p$. Then,
\[
\inf_{(h, \ell)\in \underline{\mathcal{R}_3}}\E_{\mathbf n}\left[\frac{\rewBad }{\textsc{opt}(\mathbf n)}\big\vert{\eventfive}\right] \geq \min\left\{\alpha, 1-\frac{1-p}{pm}\ell_1\right\}.
\]
\end{lemma}

To show Lemma \ref{lem:cr3low2}, we define event $\eventfour=\{n_2 \geq \frac{1-p}{p}s_2 \}$. Under this event the number of type $2$ agents is greater than or equal to the protection level assigned to this type of agents. Then, we provide the lower bound of  $\inf_{(h, \ell)\in \underline{\mathcal{R}_3}}\E_{\mathbf n}\big[\frac{\rewBad}{\textsc{opt}(\mathbf n)}\big\vert{\eventfive}\big]$ by further conditioning  on either $\eventfour$ or $\eventfourc$. For the part with $\eventfour$, we do not reject type $1$ agents unless there is no resource remaining. Therefore, the lower bound can be bounded by $\alpha$. For the part with $\eventfourc$, we directly calculate how many resources we waste, and then provide a lower bound. 

Finally, combining the results from Lemmas \ref{lem:cr3low1} and \ref{lem:cr3low2},  we have

\begin{equation} \label{eq:cr3finalbad}
\inf_{(h, \ell)\in \mathcal R_3} \E_{\mathbf n}\left[\frac{\rewBad}{\textsc{opt}(\mathbf n)}\big\vert{\eventfive}\right] \geq \min\left\{(1-\frac{1}{\ell_0^2})\alpha, 1-\frac{1-p}{pm}\ell_1 \right\}.
\end{equation}
Combining the results from Equations \eqref{eq:CR3c1}, and \eqref{eq:cr3finalgood}, we have if $(1-\frac{1}{m})^2\X \leq \min\left\{(1-\frac{1}{\ell_0^2})\alpha, 1-\frac{1-p}{pm}\ell_1 \right\}$, 
\begin{equation} \label{eq:CR3FINAL1}
\inf_{(h, \ell)\in \mathcal R_3}~\E_{\mathbf n}
    \left[\frac{\rewGen}{\textsc{opt}(\mathbf n)}\big\vert{\eventfive} \right] \geq (1-\frac{1}{m})^2\X.
\end{equation}
Otherwise, by Equations \eqref{eq:CR3c2}, and \eqref{eq:cr3finalbad}, we have
\begin{equation} \label{eq:CR3FINAL2}
\inf_{(h, \ell)\in \mathcal R_3}~\E_{\mathbf n}
    \left[\frac{\rewGen}{\textsc{opt}(\mathbf n)}\big\vert{\eventfive} \right] \geq \frac{1}{2}(1-\frac{1}{m})^2\X+\frac{1}{2}\min\{(1-\frac{1}{\ell_0^2})\alpha, 1-\frac{1-p}{pm}\ell_1\}.
\end{equation}

Therefore, combining the results in Equations \eqref{eq:CR3EQ4}, \eqref{eq:pevent5}, and \eqref{eq:CR3FINAL1}, \eqref{eq:CR3FINAL2}, we have 
\[
\inf_{(h, \ell)\in \mathcal R_3}~\E_{\mathbf n}
    \left[\frac{\rewGen}{\textsc{opt}(\mathbf n)}\right] \geq \V \cdot \min \Bigg\{(1-\frac{1}{m^2})\X, \frac{1}{2}(1-\frac{1}{m^2})\X+\frac{1}{2}\min\left\{(1-\frac{1}{\ell_0^2})\alpha, 1-\frac{1-p}{pm}\ell_1\right\} \Bigg\},
\]
where recall that $\V=1-2(1-p)^{\sqrt{m}}$.


\subsubsection{Proof of Lemma \ref{lem:copylow}}

\label{sec:proof_lem:copylow}
Recall that $\eventzero=\{n_2>m-\frac{1-p}{p}s_1\}$. We provide a lower bound of $\E_{\mathbf n}\left[\frac{\rewGood }{\textsc{opt}(\mathbf n)}\big\vert{\eventtwo}\right]$ by considering the conditional expectation on $\eventzeroc$.  
\begin{align}\nonumber
    &\inf_{(h,\ell) \in \rl}\E_{\mathbf n}\left[\frac{\rewGood }{\textsc{opt}(\mathbf n)}\big\vert{\eventtwo}\right] 
    \\&\geq \label{eq:condcr2twocopy} \inf_{(h,\ell) \in \rl}\E_{\mathbf n}\left[\frac{\rewGood }{\textsc{opt}(\mathbf n)}\big \vert{\eventtwo,\eventzero}\right]\Pr(\eventzero \vert \eventtwo).
\end{align}
In light of Equation \eqref{eq:condcr2twocopy}, we will provide a lower bound for $\E_{\mathbf n}\left[\frac{\rewGood }{\textsc{opt}(\mathbf n)}\big \vert{\eventtwo, \eventzero}\right]$ and $\Pr(\eventzero \vert \eventtwo)$.

\textit{Part 1: bounding $\E_{\mathbf n}\left[\frac{\rewGood }{\textsc{opt}(\mathbf n)}\big \vert{\eventtwo, \eventzero}\right]$.} For a fixed realization of $\mathbf s$ (or equivalently $\mathbf n$), we have
\begin{align} \notag
\rewGood&=\min\left\{n_2,\left(m-\frac{1-p}{p}s_1\right)^{+}\right\}\cdot r_2\\
&+\min\left\{m-\min\left\{n_2,\left(m-\frac{1-p}{p}s_1\right)^{+}\right\},n_1 \right\}\cdot r_1\,\label{eq:good_reward_2copy}
\end{align}
Conditional on $\eventzero=\{n_2 \leq m-\frac{1-p}{p}s_1\}$, we have $\min\{n_2,(m-\frac{1-p}{p}s_1)^{+}\}=n_2$, and $\min\{m-\min\{n_2,(m-\frac{1-p}{p}s_1)^{+}\},n_1 \}=\min\{m-n_2,n_1 \}$. Then, we have
\begin{align}\nonumber
              \frac{\rewBad}{\textsc{opt}(\mathbf n)} &\ge \frac{\rewBad}{n_1r_1+\min\{m-n_1,n_2\}r_2}\\
              &\nonumber=\frac{n_2r_2+\min\{m-n_2,n_1\}r_1}{n_1r_1+\min\{m-n_1,n_2\}r_2} \\&\nonumber \ge  \min\left\{\frac{n_1r_1+n_2r_2}{n_1r_1+n_2r_2},\frac{(m-n_2)r_1+n_2r_2}{n_1r_1+(m-n_1)r_2} \right\} \\&= \label{eq:cr2part1copy} \frac{(m-n_2)r_1+n_2r_2}{n_1r_1+(m-n_1)r_2} \geq \frac{(m-\mzero )r_1}{mr_1} =1-\frac{1}{\mzero}.
            \end{align}
            The second inequality is because when $\min\{m-n_2,n_1\}=n_1$, we have $\min\{m-n_1,n_2\}=n_2$, and when $\min\{m-n_2,n_1\}=m-n_2$, we have $\min\{m-n_1,n_2\}=m-n_1$.
            
            \textit{Part 2: bounding $\Pr(\eventzero \vert \eventtwo)$.} We want to show that $\Pr(\eventzero \vert \eventtwo)$ equals $1$ for $(h,
            \ell) \in \rl$. As $\ell\in \mzero $ for any $(h,
            \ell) \in \rl$, we have
            $n_2 \leq \mzero $.  Then, because $\eventzero = \{n_2 \le m-\frac{1-p}{p}s_1\}$, we have
            \begin{align*}
            \Pr(\eventzero \vert \eventtwo ) &= \Pr\left(\frac{1-p}{p}s_1 \leq m-n_2 \vert \eventtwo \right) \geq \Pr\left(\frac{1-p}{p}s_1 \leq m-\mzero \vert \eventtwo \right) \\&= \Pr\left(s_1 \leq \frac{p(m-\mzero )}{1-p} \vert \eventtwo \right).
            \end{align*}
            Next, we show that conditional on $\eventtwo$, the upper bound of $s_1$ is no more than $\frac{p(m-\mzero )}{1-p}$.
            If $\eventtwo$ happens, we have $s_1 \leq hp+ \sqrt{h}$. We then have \begin{align*}
               s_1 \le hp+\sqrt{h} \Rightarrow  s_1 \le \hzero p+\sqrt{\hzero } \leq s_1 \Rightarrow \frac{p(m-\mzero )}{1-p}\,, 
                \end{align*}
                where the second  inequality holds because $hp+ \sqrt{h}$ is an increasing function, and $h<\hzero $ in region $\underline{\mathcal R_2}$. These imply that $hp+\sqrt{p}\le \hzero  p+\sqrt{p_0}$.  The last inequality holds because $\frac{p(m-\mzero )}{1-p}\ge \hzero  p+\sqrt{\hzero }$. Recall that $\hzero =\hzeroExp$. This means that 
            \begin{equation} \label{eq:cr2part2copy}
                \Pr(\eventzeroc \vert \eventtwo) = 1.
            \end{equation}

 By Equations \eqref{eq:cr2part1copy} and \eqref{eq:cr2part2copy}, we have 
        \begin{equation}
            \inf_{(h,\ell) \in \rl}\E_{\mathbf n}\left[\frac{\rewGood }{\textsc{opt}(\mathbf n)}\big \vert{\eventzero}\right]\Pr(\eventzero) \geq (1-\frac{1}{\sqrt{m}})\,,
        \end{equation}

\subsubsection{Proof of Lemma \ref{lem:copyhigh}} \label{sec:proof_lem:copyhigh}  Recall that $\eventone=\{n_1>\frac{1-p}{p}s_1\}$. We bound $\E_{\mathbf n}\left[\frac{\rewGood }{\textsc{opt}(\mathbf n)} \big \vert{\eventtwo}\right] $ by considering the conditional expectation on $\eventone$ and $\eventonec$ respectively.
\begin{align} \nonumber
    \inf_{(h,\ell) \in \rh }&\E_{\mathbf n}\left[\frac{\rewGood }{\textsc{opt}(\mathbf n)} \big \vert{\eventtwo}\right]    \\
    \nonumber &\geq \inf_{(h,\ell) \in \rh } \Bigg\{ \min\Bigg\{ \E_{\mathbf n}\left[\frac{\rewGood }{\textsc{opt}(\mathbf n)} \big\vert{\eventtwo, \eventone}\right], \\ \label{eq:condcr2secondcopy}
     &\E_{\mathbf n}\left[\frac{\rewGood }{\textsc{opt}(\mathbf n)}\big\vert{\eventtwo, \eventone^{\mathsf{C}}}\right]
   \Bigg\}\Bigg\}\,
\end{align}

We will provide the lower bound of $\inf_{(h,\ell) \in \rh} \E_{\mathbf n}\big[\frac{\rewGood }{\textsc{opt}(\mathbf n)}\big\vert{\eventtwo,\eventone} \big]$ in the first part  of the proof, and the lower bound of $\inf_{(h,\ell) \in \rh} \E_{\mathbf n}\big[\frac{\rewGood }{\textsc{opt}(\mathbf n)}\big\vert{\eventtwo, \eventone^{\mathsf{C}}} \big]$ in the second part of the proof.

\textit{Part 1: bounding $\inf_{(h,\ell) \in \rh} \E_{\mathbf n}\big[\frac{\rewGood }{\textsc{opt}(\mathbf n)}\big\vert{\eventtwo, \eventone} \big]$.} We then  have
\begin{align*}
    &\inf_{(h,\ell) \in \rh} \E_{\mathbf n}\big[\frac{\rewGood }{\textsc{opt}(\mathbf n)}\big\vert{\eventtwo, \eventone} \big] \\&\geq \inf_{(h,\ell) \in \rh} \min \Bigg\{ \E_{\mathbf n}\big[\frac{\rewGood }{\textsc{opt}(\mathbf n)}\big\vert{\eventtwo,\eventone, \eventzeroc} \big],\\& \E_{\mathbf n}\big[\frac{\rewGood }{\textsc{opt}(\mathbf n)}\big\vert{\eventtwo, \eventone, \eventzero} \big] \Bigg\}\,,
\end{align*}
where we recall $\eventzeroc=\{n_2>m-\frac{1-p}{p}s_1\}$.

To provide a lower bound for $\E_{\mathbf n}\big[\frac{\rewGood }{\textsc{opt}(\mathbf n)}\big\vert{\eventtwo, \eventone, \eventzeroc} \big]$, we note that under event 
$\eventzeroc$, we have  $n_2>m-\frac{1-p}{p}s_1$, and hence   
            \begin{equation} 
            \E_{\mathbf n}\left[\frac{\rewGood }{\textsc{opt}(\mathbf n)}\big\vert{\eventtwo, \eventone, \eventzeroc} \right]\geq \frac{(m-\frac{1-p}{p}s_1)^{+}r_2+\min\{m,\frac{1-p}{p}s_1\}r_1}{mr_1} \geq \frac{(m-n_2) r_1}{mr_1} \geq 1-\frac{1}{\mzero}.
            \end{equation}
            
Then, we give a lower bound to $\E_{\mathbf n}\big[\frac{\rewGood }{\textsc{opt}(\mathbf n)}\big\vert{\eventtwo, \eventone, \eventzero} \big]$. Conditional on $\eventone$ and  $\eventzero$, we have $\frac{1-p}{p}s_1 \geq m-n_2$ and $\frac{1-p}{p}s_1 \leq n_1$. That is, the number of type $1$ agents in the online arrival sequence is greater than or equal to the protection level of Algorithm \ref{alg:ptas} for type $1$ agents (i.e., $\frac{1-p}{p}s_1 \leq n_1$), and in addition, the algorithm ends up rejecting some of type $2$ agents as $\frac{1-p}{p}s_1 \geq m-n_2$. Therefore, 
we have
\[
\E_{\mathbf n}\left[\frac{\rewGood }{\textsc{opt}(\mathbf n)}\big\vert{\eventtwo, \eventone, \eventzero} \right] \geq \frac{(m-\frac{1-p}{p}s_1)^{+}r_2+\min\{m,\frac{1-p}{p}s_1\}r_1}{mr_1} \geq \frac{(m-n_2)r_1}{mr_1} \geq 1-\frac{1}{\mzero}\,,
\]
where the last inequality holds because $n_2\le \ell \le \mzero $ for any $(h,\ell)\in \rh$.

\textit{Part 2: bounding $\inf_{(h,\ell) \in \rh} \E_{\mathbf n}\big[\frac{\rewGood }{\textsc{opt}(\mathbf n)}\big\vert{\eventtwo, \eventone^{\mathsf{C}}} \big]$. } Recall that for a fixed realization of $\mathbf s$ (or equivalently $\mathbf n$), 
\begin{align*} \notag
\rewGood&=\min\left\{n_2,\left(m-\frac{1-p}{p}s_1\right)^{+}\right\}\cdot r_2\\
&+\min\left\{m-\min\left\{n_2,\left(m-\frac{1-p}{p}s_1\right)^{+}\right\},n_1 \right\}\cdot r_1\,.
\end{align*}

Conditional on $\eventone^{\mathsf{C}} = \{n_1 \leq \frac{1-p}{p}s_1 \}$, we do not receive enough number of type $1$ agents and thus, we have
\[
\min\left\{m-\min\left\{n_2,\left(m-\frac{1-p}{p}s_1\right)^{+}\right\},n_1 \right\}\geq \min\left\{\max\left\{m-n_2,\frac{1-p}{p}s_1\right\},n_1 \right\} =n_1\,.
\]
That is, Algorithm \ref{alg:ptas}  accepts all type $1$ agents in the online arrival sequence. Then, we have
\begin{align*}
\E_{\mathbf n}\Big[\frac{\rewGood }{\textsc{opt}(\mathbf n)}\big\vert{\eventtwo, \eventone^{\mathsf{C}}} \Big] &\geq \E_{\mathbf n}\Big[\frac{\min\left\{n_2,\left(m-\frac{1-p}{p}s_1\right)^{+}\right\}\cdot r_2+n_1r_1 }{mr_1}\big\vert{\eventtwo, \eventone^{\mathsf{C}}} \Big] \\&\geq \frac{1}{m}\E_{\mathbf n}\big[ n_1 \big\vert{\eventtwo, \eventone^{\mathsf{C}}} \big].
\end{align*}
We will show that $\E_{\mathbf n}\big[ n_1 \big\vert{\eventtwo, \eventone^{\mathsf{C}}} \big] \ge \mOne $, where $\mOne  = \mOneExp$ and $\beta=0.4215 \cdot \frac{p^2+(1-p)^2}{p(1-p)}$.

 Since $\eventone^{\mathsf{C}} = \{n_1 \leq \frac{1-p}{p}s_1 \}= \{h(1-p)-n_1 \geq 0 \}$,  we have
\begin{align*}
&\E_{\mathbf{n}}[n_1 \vert \eventonec, \eventtwo] \\&= h(1-p)-\E_{\mathbf{n}}[h(1-p)-n_1 \big \vert \{h(1-p)-n_1 \geq 0\}, \eventtwo] \\&=   h(1-p)-\E_{\mathbf{n}}[h(1-p)-n_1 \big \vert \{h(1-p)-n_1 \geq 0\} \cap \{h(1-p)-\sqrt{h} \leq n_1 \leq h(1-p)+\sqrt{h}\}] \\&= h(1-p)-\E_{\mathbf{n}}[h(1-p)-n_1 \big \vert \{0 \leq h(1-p)-n_1 \leq \sqrt{h}\}] \\&\geq h(1-p)-\E_{\mathbf{n}}[h(1-p)-n_1 \big \vert \{h(1-p)-n_1 \geq 0\}] \\&=   \E_{\mathbf{n}}[n_1 \vert \eventonec].
\end{align*}
The inequality is because  
\begin{align*}
&\E_{\mathbf{n}}[h(1-p)-n_1 \big \vert \{0 \leq h(1-p)-n_1 \leq \sqrt{h}\}] \\&\leq \E_{\mathbf{n}}[h(1-p)-n_1 \big \vert \{0 \leq h(1-p)-n_1 \leq \sqrt{h}\}]+\E_{\mathbf{n}}[h(1-p)-n_1 \big \vert \{h(1-p)-n_1 \geq \sqrt{h}\}] \\&= \E_{\mathbf{n}}[h(1-p)-n_1 \big \vert \{h(1-p)-n_1 \geq 0\}].
\end{align*}

Next, we have 
\begin{align} \notag
    \E_{\mathbf n}\big[ n_1 \big\vert{\eventone^{\mathsf{C}}} \big] &= h(1-p)-\E_{\mathbf n}\big[ h(1-p)-n_1 \big\vert{h(1-p)-n_1 \geq 0} \big]\\ \notag
       &\ge h(1-p)-\E_{\boldsymbol {s}}
    \left[|h(1-p)-n_1|\right]\big| -  \frac{\beta}{\sqrt{h}}\\ \notag
    &\ge h(1-p)-\sqrt{p(1-p)h}-  \frac{\beta}{\sqrt{h}}\\ \notag
      &\ge  \hzero (1-p)-\sqrt{p(1-p)\hzero }-  \frac{\beta}{\sqrt{\hzero }}
  \\ &= \notag
   \mOne \,,
\end{align} 
where the first inequality holds because by 
 \cite{nagaev2011bound}, we have
$
    \Big|\E_{\boldsymbol {n}}
    \Big[(h(1-p)-n_1) \big \vert \{h(1-p)-n_1 \geq 0\}\Big]\,-\E_{\boldsymbol {n}}
    \left[|h(1-p)-n_1|\right]\Big| \le \frac{\beta}{\sqrt{h}}
$. The second inequality follows from  \cite{berend2013sharp} that shows  
$
\E_{\boldsymbol n}\left[|h(1-p)-n_1|\right]\, \leq \sqrt{p(1-p)h}.
$ The third inequality holds because when $h>\hzero $, $h(1-p)-\sqrt{p(1-p)h}- \frac{\beta}{\sqrt{h}}$ is a non-decreasing function for any $p$. To see why note that  \begin{equation} \label{eq:finalcr2copy}
    \frac{\partial (h(1-p)-\sqrt{p(1-p)h}- \frac{\beta}{\sqrt{h}})}{\partial h} \geq (1-p)-\frac{1}{2}\sqrt{p(1-p)}\frac{1}{\sqrt{h}}\ge 0\,,
\end{equation}
where the last inequality holds for any  $h \geq \frac{1}{4}\frac{p}{1-p}$, and we have $\hzero  > \frac{m-\sqrt{m}}{1-p} > \frac{1}{4}\frac{p}{1-p}$ for all $m \geq 2$. 

\subsection{Proof of Lemma \ref{lem:cr3low1}}
As $\ell_0=\ellzeroExp$, by Hoeffding's inequality, we have
\[
\Pr(n_2 \geq m) = \Pr(\text{Bin}(\ell_0,1-p) \geq m) \geq 1-\frac{1}{\ell_0^2}.
\]
Given that $n_2 \geq m$, the worst case is where the decision-maker accepts $m$ type $2$ agents, and the CR is lower bounded by $\frac{mr_2}{mr_1}=\alpha$.

\subsubsection{Proof of Lemma \ref{lem:cr3low2}}
We have 
\begin{align}\nonumber
    \inf_{(h, \ell)\in \underline{\mathcal{R}_3}}\E_{\mathbf n}\left[\frac{\rewBad }{\textsc{opt}(\mathbf n)}\right]= \nonumber &\inf_{(h, \ell)\in \underline{\mathcal{R}_3}} \min \Bigg\{\E_{\mathbf n}\left[\frac{\rewBad }{\textsc{opt}(\mathbf n)}\big \vert \{n_2 \geq m \}\right], \\&\label{eq:lowtype} \E_{\mathbf n}\left[\frac{\rewBad }{\textsc{opt}(\mathbf n)}\big \vert \{n_2 < m \}\right] \Bigg\}
\end{align}

We have discussed in  the proof of Lemma \ref{lem:cr3low1} that $\E_{\mathbf n}\left[\frac{\rewBad }{\textsc{opt}(\mathbf n)}\big \vert \{n_2 \geq m \}\right] \geq \alpha$. Then, given Equation \eqref{eq:lowtype}, we provide a lower bound to $\E_{\mathbf n}\left[\frac{\rewBad }{\textsc{opt}(\mathbf n)}\big \vert \{n_2 < m \}\right]$.

Conditional on event $\{n_2 < m \}$, the algorithm accepts all $n_2$ type $2$ agents. If $n_2$ is larger than the protection level $\min\{m,\frac{1-p}{p}s_2\}$, the algorithm accepts as many type $1$ agents as possible with the remaining $m-n_2$ units of the resource. If $n_2$ is small than the protection level $\min\{m,\frac{1-p}{p}s_2\}$, the algorithm rejects some type $1$ agents because we can accept at most $(m-\frac{1-p}{p}s_2)^{+}$ type $1$ agents. Then, we have
\[
\E_{\mathbf n}\left[\frac{\rewBad }{\textsc{opt}(\mathbf n)}\big \vert \{n_2 < m \}\right]\geq \E_{\mathbf n}\left[\frac{n_2r_2+\min\{n_1,m-\max\{n_2,\min\{m,\frac{1-p}{p}s_2\}\} \}r_1 }{\min\{m,n_1\}r_1+\min\{n_2,(m-n_1)^{+} \}r_2}\right].
\]

We define the event $\eventfour=\{n_2 \geq \frac{1-p}{p}s_2 \}$. Then, we have 
\begin{align}\nonumber
    &\E_{\mathbf n}\left[\frac{n_2r_2+\min\{n_1,m-\max\{n_2,\min\{m,\frac{1-p}{p}s_2\} \} \}r_1 }{\min\{m,n_1\}r_1+\min\{n_2,(m-n_1)^{+} \}r_2}\right] \\&\geq \nonumber \min \Bigg\{ \E_{\mathbf n}\left[\frac{n_2r_2+\min\{n_1,m-\max\{n_2,\min\{m,\frac{1-p}{p}s_2\} \} \}r_1 }{\min\{m,n_1\}r_1+\min\{n_2,(m-n_1)^{+} \}r_2} \big \vert \eventfour \right], \\& \label{eq:CR3bad1} \E_{\mathbf n}\left[\frac{n_2r_2+\min\{n_1,m-\max\{n_2,\min\{m,\frac{1-p}{p}s_2\} \} \}r_1 }{\min\{m,n_1\}r_1+\min\{n_2,(m-n_1)^{+} \}r_2} \big \vert \eventfourc \right] \Bigg\}.
\end{align}
Conditioned on $\eventfour$, for any realization $\mathbf{s}$, we have 
\begin{align*}
\frac{n_2r_2+\min\{n_1,m-\max\{n_2,\min\{m,\frac{1-p}{p}s_2\} \} \}r_1 }{\min\{m,n_1\}r_1+\min\{n_2,(m-n_1)^{+} \}r_2} &= \frac{n_2 r_2+\min\{n_1,(m-n_2)\} r_1}{\min\{m,n_1\}r_1+\min\{n_2,(m-n_1)^{+} \}r_2} \\& \geq \frac{mr_2}{mr_1} = \alpha.
\end{align*}

When $\eventfour$ happens, the number of arriving type $2$ agents is larger than the protection level, and the algorithm accepts all the arriving agents until there is no resource left. If $\textsc{OPT}$ does not use all resources, then the algorithm accepts everyone, and hence the above ratio will be one. If $\textsc{OPT}$ uses all resources, then the algorithm also uses all resources, and thus, this bound clearly holds.   

 Conditioned on $\eventfourc$, we have $n_2 < \frac{1-p}{p}s_2$. Then, 
\begin{align*}
\frac{n_2r_2+\min\{n_1,m-\max\{n_2,\min\{m,\frac{1-p}{p}s_2\} \} \}r_1 }{\min\{m,n_1\}r_1+\min\{n_2,(m-n_1)^{+} \}r_2} &= \frac{n_2r_2+\min\{n_1,(m-\frac{1-p}{p}s_2)^{+} \}r_1 }{\min\{m,n_1\}r_1+\min\{n_2,(m-n_1)^{+} \}r_2} \\&\geq \frac{n_2r_2+(m-\frac{1-p}{p}s_2 )^{+}r_1 }{mr_1}.
\end{align*}
The inequality is because if $n_1<m-\frac{1-p}{p}s_2$, then the algorithm accepts $n_2$ type $2$ agents and $n_1$ type $1$ agents, and the ratio is $1$. Then, we have
\begin{align} \nonumber
\E_{\mathbf n}\left[\frac{n_2r_2+\min\{n_1,m-\max\{n_2,\min\{m,\frac{1-p}{p}s_2\} \} \}r_1 }{\min\{m,n_1\}r_1+\min\{n_2,(m-n_1)^{+} \}r_2} \big \vert \eventfourc \right] &\geq \nonumber \E_{\mathbf n}\left[\frac{n_2r_2+(m-\frac{1-p}{p}s_2 )^{+}r_1 }{mr_1} \big \vert \eventfourc \right] \\& =\nonumber \E_{\mathbf n}\left[\frac{(m-\frac{1-p}{p}s_2 )^{+}r_1 }{mr_1} \big \vert \eventfourc \right] \\&\geq \nonumber \E_{\mathbf n}\left[\frac{(m-\frac{1-p}{p}s_2 )r_1 }{mr_1} \big \vert \eventfourc \right] \\&= 1-\frac{1-p}{pm}\E_{\mathbf n}\left[s_2 \big \vert \eventfourc \right].
\end{align}

Since $\eventfourc = \{s_2 \geq \ell p \}$,  we have
\begin{align} \notag
    \E_{\mathbf n}\left[s_2 \big \vert s_2 \geq \ell p \right] &= \E_{\mathbf n}\left[s_2-\ell p \big \vert s_2 \geq \ell p \right] + \ell p \\&\leq \nonumber \E_{\mathbf n}\left[\vert s_2-\ell p \vert \right] +\frac{\beta}{\sqrt{\ell}} +\ell p\\ \notag
       &\le \sqrt{p(1-p)\ell} +\frac{\beta}{\sqrt{\ell}} +\ell p\\ \notag
    &\le  \sqrt{p(1-p)\ell_0} +\frac{\beta}{\sqrt{m}} +\ell_0 p
  \\ &= \ell_1 \notag \,,
\end{align} 
where the first inequality holds because by 
 \cite{nagaev2011bound}, we have
$
    \Big|\E_{\mathbf n}\left[s_2-\ell p \big \vert s_2 \geq \ell p \right]-\E_{\mathbf n}\left[\vert s_2-\ell p \vert \right] \Big|  \le \frac{\beta}{\sqrt{\ell}}
$. The second inequality follows from  \cite{berend2013sharp} that shows  
$
\E_{\mathbf n}\left[\vert s_2-\ell p \vert \right] \leq \sqrt{p(1-p)\ell}.
$.

Finally, we have
\[
\inf_{(h, \ell)\in \underline{\mathcal{R}_3}}\E_{\mathbf n}\left[\frac{\rewBad }{\textsc{opt}(\mathbf n)}\right] \geq \min\left\{\alpha, 1-\frac{1-p}{pm}\ell_1\right\}.
\]

    \section{Proof of Proposition \ref{thm:3}} \label{subsec:proofthm3} 
To show the asymptotic result, we only need to consider the case that $m \geq \mTwo$. We can write the CR as
\[
    \text{CR}_{A} \ge  \min\left\{\CROne , \CRTwo , \overline{\text{CR}_3}\right\}.
    \]
To show the result, we will show how $\CROne , \CRTwo $, and $\overline{\text{CR}_3}$ scale with $m$. 
First it is easy to see that $\CROne =\CROneExp$ scales with $1- \Theta(1/(p\sqrt{m}))$. Similarly, it is easy to see that $\CRTwo = \CRTwoExp$ scales with $1- \Theta(1/\sqrt{m})$. This is because
$\hzero={\frac{m-\mzero }{1-p}+\frac{1}{2p^2}+\sqrt{\frac{1}{4p^4}+\frac{(m-\mzero )}{p^2(1-p)}}}$, and $\mOne =\mOneExp = m-(1+\sqrt{p})\sqrt{m}-o(\sqrt{m})$. Then, we have $\frac{h_1}{m}=1- \Theta(1/\sqrt{m})$.

To complete the proof, we show how 
\begin{align*}
\overline{\text{CR}_3}&=\CRThreeExpOne
\end{align*}
scales with $m$. Observe that $(1-\frac{1}{m})^2(1-\frac{1}{(\mzero p-m^{1/4})^2})^2= 1- \Theta(1/(p^2m))$.  When $p=\omega(1/\sqrt{m})$,  $\Theta(1/(p^2m))$ is dominated by $\Theta(1/(p\sqrt{m}))$. Further, by the analysis for $\CRTwo$, we know that $\X= \min\left\{1-\frac{1}{\sqrt{m}}, \frac{\mOne }{m} \right\}$ also scales with $1- \Theta(1/\sqrt{m})$. Putting these together, we conclude that $\text{CR}_{A} = 1-\Theta(1/(p\sqrt{m}))$, as desired.

\section{Proof of Theorem \ref{thm:upperbound}}\label{sec:proof:upperbound} 
We construct the following input family $\textbf{F}$: let  $\underline{h}=0$ and  $\bar{h}=pm$. The input family $\textbf{F}$  contains all $(h, \ell)$ such that $h \in [\underline{h},\bar{h}]$ and $\ell = \frac{10000\cdot m}{p}$.  For any $h \in [\underline{h},\bar{h}]$, we then denote $I_h$ as a random arrival sequence under which $n_2$ type $2$ agents arrive followed by  $n_1$ type $1$ agents, where we recall that   $n_1 \sim \text{Bin}(h,1-p)$, $n_2 \sim \text{Bin}(\ell,1-p)$. We denote $I_h$ as the random arrival sequences, where in $I_h$, $n_2$ type $2$ agents arrive followed by  $n_1$ type $1$ agents.

Our goal is to characterize an upper bound on the CR of any   deterministic algorithm under the family $\textbf{F}$. Consider a specific input $h\in [\underline h, \bar h]$. The decision maker gets the sample $(s_1,s_2)$ at the beginning of the allocation period, where $s_1\sim\text{Bin}(h,p)$, and $s_2\sim\text{Bin}(\ell,p)$. As we select $
\ell$ to be a super large constant, we have with probability $1$, $s_2 \gg m$. Therefore, under the family $\textbf{F}$, the decision-maker knows that there will be more than $m$ type $2$ agents showing up. Based on the structure that all type $2$ agents arrive before all type $1$ agents, any deterministic algorithm has to decide how many type $2$ agents they accept if they observe $(s_1,s_2)$. As $s_2$ is always larger than $m$, it does not  impact the decisions. Put differently, any deterministic algorithm can be represented by a mapping $z(\cdot)$: $[0, \bar{h}] \to [0, m]$, such that once they observe $(s_1,s_2)$, they will accept $m-z(s_1)$ type $2$ agents. 
 (As the number of type $1$ agents, $h$, is less than $m$, the number of type $1$ agents in the online arrival sequence $n_1<m$.)  Let $A$ be an algorithm associated with mapping $z$. We define the following loss function for algorithm $A$ under arrival sequence $I$ and sample information $\boldsymbol{s}$: 
\begin{equation} \label{eq:loss}
    \L_{A}\left(I, \boldsymbol s \right)=\textsc{opt}(I)-\rew_{A }\left(I, \boldsymbol s \right).
\end{equation}
(Note that since $r_1$ and $r_2$ are known, we can replace $\boldsymbol \psi= (\mathbf{s}, \boldsymbol{\rho})$ with $\mathbf{s}$.)
In the rest of this proof, we use the shorthand of 
 $\L_{}\left(I, \boldsymbol {s}, z \right)$ and  $\rew\left(I, \boldsymbol s, z\right)$ in place of 
 $\L_{A}(I, \mathbf{s})$ and $\rew_{A }\left(I, \boldsymbol {s}\right)$, respectively.  
Then, the CR of a deterministic  algorithm $A$ with mapping $z(\cdot)$ on family ${\mathbf{F}}$  is at least 
\begin{equation}
    \text{CR}(z) =  1-\sup_{(h,\ell) \in \mathbf{F}}~\E_{\boldsymbol {s}}
    \left[\frac{\L\left(I, \boldsymbol {s}, z \right)}{\textsc{opt}(I)}\right]\,. \label{eq:CR_F}
\end{equation}

In what follows, we provide an upper bound on  $\text{CR}(z)$ in Equation \eqref{eq:CR_F}.  To do so,  we consider  the following two cases:

\begin{itemize}[leftmargin=*]
    \item \textit{Case 1: $n_1 \geq z(s_1)$.}  In this case, the number of type $1$ agent arriving online (which is a random variable) is larger than the remaining number of resources after we accept $m-z(s_1)$ type $2$ agents. Therefore, the decision-maker accepts $m-z(s_1)$ type $2$ agents and $z(s_1)$ type $1$ agents. The loss (i.e., the difference between the optimal clairvoyant cumulative expected reward and what the decision-maker gets which is defined in Equation \eqref{eq:loss}) is then $(n_1-z(s_1))(r_1-r_2)$. 
    
    \item \textit{Case 2: $n_1 < z(s_1)$.} In this case, the number of type $1$ agent arriving online is less than the remaining number of resources after we accept $m-z(s_1)$ type $2$ agents. Therefore, the decision-maker accepts $m-z(s_1)$ type $2$ agents and $n_1$ type $1$ agents. Therefore, we waste $z(s_1)-n_1$ units of resources. In this case, the loss, defined in Equation \eqref{eq:loss}, is $(z(s_1)-n_1)r_2$.

\end{itemize}

Putting these two cases together, we have 
\begin{align*}
    \L\left(I, \boldsymbol s, z\right)&=(n_1-z(s_1))(r_1-r_2)\boldsymbol{1}(n_1 \geq z(s_1))+(z(s_1)-n_1)r_2\boldsymbol{1}(n_1 < z(s_1)) \\&= (n_1-z(s_1))r_1\boldsymbol{1}(n_1 \geq z(s_1))+(z(s_1)-n_1)r_2\,.
\end{align*}
Hence, an upper bound on the CR of any deterministic algorithm on family ${\mathbf{F}}$, denoted by $\overline{\text{CR}}$,  is at least  
\begin{align} \notag
    \overline{\text{CR}} &= \max_{z}\text{CR}(z) \\ \notag 
    &= 1-\min_{z}\max_{h \in [\underline{h},\bar{h}]}\E_{\boldsymbol s}
    \left[\frac{\L\left(I, \boldsymbol {s}, z\right)}{\textsc{opt}(I)}\right]\,
    \\&= 1-\min_z\max_{h \in [\underline{h},\bar{h}]}\E_{\boldsymbol s}
    \left[\frac{(n_1-z(s_1))r_1\boldsymbol{1}(n_1 \geq z(s_1))+(z(s_1)-n_1)r_2}{n_1r_1+(m-n_1)r_2}\right]\,, \label{eq:CR_z}
\end{align}
where $\text{CR}(z)$ is defined in Equation \eqref{eq:CR_F}.

Although we can solve this minmax optimization problem by some solver because the objective function is convex in every dimension, we cannot calculate the value of the optimal solution in a closed form.
Therefore, we develop the analysis below to show that $\overline{\text{CR}}$ is $1-\Theta(\frac{1}{p\sqrt{m}})$.

In our analysis, we first focus on a specific mapping:  $z(s_1)=\frac{1-p}{p}s_1$ for any $s_1 \in [0,\bar{h}]$. (Since $s_1 \le \bar h = mp$, we have $z(s_1)\le m$.)  For this mapping, we  calculate the CR of its associated  deterministic algorithm. Second, we compare any other feasible mappings with this mapping to show our result. We want to highlight that this mapping (i.e., $z(s_1)=\frac{1-p}{p}s_1$) is not the optimal mapping:  The optimal mapping should depend on $r_1$ and $r_2$, but we cannot find its closed form. This mapping is considered as a benchmark with which we  compare all other mappings.

\begin{lemma} [The CR under Mapping $z(\cdot)$] \label{lem:specific z}
Consider the following mapping $z: [0, \bar{h}] \to [0, m]$ such that  $z(s_1)=\frac{1-p}{p}s_1$ for any $s_1 \in [0, \bar{h}]$. Then, for any value of $m$ 
\[\max_{h \in [\underline{h},\bar{h}]}\E_{\boldsymbol s}
    \left[\frac{\L\left(I, \boldsymbol {s}, z\right)}{\textsc{opt}(I)}\right] ~\le~ \frac{r_1\sqrt{1-p}}{r_2}\frac{1}{p\sqrt{m}}  +\frac{\sqrt{2}\beta r_1}{p\sqrt{p}r_2} \frac{1}{m\sqrt{m}}\,,\]
    and 
    \[\max_{h \in [\underline{h},\bar{h}]}\E_{\boldsymbol s}
    \left[\frac{\L\left(I, \boldsymbol {s}, z\right)}{\textsc{opt}(I)}\right] ~\ge~ \frac{1}{2\sqrt{2}}\sqrt{1-p}\frac{1}{p\sqrt{m}}- \frac{\beta}{p\sqrt{p}}\frac{1}{m\sqrt{m}}\,,\]
    where  $\beta=0.4215\cdot \frac{p^2+(1-p)^2}{p(1-p)}$. This implies that  $\text{CR}(z) = 1- \Theta(\frac{1}{p\sqrt{m}})$.  
\end{lemma}

By Lemma \ref{lem:specific z}, we have under that particular mapping $z$, we have $\text{CR}(z)=1-\Theta(\frac{1}{p\sqrt{m}})$ over input family $\mathbf{F}$.
The next lemma shows that there does not exist any other mapping $\widehat z(\cdot)$ such that the CR over $\mathbf F$ is $1-o(\frac{1}{p\sqrt{m}})$.

\begin{lemma}[Other Mappings] \label{lem:other_mappings} 
Consider the setting in Section \ref{sec:proof:upperbound} and the input family $\mathbf F$. For this setting, there does not exist any  mapping $\widehat z(\cdot)$ such that the competitive ratio over $\mathbf F$ is $1-o(\frac{1}{p\sqrt{m}})$.
\end{lemma}

To show this, we compare all other mappings $\widehat z: [0, \bar{h}] \to [0, m]$,  with $z(s_1)=\frac{1-p}{p}s_1$, and using contradiction, we assume that there exists a mapping $\widehat z$ which achieves the CR of $1-o(\frac{1}{p\sqrt{m}})$. As $z$ achieves a CR of $1-\Theta(\frac{{1}}{p\sqrt{m}})$, the difference between the CR under $\widehat z$ and $z$ is at least $\Theta(\frac{{1}}{p\sqrt{m}})$. This gives constraints that $\widehat z$ and $z$ should be different on a set, and the difference between $\widehat z$ and $z$ on the points in this set is lower bounded. Then, we further claim that for a mapping $\widehat z$ which has a certain difference compared to $z$, the CR is at least $1-\Omega(\frac{1}{p\sqrt{m}})$, which gives a contradiction and shows the best  algorithm (mapping) can only achieve an asymptotic  CR of  $1-\Theta(\frac{{1}}{p\sqrt{m}})$.

\subsection{Proof of Lemma \ref{lem:specific z}} 
The proof has two parts. In the first part, we present an upper bound on $\max_{h \in [\underline{h},\bar{h}]}\E_{\boldsymbol s}
    \left[\frac{\L\left(I, \boldsymbol s, z\right)}{\textsc{opt}(I)}\right]$ and in the second part, we present a lower bound on the same quantity. 
    
\textit{First part: upper bound.} By definition,  $z(s_1)=\frac{1-p}{p}s_1$ for any $s_1 \in [0, \bar h]$, and hence we  have $n_1-z(s_1)=h-s_1-\frac{1-p}{p}s_1=h-\frac{1}{p}s_1$. Similarly, $z(s_1)-n_1=\frac{1}{p}s_1-h$. Then, by Equation \eqref{eq:CR_z}, we have
\begin{align*}
  \max_{h \in [\underline{h},\bar{h}]}\E_{\boldsymbol s}
    \left[\frac{\L\left(I, \boldsymbol s, z\right)}{\textsc{opt}(I)}\right] &\leq \max_{h \in [\underline{h},\bar{h}]}\E_{\boldsymbol s}
    \left[\frac{(n_1-z(s_1))r_1\boldsymbol{1}(n_1 \geq z(s_1))+(z(s_1)-n_1)r_2}{mr_2}\right]\, \\&= \frac{1}{m}\max_{h \in [\underline{h},\bar{h}]}\E_{\boldsymbol s}
    \left[(h-\frac{1}{p}s_1)\frac{r_1}{r_2}\boldsymbol{1}(h-\frac{1}{p}s_1 \geq 0)+(\frac{1}{p}s_1-h)\right]\,\\&=\frac{1}{m}\max_{h \in [\underline{h},\bar{h}]}\E_{\boldsymbol s}
    \left[(h-\frac{1}{p}s_1)\frac{r_1}{r_2}\boldsymbol{1}(h-\frac{1}{p}s_1 \geq 0)\right]\,\\&= \frac{r_1}{mpr_2}\max_{h \in [\underline{h},\bar{h}]}\E_{\boldsymbol s}
    \left[(hp-s_1)\boldsymbol{1}(hp-s_1 \geq 0)\right]\,.
\end{align*}
The second to last equation  follows because $\E_{\boldsymbol s}[\frac{1}{p}s_1-h]=0$. Recall that $s_1 \sim \text{Bin}(h, p)$. Next, we bound $\E_{\boldsymbol s}
    \left[(hp-s_1)\boldsymbol{1}(hp-s_1 \geq 0)\right]$. We will show that 
    \begin{align} \label{eq:one_side}
\max_{h \in [\underline{h},\bar{h}]}\E_{\boldsymbol s}
    \left[(hp-s_1)\boldsymbol{1}(hp-s_1 \geq 0)\right]\,\le\frac{1}{2}\sqrt{1-p}\sqrt{m}+\frac{\sqrt{2}\beta}{\sqrt{pm}},
\end{align}
where $\beta=0.4215\cdot \frac{p^2+(1-p)^2}{p(1-p)}$.
This confirms that  \begin{align*}
\max_{h \in\underline{h},\bar{h}]}\E_{\boldsymbol s}
    \left[\frac{\L\left(I, \boldsymbol s, z\right)}{\textsc{opt}(I)}\right] &\le \frac{r_1}{mpr_2}  \left(\frac{1}{2}\sqrt{1-p}\sqrt{m}+\frac{\sqrt{2}\beta}{\sqrt{pm}} \right)\\
    &= \frac{r_1\sqrt{1-p}}{r_2 p \sqrt{m}}  +\frac{\sqrt{2}\beta r_1}{p\sqrt{p}r_2} \frac{1}{m\sqrt{m}}\,,
    \end{align*}
which is the desired result.

It remains to show Equation \eqref{eq:one_side}. 
By  \cite{nagaev2011bound}, we have
\begin{equation} \label{eq:berryesseen}
    \big|\E_{\boldsymbol {s}}
    \left[(hp-s_1)\boldsymbol{1}(hp-s_1 \geq 0)\right]\,-\frac{1}{2}\E_{\boldsymbol {s}}
    \left[|hp-s_1|\right]\big| \le \frac{\beta}{\sqrt{h}}.
\end{equation}
In addition, \cite{berend2013sharp} provides a sharp estimate of $\E_{\boldsymbol s}\left[|hp-s_1|\right]\,$:
\begin{equation} \label{eq:sharp}
\frac{1}{\sqrt{2}}\sqrt{\frac{1-p}{p}h} \leq \E_{\boldsymbol s}\left[|hp-s_1|\right]\, \leq \sqrt{\frac{1-p}{p}h}.
\end{equation}
This leads to
\begin{align}
\E_{\boldsymbol s}
    \left[(hp-s_1)\boldsymbol{1}(hp-s_1 \geq 0)\right]\,&~\le~ 
    \frac{1}{2}\E_{\boldsymbol {s}}
    \left[|hp-s_1|\right] +  \frac{\beta}{\sqrt{h}}
   ~ \leq~ \frac{1}{2}\sqrt{\frac{1-p}{p}h}+   \frac{\beta}{\sqrt{h}}.
\end{align} 
Therefore, for large enough $m$ , we have
\begin{equation}
    \max_{h \in [\underline{h},\bar{h}]}\E_{\boldsymbol s}
    \left[(hp-s_1)\boldsymbol{1}(hp-s_1 \geq 0)\right]\,\leq \sqrt{\frac{1-p}{p}\bar{h}}+0.4215\frac{\beta}{\sqrt{\underline{h}}}=\frac{1}{2}\sqrt{1-p}\sqrt{m}+\frac{\sqrt{2}\beta}{\sqrt{pm}}\,,
\end{equation}
which is the desired result.

\textit{Second part: lower bound.} Here, we present a lower bound on $\max_{h \in [\underline{h},\bar{h}]}\E_{\boldsymbol s}
    \left[\frac{L\left(I, \boldsymbol s, z\right)}{\textsc{opt}(I)}\right]$, following similar steps in the first part. By definition, 
\begin{align*}
  \max_{h \in [\underline{h},\bar{h}]}\E_{\boldsymbol s}
    \left[\frac{\L\left(I, \boldsymbol s, z\right)}{\textsc{opt}(I)}\right]\,&\ge \max_{h \in [\underline{h},\bar{h}]}\E_{\boldsymbol s}
    \left[\frac{(n_1-z(s_1))r_1\boldsymbol{1}(n_1 \geq z(s_1))+(z(s_1)-n_1)r_2}{mr_1}\right]\, 
    \\&\ge  \frac{1}{m}\max_{h \in [\underline{h},\bar{h}]}\E_{\boldsymbol s}
    \left[(h-\frac{1}{p}s_1)\boldsymbol{1}(h-\frac{1}{p}s_1 \geq 0)+(\frac{1}{p}s_1-h)\right]\,\\
    &=\frac{1}{m}\max_{h \in [\underline{h},\bar{h}]}\E_{\boldsymbol s}
    \left[(h-\frac{1}{p}s_1)\frac{r_1}{r_2}\boldsymbol{1}(h-\frac{1}{p}s_1 \geq 0)\right]\,\\
    &= \frac{1}{mp}\max_{h \in [\underline{h},\bar{h}]}\E_{\boldsymbol s}
    \left[(hp-s_1)\boldsymbol{1}(hp-s_1 \leq 0)\right]\, \\&\ge 
    \frac{1}{mp}\max_{h \in [\underline{h},\bar{h}]}\left(\frac{1}{2}\E_{\boldsymbol {s}}
    \left[|hp-s_1|\right]\big| -  \frac{\beta}{\sqrt{h}} \right)\\&
    \ge \frac{1}{mp}\max_{h \in [\underline{h},\bar{h}]}\left(\frac{1}{2\sqrt{2}}\sqrt{\frac{1-p}{p}h}-   \frac{\beta}{\sqrt{h}} \right) \\&=\frac{1}{mp}\left(\frac{1}{2\sqrt{2}}\sqrt{(1-p)m}-   \frac{\beta}{\sqrt{pm}} \right) \\&=\frac{1}{2\sqrt{2}}\sqrt{(1-p)}\frac{1}{p\sqrt{m}}-   \frac{\beta}{p\sqrt{p}}\frac{1}{m\sqrt{m}}. 
\end{align*}
The second to last inequality is due to Equation \eqref{eq:berryesseen}. The last inequality is due to Equation \eqref{eq:sharp}.

\subsection{Proof of Lemma \ref{lem:other_mappings}}

Consider any arbitrary mapping $\widehat z(\cdot)$. By Equation \eqref{eq:CR_z}, we have
\begin{align*}
    \text{CR}(\widehat z) &= 1-\max_{h \in [\underline{h},\bar{h}]}\E_{\boldsymbol s}
    \left[\frac{\L\left(I, \boldsymbol s,\widehat z\right)}{\textsc{opt}(I)}\right]\, \\&\leq 1-\max_{h \in [\underline{h},\bar{h}]}\E_{\boldsymbol s}
    \left[\frac{(n_1-\widehat z(s_1))r_1\boldsymbol{1}(n_1 \geq \widehat z(s_1))+(\widehat z(s_1)-n_1)r_2}{mr_2}\right]\, \\&= 1-\frac{1}{m}\max_{h \in [\underline{h},\bar{h}]}\E_{\boldsymbol s}
    \left[(n_1-\widehat z(s_1))\frac{r_1}{r_2}\boldsymbol{1}(n_1 \geq \widehat z(s_1))+(\widehat z(s_1)-n_1)\right]\,.
\end{align*}

Let $\widehat y(s_1)=s_1+\widehat z(s_1)$ for any $s_1 \in [0, \bar h]$. Then, since $n_1+ s_1 = h$, we have 
\[
\text{CR} (\widehat z) ~\leq ~ 1-\frac{1}{m}\max_{h \in [\underline{h},\bar{h}]}\E_{\boldsymbol s}
    \left[(h-\widehat y(s_1))\frac{r_1}{r_2}\boldsymbol{1}(h \geq \widehat y(s_1))+(\widehat y(s_1)-h)\right]\,.
\]
Our goal to show that $\text{CR} (\widehat z) = 1- \Omega(\frac{1}{p\sqrt{m}})$. We show this by contradiction. Contrary to our claim, suppose that there exists a feasible mapping $y_0(\cdot)$ with $y_0(s_1) = s_1+ z_0(s_1)$ under which 
\begin{equation} \label{eq:cond3}
\lim_{m \to \infty}\max_{h \in [\underline{h},\bar{h}]}\frac{\E_{\boldsymbol s}
    \left[(h-y_0(s_1))\frac{r_1}{r_2}\boldsymbol{1}(h \geq y_0(s_1))+(y_0(s_1)-h)\right]\,}{p\sqrt{m}}=0\,.
\end{equation}
Note that if  Equation \eqref{eq:cond3} holds, we have $\text{CR}(z_0) = 1- o(\frac{1}{p\sqrt{m}})$.    
Now observe that $\big((h-y_0(s_1))\frac{r_1}{r_2}\boldsymbol{1}(h \geq y_0(s_1))\big) \geq 0$, and hence if   Equation \eqref{eq:cond3} holds, we must  have  
\begin{equation} \label{eq:cond8}
   \lim_{m \to \infty} \max_{h \in [\underline{h},\bar{h}]}\frac{\E_{\boldsymbol s }
    \left[y_0(s_1)-h\right]\,}{p\sqrt{m}} \leq 0.
\end{equation}
Next we consider two cases. In the first case, we assume that $ \lim_{m \to \infty} \max_{h \in [\underline{h},\bar{h}]}\frac{\E_{\boldsymbol s}
    \left[y_0(s_1)-h\right]\,}{p\sqrt{m}} < 0$ and we reach a contradiction under this assumption. In the second case, assume that  $ \lim_{m \to \infty} \max_{h \in [\underline{h},\bar{h}]}\frac{\E_{\boldsymbol s}
    \left[y_0(s_1)-h\right]\,}{p\sqrt{m}} = 0$ and again reach a contradiction. 

  
\textit{Case 1:  $ \lim_{m \to \infty} \max_{h \in [\underline{h},\bar{h}]}\frac{\E_{\boldsymbol s}
    \left[y_0(s_1)-h\right]\,}{p\sqrt{m}} < 0$.} Suppose that $
  \lim_{m \to \infty}\max_{h \in [\underline{h},\bar{h}]}\frac{\E_{\boldsymbol s}
    \left[y_0(s_1)-h\right]\,}{p\sqrt{m}}=\eta<0
  $. Then, we show that $y_0(\cdot)$ does not satisfy Equation \eqref{eq:cond3}, which is a contradiction. The left hand side of Equation \eqref{eq:cond3} can be written as
\begin{align*}
    \lim_{m \to \infty}\max_{h \in [\underline{h},\bar{h}]}& \frac{\E_{s_1 \sim\text{Bin}(h, p)}
    \left[(h-y_0(s_1))\frac{r_1}{r_2}\boldsymbol{1}(h \geq y_0(s_1))+(y_0(s_1)-h)\right]\,}{p\sqrt{m}} 
    \\
    &\geq \lim_{m \to \infty}\max_{h \in [\underline{h},\bar{h}]}\frac{
    \E_{s_1 \sim\text{Bin}(h, p)}
    \left[(h-y_0(s_1))\frac{r_1}{r_2}+(y_0(s_1)-h)\right]\,}{p\sqrt{m}} = -\frac{r_1}{r_2}\eta+\eta = -\frac{r_1-r_2}{r_2}\eta >0,
\end{align*}
which implies that Equation \eqref{eq:cond3} does not hold and hence it is a contradiction. 

\textit{Case 2:  $ \lim_{m \to \infty} \max_{h \in [\underline{h},\bar{h}]}\frac{\E_{\boldsymbol s}
    \left[y_0(s_1)-h\right]\,}{p\sqrt{m}} = 0$.}
Here, we compare $y_0(\cdot)$ and $y(\cdot)$ where we recall $y(s_1)= s_1+z(s_1)=\frac{1}{p} s_1$ for any $s_1 \in [0, \bar h]$.  We then show a contradiction. In Lemma \ref{lem:specific z}, we have shown that if $y(s_1)=\frac{1}{p}s_1$ (or equivalently $z(s_1)=\frac{1-p}{p}s_1$), then for any $h \in [\underline{h},\bar{h}]$, we have 
\begin{align} \notag
\lim_{m \to \infty}&\frac{\E_{s_1\sim \text{Bin}(h, p)}
    \left[(h-y(s_1))\frac{r_1}{r_2}\boldsymbol{1}(h \geq y(s_1))+(y(s_1)-h)\right]\,}{p\sqrt{m}}\\
    &= 
\lim_{m \to \infty}\frac{\E_{s_1\sim \text{Bin}(h, p)}
    \left[(h-\frac{s_1}{p})\frac{r_1}{r_2}\boldsymbol{1}(h \geq \frac{s_1}{p})+(\frac{s_1}{p}-h)\right]\,}{p\sqrt{m}} > 0.\label{eq:cond2}
\end{align}
   
Let Equation \eqref{eq:cond3} takes its maximum at $h=\hzero \in [\underline h, \bar h]$, and notice that $\hzero$ is a function with respect to $m$.  Also, note that Equation \eqref{eq:cond2} holds for any $h$, including $\hzero $. Hence, 
by subtracting Equation \eqref{eq:cond2} from Equation \eqref{eq:cond3} evaluated at $\hzero $, we have
\begin{align}  \notag
    \lim_{m \to \infty}\frac{1}{p\sqrt{m}}&\Bigg(\E_{s_1 \sim\text{Bin}(\hzero , p)}
    \left[(\hzero p -s_1)\frac{r_1}{r_2}\boldsymbol{1}(\hzero p  \geq s_1)+(s_1-\hzero p)\right]\,\\
    &-\E_{s_1 \sim \text{Bin}(\hzero , p)}
    \left[(\hzero p -y_0(s_1) p)\frac{r_1}{r_2}\boldsymbol{1}(\hzero  \geq y_0(s_1))+(y_0(s_1)-\hzero )\right]\,\Bigg)>0.\label{eq:cond4}
\end{align}
Now recall that under this case, 
we have  $ \lim_{m \to \infty} \frac{\E_{s_1\sim \text{Bin}(\hzero , p)}
    \left[y_0(s_1) p -\hzero p \right]\,}{p\sqrt{m}} = 0$ where $s_1 \sim \text{Bin}(\hzero , p)$. This implies that 
 $\lim_{m \to \infty}\frac{\E_{s_1\sim \text{Bin}(\hzero , p)}
    \left[s_1-y_0(s_1)p \right]\,}{\sqrt{m}p}= 0$, as $\E[s_1] = \hzero p$.  
Therefore, Equation \eqref{eq:cond4} can be written as 
\begin{equation} \label{eq:cond5}
    \lim_{m \to \infty} \frac{1}{p\sqrt{m}} \left( \E_{s_1 \sim \text{Bin}(\hzero , p)}
    \left[(\hzero p -s_1)\frac{r_1}{r_2}\boldsymbol{1}(\hzero p  \geq s_1)-(\hzero p -y_0(s_1) p)\frac{r_1}{r_2}\boldsymbol{1}(\hzero  \geq y_0(s_1))\right]\, \right)>0.
\end{equation}

Next, let sequence $q_d(m)=\frac{\vert y_0(d) p -d \vert}{p\sqrt{m}}$ for $d\in [0,\bar{h}] $, We define set $\mathcal D$ as follows 
\[
\mathcal D=\{d\in [0,\bar{h}] : \lim_{m \to \infty} q_{d}(m)>0 \}.
\]

By definition $q_d(m)$ is non-negative for all $d$, and hence we define the complement of set $\mathcal D$ as follows 
\begin{equation} \label{eq:def_D_c}
\mathcal D^{\mathsf{C}}=\{d\in [0,\bar{h}]: \lim_{m \to \infty} q_{d}=0 \}. 
\end{equation}

Then, we can write the inner term of Equation \eqref{eq:cond5} as
\begin{align*}
    \E_{s_1 \sim \text{Bin}(\hzero , p)}&
    \left[(\hzero p  -s_1)\frac{r_1}{r_2}\boldsymbol{1}(\hzero p   \geq s_1)-(\hzero p -y_0(s_1) p)\frac{r_1}{r_2}\boldsymbol{1}(\hzero  \geq y_0(s_1))\right]\,\\
    &=\sum_{d \in \mathcal D}\Pr(s_1=d)\left[(\hzero p -d)\frac{r_1}{r_2}\boldsymbol{1}(\hzero p \geq d)-(\hzero p -y_0(d)p)\frac{r_1}{r_2}\boldsymbol{1}(\hzero  \geq y_0(d))\right]\, \\&+ \sum_{d \in \mathcal D^{\mathsf{C}}}\Pr(s_1=d)\left[(\hzero p  -d)\frac{r_1}{r_2}\boldsymbol{1}(\hzero p \geq d)-(\hzero p -y_0(d) p)\frac{r_1}{r_2}\boldsymbol{1}(\hzero  \geq y_0(d))\right]\,.
\end{align*}

In Step (i), we show that  \begin{align}\lim_{m\rightarrow \infty}\frac{1}{p\sqrt{m}}\sum_{d \in \mathcal D^{\mathsf{C}}}\Pr(s_1=d)\left[(\hzero p -d)\frac{r_1}{r_2}\boldsymbol{1}(\hzero p  \geq d)-(\hzero p -y_0(d)p)\frac{r_1}{r_2}\boldsymbol{1}(\hzero  \geq y_0(d))\right]~=~ 0\,.\end{align} 

In Step (ii), we show that 
\[
\lim_{m \to \infty}\sum_{d \in \mathcal D}\Pr(s_1=d)>0.
\]
By these two equations, we can make a contradiction of Equation \eqref{eq:cond8}.

\textit{Step (i): }Observe that $\sum_{d \in \mathcal D^{\mathsf{C}}}\Pr(s_1=d) \leq 1$, and hence  we have
\begin{align} \notag
&\lim_{m \to \infty}\frac{1}{p\sqrt{m}}\sum_{d \in \mathcal D^{\mathsf{C}}}\Pr(s_1=d)\left[(\hzero p -d)\frac{r_1}{r_2}\boldsymbol{1}(\hzero p \geq d)-(\hzero -y_0(d))\frac{r_1}{r_2}\boldsymbol{1}(\hzero p  \geq y_0(d)p)\right]\,\\&\leq \lim_{m \to \infty}\frac{1}{p\sqrt{m}}\frac{r_1}{r_2}\max_{d \in \mathcal D^{\mathsf{C}}}\left[(\hzero p -d)\boldsymbol{1}(\hzero p \geq d)-(\hzero p -y_0(d) p)\boldsymbol{1}(\hzero   \geq y_0(d) )\right]\,. \label{eq:D_c} 
\end{align}

The inner term of the above equation can then be written as 
\begin{align*} \notag
    (\hzero p &-d)\boldsymbol{1}(\hzero p  \geq d)-(\hzero p -y_0(d)p)\boldsymbol{1}(\hzero p \geq y_0(d)p) \\&= \left((\hzero p -d)\boldsymbol{1}(\hzero p \geq d)-(\hzero p -y_0(d)p)\textbf{1}_{\hzero p \geq y_0(d)p}\right)\boldsymbol{1}(d<y_0(d) p)\\
    &+\left((\hzero p -d)\boldsymbol{1}(\hzero p \geq d)-(\hzero p -y_0(d) p)\boldsymbol{1}(\hzero p \geq y_0(d)p)\right)\boldsymbol{1}(d \geq y_0(d) p)\,.
\end{align*}

Notice that
\begin{align*}
    &\left( (\hzero p -d)\boldsymbol{1}(\hzero p \geq d)-(\hzero p -y_0(d)p)\boldsymbol{1}(\hzero  \geq y_0(d))\right) \boldsymbol{1}(d<y_0(d) p)\\
    &= \left((\hzero p-d-\hzero p +y_0(d)p)\boldsymbol{1}(\hzero p  \geq y_0(d)p)+(\hzero p -d)\boldsymbol{1}(d \leq \hzero p \leq y_0(d) p)\right)\boldsymbol{1}(d<y_0(d)p) \\&= \left((y_0(d) p-d)\boldsymbol{1}(\hzero p  \geq y_0(d)p)+(\hzero p -d)\boldsymbol{1}(d\leq \hzero p \leq y_0(d)p)\right)\boldsymbol{1}(d<y_0(d)p) \\&\leq \left(\vert y_0(d) p-d \vert+(y_0(d) p-d)\right)\boldsymbol{1}(d<y_0(d) p) \\&\leq 2\left(\vert y_0(d) p-d \vert\right)\boldsymbol{1}(d<y_0(d)p).
\end{align*}

Similarly, we have 
\[
\left((\hzero p -d)\boldsymbol{1}(\hzero p \geq d)-(\hzero p -y_0(d) p)\boldsymbol{1}(\hzero  \geq y_0(d))\right)\boldsymbol{1}(d \geq y_0(d)p) \leq 2\left(\vert y_0(d) p-d \vert\right)\boldsymbol{1}(d \geq y_0(d)p).
\]

Therefore, the inner term of Equation \eqref{eq:D_c} is upper bounded by
\begin{equation} \label{eq:cond17}
    (\hzero p-d)\boldsymbol{1}(\hzero p \geq d)-(\hzero p -y_0(d)p)\boldsymbol{1}(\hzero p \geq y_0(d)p) \leq 2\vert y_0(d) p-d \vert.
\end{equation}
Using this in Equation \eqref{eq:cond17} gives us
\begin{align*}
    \lim_{m \to \infty}&\frac{1}{p\sqrt{m}}\frac{r_1}{r_2}\max_{d \in \mathcal D^{\mathsf{C}}}\left[(\hzero p -d)\boldsymbol{1}(\hzero p  \geq d)-(\hzero p -y_0(d) p)\boldsymbol{1}(\hzero  \geq y_0(d))\right]\, \\&\leq \lim_{m \to \infty}\frac{1}{p\sqrt{m}}\frac{r_1}{r_2}\max_{d \in \mathcal D^{\mathsf{C}}}\left[2\vert y_0(d) p-d\vert\right]\, \\&= \lim_{m \to \infty}\frac{2r_1}{r_2}\max_{d \in \mathcal D^{\mathsf{C}}}\frac{\vert y_0(d) p -d \vert}{p\sqrt{m}}\,, \\&=0\,,
\end{align*}
which is what we wanted to show. Note that the last equation holds because of definition of $\mathcal D^{\mathsf{C}}$ in Equation \eqref{eq:def_D_c}.

\textit{Step (ii): } By Equation \eqref{eq:cond5} and the inequality above we have 
\begin{equation} \label{eq:cond9}
\lim_{m \to \infty}\frac{1}{p\sqrt{m}}\sum_{d \in \mathcal D}\Pr(s_1=d)\left[(\hzero p -d)\frac{r_1}{r_2}\boldsymbol{1}(\hzero p \geq d)-(\hzero p -y_0(d) p)\frac{r_1}{r_2}\boldsymbol{1}(\hzero  \geq y_0(d))\right]\,>0.
\end{equation}

Notice that for $d \in \mathcal{D}$, we have $\lim_{m \to \infty}\frac{|y_0(d) p-d|}{p\sqrt{m}}>0$. We define $c_1,c_2>0$ be the constant such that 
\[
\sup_{d \in \mathcal{D}}\lim_{m \to \infty}\frac{|y_0(d) p-d|}{p\sqrt{m}}=c_1\qquad \text{and} \qquad
\inf_{d \in \mathcal{D}}\lim_{m \to \infty}\frac{|y_0(d) p-d|}{p\sqrt{m}}=c_2.
\]
Recall that
$
\mathcal D=\{d\in \{0, 1, \ldots,\} : \lim_{m \to \infty} q_{d}(m)>0 \}
$, where $q_d(m)=\frac{\vert y_0(d) p-d \vert}{p\sqrt{m}}$ for $d\in \{0, 1, \ldots\} $.


By Equation \eqref{eq:cond17}, we have 
\[
(\hzero p -d)\frac{r_1}{r_2}\boldsymbol{1}(\hzero p  \geq d)-(\hzero p-y_0(d)p)\frac{r_1}{r_2}\boldsymbol{1}(\hzero  \geq y_0(d)) \leq 2\frac{r_1}{r_2}\vert y_0(d) p-d \vert.
\]

Therefore,  Equation \eqref{eq:cond9} can be upper bounded as 
\begin{equation} \label{eq:cond19}
0<\lim_{m \to \infty}\frac{1}{p\sqrt{m}}\sum_{d \in \mathcal D}\Pr(s_1=d)2\frac{r_1}{r_2}\vert y_0(d)p-d \vert 
\le 2c_1\frac{r_1}{r_2} \lim_{m \to \infty}\sum_{d \in \mathcal D}\Pr(s_1=d) ,
\end{equation}
which implies that 
\begin{equation} \label{eq:cond10}
\lim_{m \to \infty}\sum_{d \in \mathcal D}\Pr(s_1=d)>0.
\end{equation}

We are now ready to show the contraction. Recall that in the current case, we assumed that $\lim_{m \to \infty}\frac{\E_{\boldsymbol s}
    \big(y_0(s_1)-\hzero \big)}{p\sqrt{m}} =0$, which implies that $\lim_{m \to \infty}\frac{\E_{\boldsymbol s}
    \vert y_0(s_1)-\hzero \vert}{p\sqrt{m}} =0$. As in Steps (i) and (ii), we have shown that $\lim_{m\rightarrow \infty}\frac{1}{p\sqrt{m}}\sum_{d \in \mathcal D^{\mathsf{C}}}\Pr(s_1=d)\left[(\hzero p -d)\frac{r_1}{r_2}\boldsymbol{1}(\hzero p  \geq d)-(\hzero p -y_0(d)p)\frac{r_1}{r_2}\boldsymbol{1}(\hzero  \geq y_0(d))\right]~=~ 0\,$ and $\lim_{m \to \infty}\sum_{d \in \mathcal D}\Pr(s_1=d)>0$, by the following chain of equations we have the desired contradiction. 
    
\begin{align*}
    \lim_{m \to \infty}\frac{\E_{\boldsymbol s}
    \vert y_0(s_1)-\hzero \vert}{p\sqrt{m}}\, &= \lim_{m \to \infty} \frac{1}{p\sqrt{m}}\left(\sum_{d \in \mathcal D}\Pr(s_1=d)\vert y_0(d)p-\hzero p \vert+\sum_{d \in \mathcal D^{\mathsf{C}}}\Pr(s_1=d) \vert y_0(d) p-\hzero p \vert \right) \\&=  \lim_{m \to \infty} \frac{1}{p\sqrt{m}}\Bigg(\sum_{d \in \mathcal D}\Pr(s_1=d)\vert y_0(d) p-d+d-\hzero p \vert\\&+\sum_{d \in \mathcal D^{\mathsf{C}}}\Pr(s_1=d)\vert y_0(d) p-d+d-\hzero p \vert \Bigg) \\&=  \lim_{m \to \infty} \frac{1}{p\sqrt{m}}\left(\sum_{d \in \mathcal D}\Pr(s_1=d)\vert y_0(d) p-d+d-\hzero p \vert+\sum_{d \in \mathcal D^{\mathsf{C}}}\Pr(s_1=d)(d-\hzero p) \right) \\&\geq \lim_{m \to \infty} \frac{1}{\sqrt{m}}\left(\sum_{d \in [\hzero ]}(d-\hzero p)+\sum_{d \in \mathcal D}\Pr(s_1=d)\vert y_0(d) p-d \vert  \right) \\&= \E_{\boldsymbol s}
    \left[s_1-\hzero p\right]\,+\lim_{m \to \infty} \frac{1}{p\sqrt{m}}\sum_{d \in \mathcal D}\Pr(s_1=d)\vert y_0(d) p-d \vert \\&\geq \lim_{m \to \infty}\sum_{d \in \mathcal D}\Pr(s_1=d)c_2 \\&> 0\,.
\end{align*}
The first inequality  holds because $
\mathcal D^{\mathcal C}=\{d\in [0,\bar{h}] : \lim_{m \to \infty} q_{d}(m)=0 \}
$, where $q_d(m)=\frac{\vert y_0(d) p-d \vert}{p\sqrt{m}}$ for $d\in [0,\bar{h}] $. The third equation is due to definition of $\mathcal{D}^{\mathsf{C}}$. The fourth inequality is because $\vert y_0(d) p-d+d-\hzero p \vert \geq \vert y_0(d) p-d \vert + (d-\hzero p) $. The last inequality is from Equation \eqref{eq:cond10}. This is a contradiction to $ \lim_{m \to \infty} \max_{h \in [\underline{h},\bar{h}]}\frac{\E_{\boldsymbol s}
    \left[y_0(s_1)-h\right]\,}{p\sqrt{m}} = 0$. 

Put the two cases together, we conclude that there does not exist such mapping $y_0(\cdot)$ with which  the competitive ratio on family of $\mathbf{F}$ is $1-o(\frac{1}{p\sqrt{m}})$.

\section{Proof of Theorem \ref{thm:upperboundrandom}}\label{sec:proof:thm:upperboundrandom}
We construct the following random input distribution by taking advantage of input family $\textbf{F}$ that we used in the proof of Theorem \ref{thm:upperbound}. Recall that under a given input $h,\ell$ in  family $\textbf{F}$, $I_h$ is the random arrival sequence such that $n_2$ type $2$ agents arrive first followed by $n_1$ type $1$ agents, where $n_1 \sim \text{Bin}(h,1-p)$, $n_2 \sim \text{Bin}(\ell,1-p)$, and $\ell=\frac{10000\cdot m}{p}$, and $h\in [\underline h, \bar h]$ with $\underline{h}=\frac{1}{2}pm$ and  $\bar{h}=pm$.
 Then, in our random input distribution, we choose one of the feasible $h \in [\underline h, \bar h]$ uniformly at random; that is $\Pr(\text{Uniform}[\underline h, \bar h]=h) = \frac{1}{\bar{h}-\underline{h}+1}$ for any $h \in [\underline{h},\bar{h}]$.

Observe that due to our sampling procedure, even when the 
input $h\in [\underline h, \bar h]$ is realized, the online arrival sequence $I =I_h$ is still random. That prevents us from 
 applying the Von Neuman/Yao principle \cite{seiden2000guessing} to our setting. {At least in a footnote, clarify what the Von Neuman/Yao principle is and why we cannot use it in our setting.} Nonetheless, we derive a result similar to the Von Neuman/Yao principle that can be applied to our setting.

\begin{lemma} [Von Neuman/Yao principle: from deterministic to randomized inputs] \label{lem:yao} Consider a setting where the adversary chooses a \emph{meta input} $X\in\mathcal{X}$, and then based on the meta input $X$, a random \emph{input} $\widetilde X$ from a distribution $D(X)$ is realized.
For any random or deterministic algorithm $A$ and meta input $X \in \mathcal{X}$, let $C_A(\widetilde X) \in [0,1]$ be the (realized) reward of  algorithm $A$ under input $\widetilde X$ over the reward of the optimal in-hindsight algorithm that knows $\widetilde{X}$ in advance, where $\widetilde{X}\sim D(X)$. Then, $\E_{\widetilde X \sim D(X)}[C_A(\widetilde X)]$ is the competitive ratio of  algorithm $A$ under meta input $X$, where the expectation is with respect to $\tilde X$ and any randomness in algorithm $A$.  

Now, let $\mathcal A_d$ be the set of all deterministic algorithms.  
Let $P$ be a probability distribution over any deterministic algorithms.  Define $A\sim P$  as an algorithm chosen according to probability distribution $P$. \footnote{By \cite{yao1977probabilistic}, \cite{ball2009toward}, any randomized algorithm may
be viewed as a random choice $A\sim P$ among deterministic algorithms.} 
Further, let $Q$ be a probability distribution over the meta inputs in  $\mathcal {X}$, and let ${X}\sim Q$ denote a random meta input chosen according to $Q$. Then, for any distribution $P$ over  algorithms, distribution $Q$ over meta inputs, we have
\begin{align}
    \min_{X \in \mathcal{X}}~\E_{A\sim P}[\E_{\widetilde X\sim D(X)}[C_A(\widetilde X)]] ~\leq~ \max_{A \in \mathcal{A}_d}~\E_{X\sim Q}[\E_{\widetilde X\sim D(X)}[C_{A}(\widetilde{X}))]]. \label{eq:yao}
\end{align}
\end{lemma}

 The proof can be found in Section \ref{sebsec:yao}. We now ready to show the result using Lemma \ref{lem:yao}. 
Here, our meta inputs are $(h,\ell) \in \textbf{F}$ and our inputs are online arrival sequence $I_h$. Furthermore, the distribution over our meta inputs, i.e., $Q$ in Lemma \ref{lem:yao}, is uniform distribution. Recall that 
$\Pr(\text{Uniform}[\underline h, \bar h]=h) = \frac{1}{\bar{h}-\underline{h}+1}$ for any $h \in [\underline{h},\bar{h}]$.

Our goal here is to upper bound the competitive ratio of any randomized   algorithm on family $\mathbf F$, defined above. That is, for any distribution $P$ over any feasible  deterministic  algorithms $A$, we would like to bound the following quantity 
\[\text{CR}_{\text{rand}}:=\min_{\widehat{I}\in \mathbf F}~\E_{A\sim P}[\E_{I}[\text{CR}_{ A}(I)]  \,.\]
By Lemma \ref{lem:yao}, we have
\[\text{CR}_{\text{rand}} \le \max_{A \in \mathcal{A}_d}~\E_{\widehat I\sim Q}[\E_{I }[\text{CR}_{A}(I)]] = \max_{A \in \mathcal{A}_d}~\E_{\widehat I\sim Q}\left[\E_{I }\left[\frac{\rew_{A}(I)}{\textsc{opt}(I)}\right]\right] \,,\]
where $\mathcal{A}_d$ is the set of all deterministic algorithms. Since $Q$ is a uniform distribution over any 
$h\in [\underline h, \bar h]$, we have 
\begin{equation}
     \text{CR}_{\text{rand}} \leq \frac{1}{\bar{h}-\underline{h}+1}\cdot  \max_{A \in \mathcal{A}_d}\sum_{h=\underline{h}}^{\bar{h}}\E_{s_1 \sim \text{Bin}(h, p)}\left[\frac{\rew_{A}\left({I_h}, s_1\right)}{\textsc{opt}({I_h})}\right]\,,
\end{equation}
where $I_h$ is an arriving instance that a huge number of type $2$ agents arrive first, and followed by $h$ type $1$ agents. $\widehat I_h$ is a random variable which is consisted of $I_h$ with a uniform distribution on $h \in [\underline{h},\bar{h}]$.
 Recall that under $\widehat I_h$, $\rew_{A}$ does not depend on $s_2\gg m$, and hence we denote 
$\rew_{A}(I, \boldsymbol{s})$ with $\rew_{A}(I, s_1)$. 
We would like to show that for any deterministic algorithm $A \in \mathcal{A}_d$, we have
\[\frac{1}{\bar{h}-\underline{h}+1}\cdot  \sum_{h=\underline{h}}^{\bar{h}}\E_{s_1 \sim \text{Bin}(h, p)}\left[\frac{\rew_{A}\left({I_h}, s_1\right)}{\textsc{opt}({I_h})}\right] = 1-\Theta(\frac{1}{p\sqrt{m}})\,.
\]

By Equation \eqref{eq:cond3} in the proof of Theorem \ref{thm:upperbound}, it suffices to show that  there does not exist a mapping $y_0(\cdot): \{0,\ldots,\bar{h}\} \to \{0,1, \ldots, m\}$, such that 
\begin{equation} \label{eq:cond11_2}
    \lim_{m \to \infty}\frac{1}{\bar{h}-\underline{h}+1}\sum_{h=\underline{h}}^{\bar{h}}\frac{\E_{s_1 \sim \text{Bin}(h, p)}
    \left[(h-y_0(s_1))\frac{r_1}{r_2}\boldsymbol{1}(h \geq y_0(s_1))+(y_0(s_1)-h)\right]\,}{p\sqrt{m}}=0\,.
\end{equation}
Recall that under family $\mathbf F$, any deterministic algorithm can be presented with such a mapping.   

We show this by contradiction. Contrary to our claim, suppose there exists a mapping  $y_0(\cdot)$ such that Equation \eqref{eq:cond11_2} holds.  This and the fact that  $(h-y_0(s_1))\frac{r_1}{r_2}\boldsymbol{1}(h \geq y_0(s_1)) \geq 0$, imply that 
\begin{equation} \label{eq:cond12}
   \lim_{m \to \infty} \frac{1}{\bar{h}-\underline{h}+1}\sum_{h=\underline{h}}^{\bar{h}}\frac{\E_{s_1 \sim \text{Bin}(h, p)}
    \left[y_0(s_1)-h\right]\,}{p\sqrt{m}} \leq 0.
\end{equation}

Next we consider two cases. In the first case, we assume that 

$ \lim_{m \to \infty} \frac{1}{\bar{h}-\underline{h}+1}\sum_{h=\underline{h}}^{\bar{h}}\frac{\E_{s_1 \sim \text{Bin}(h, p)}
    \left[y_0(s_1)-h\right]\,}{p\sqrt{m}} < 0$ and we reach a contradiction under this assumption. In the second case, assume that  $ \lim_{m \to \infty} \frac{1}{\bar{h}-\underline{h}+1}\sum_{h=\underline{h}}^{\bar{h}}\frac{\E_{s_1 \sim \text{Bin}(h, p)}
    \left[y_0(s_1)-h\right]\,}{p\sqrt{m}}= 0$ and again reach a contradiction.

\textit{Case 1:  $ \lim_{m \to \infty} \frac{1}{\bar{h}-\underline{h}+1}\sum_{h=\underline{h}}^{\bar{h}}\frac{\E_{s_1 \sim \text{Bin}(h, p)}
    \left[y_0(s_1)-h\right]\,}{p\sqrt{m}} < 0$.}  Suppose that $
  \lim_{m \to \infty} \frac{1}{\bar{h}-\underline{h}+1}\sum_{h=\underline{h}}^{\bar{h}}\frac{\E_{s_1 \sim \text{Bin}(h, p)}
    \left[y_0(s_1)-h\right]\,}{p\sqrt{m}}=\eta<0
  $. Then, we show that $y_0$ does not satisfy Equation \eqref{eq:cond11_2}, which is a contradiction. The left hand side of Equation \eqref{eq:cond11_2} can be written as
\begin{align*}
    &\lim_{m \to \infty}\frac{1}{\bar{h}-\underline{h}+1}\sum_{h=\underline{h}}^{\bar{h}}\frac{\E_{s_1 \sim \text{Bin}(h, p)}
    \left[(h-y_0(s_1))\frac{r_1}{r_2}\boldsymbol{1}(h \geq y_0(s_1))+(y_0(s_1)-h)\right]\,}{p\sqrt{m}} \\&\geq \lim_{m \to \infty}\frac{1}{\bar{h}-\underline{h}+1}\sum_{h=\underline{h}}^{\bar{h}}\frac{\E_{s_1 \sim \text{Bin}(h, p)}
    \left[(h-y_0(s_1))\frac{r_1}{r_2}+(y_0(s_1)-h)\right]\,}{p\sqrt{m}} \\&= -\frac{r_1}{r_2}\eta+\eta \\&= -\frac{r_1-r_2}{r_2}\eta \\&>0,
\end{align*}
which contradicts  Equation \eqref{eq:cond11_2}.

\textit{Case 2:  $ \lim_{m \to \infty} \frac{1}{\bar{h}-\underline{h}}\sum_{h=\underline{h}}^{\bar{h}}\frac{\E_{s_1 \sim \text{Bin}(h, p)}
    \left[y_0(s_1)-h\right]\,}{p\sqrt{m}} = 0$.} Similar to the proof of Theorem \ref{thm:upperbound}, 
we compare $y_0(\cdot)$ and $y(\cdot)$, where $y(s_1)=\frac{s_1}{p}$ for any $s_1 \in \{0, 1, \ldots, h\}$, and then show a contradiction. Equation \eqref{eq:cond2} in the proof of Theorem \ref{thm:upperbound}  has shown that if $y(s_1)=\frac{s_1}{p}$, then for any $h \in [\underline{h},\bar{h}]$,
\begin{align}
\lim_{m \to \infty}\frac{\E_{s_1 \sim \text{Bin}(h, p)}
    \left[(h-y(s_1))\frac{r_1}{r_2}\boldsymbol{1}(h \geq \frac{s_1}{p})+(y(s_1)-h)\right]\,}{p\sqrt{m}} > 0. \label{eq:cond2_2}
\end{align}
   
By subtracting Equation  \eqref{eq:cond2_2} from Equation \eqref{eq:cond11_2}, we have
\begin{align} \notag 
    \lim_{m \to \infty}\frac{1}{p\sqrt{m}}\frac{1}{\bar{h}-\underline{h}}\sum_{h=\underline{h}}^{\bar{h}}\Bigg(\E_{\boldsymbol s}
    \left[(h p-s_1)\frac{r_1}{r_2}\boldsymbol{1}(hp \geq s_1)+(s_1-hp)\right]\,\\
    -\E_{\boldsymbol s}
    \left[(hp-y_0(s_1)p)\frac{r_1}{r_2}\boldsymbol{1}(hp \geq y_0(s_1)p)+(y_0(s_1)p-hp)\right]\,\Bigg)>0. \label{eq:cond14}
\end{align}

As in this case, we have $\lim_{m \to \infty}\frac{1}{\bar{h}-\underline{h}+1}\sum_{h=\underline{h}}^{\bar{h}}\frac{\E_{\boldsymbol s}
    \left[s_1-y_0(s_1)p \right]\,}{p\sqrt{m}}= 0$, Equation \eqref{eq:cond14} can be written as 
\begin{equation} \label{eq:cond15}
    \lim_{m \to \infty} \frac{1}{p\sqrt{m}} \left(\frac{1}{\bar{h}-\underline{h}+1}\sum_{h=\underline{h}}^{\bar{h}} \E_{\boldsymbol s}
    \left[(hp-s_1)\frac{r_1}{r_2}\boldsymbol{1}(hp \geq s_1)-(hp-y_0(s_1)p)\frac{r_1}{r_2}\boldsymbol{1}(hp \geq y_0(s_1)p)\right]\, \right)>0.
\end{equation}

Using the same set $\mathcal{D}$ defined in the proof of Theorem \ref{thm:upperbound} and following similar steps,  we can show that
\begin{equation}
\lim_{m \to \infty}\frac{1}{p\sqrt{m}} \left(\frac{1}{\bar{h}-\underline{h}+1}\sum_{h=\underline{h}}^{\bar{h}}\sum_{d \in \mathcal D^{\mathsf{C}}\cap [h]}\Pr(s_1=d)\left[(hp-d)\frac{r_1}{r_2}\boldsymbol{1}(hp \geq d)-(hp-y_0(d)p)\frac{r_1}{r_2}\boldsymbol{1}(hp \geq y_0(d)p)\right]\right)\,=0,
\end{equation}

Therefore, we have
\begin{equation} \label{eq:cond18}
\lim_{m \to \infty}\frac{1}{p\sqrt{m}} \left(\frac{1}{\bar{h}-\underline{h}+1}\sum_{h=\underline{h}}^{\bar{h}}\sum_{d \in \mathcal D\cap [h]}\Pr(s_1=d)\left[(hp-d)\frac{r_1}{r_2}\boldsymbol{1}(hp \geq d)-(hp-y_0(d)p)\frac{r_1}{r_2}\boldsymbol{1}(hp \geq y_0(d)p)\right]\right)\,>0.
\end{equation}
By the above inequality and following similar steps in the proof of Theorem \ref{thm:upperbound}, we can show that

\begin{equation} \label{eq:cond21}
\lim_{m \to \infty}\frac{1}{\bar{h}-\underline{h}+1}\sum_{h=\underline{h}}^{\bar{h}}\sum_{d \in \mathcal D\cap [h]}\Pr(s_1=d)>0.
\end{equation}

We are now ready to show the contradiction  using the assumption in this case, i.e., \[\lim_{m \to \infty}\frac{1}{\bar{h}-\underline{h}+1}\sum_{h=\underline{h}}^{\bar{h}}\frac{\E_{\boldsymbol s}
    \left[y_0(s_1)-h\right]\,}{p\sqrt{m}}=0\,.\] We have that
\begin{align*}
    &\lim_{m \to \infty}\frac{1}{\bar{h}-\underline{h}+1}\sum_{h=\underline{h}}^{\bar{h}}\frac{\E_{\boldsymbol s}
    \left[y_0(s_1)-h\right]\,}{p\sqrt{m}} \\&= \lim_{m \to \infty} \frac{1}{p\sqrt{m}}\frac{1}{\bar{h}-\underline{h}+1}\sum_{h=\underline{h}}^{\bar{h}}\left(\sum_{d \in \mathcal D\cap [h]}\Pr(s_1=d)(y_0(d)p-hp)+\sum_{d \in \mathcal D^{\mathsf{C}}\cap [h]}\Pr(s_1=d)(y_0(d)-h)\right) \\&=  \lim_{m \to \infty} \frac{1}{p\sqrt{m}}\frac{1}{\bar{h}-\underline{h}+1}\sum_{h=\underline{h}}^{\bar{h}}\left(\sum_{d \in \mathcal D \cap [h]}\Pr(s_1=d)(y_0(d)p-d+d-hp)+\sum_{d \in \mathcal D^{\mathsf{C}}\cap [h]}\Pr(s_1=d)(\frac{d}{p}-h) \right) \\&= \lim_{m \to \infty} \frac{1}{p\sqrt{m}}\frac{1}{\bar{h}-\underline{h}+1}\sum_{h=\underline{h}}^{\bar{h}}\left(\sum_{d \in [h]}(d-hp)+\sum_{d \in \mathcal D \cap [h]}\Pr(s_1=d)(y_0(d)p-d) \right) \\&= \lim_{m\rightarrow \infty} \frac{1}{\bar{h}-\underline{h}+1}\sum_{h=\underline{h}}^{\bar{h}}\E_{d \sim \text{Bin}(h,p)}
    \left[d-hp\right]\,+\lim_{m \to \infty} \frac{1}{p\sqrt{m}}\frac{1}{\bar{h}-\underline{h}+1}\sum_{h=\underline{h}}^{\bar{h}}\sum_{d \in \mathcal D \cap [h]}\Pr(s_1=d)(y_0(d)p-d) \\&=\lim_{m \to \infty} \frac{1}{p\sqrt{m}}\frac{1}{\bar{h}-\underline{h}+1}\sum_{h=\underline{h}}^{\bar{h}}\sum_{d \in \mathcal D \cap [h]}\Pr(s_1=d)(y_0(d)p-d) \\&\geq \lim_{m \to \infty} \frac{1}{\bar{h}-\underline{h}+1}\sum_{h=\underline{h}}^{\bar{h}}\sum_{d \in \mathcal D \cap [h]}\Pr(s_1=d)c_2\\&> 0,
\end{align*}
The last inequality is from Equation \eqref{eq:cond21}. This is a contradiction to the assumption in case 2.  Therefore, there does not exist a mapping $y_0(\cdot)$ with which we can achieve a competitive ratio of  $1-o(\frac{1}{p\sqrt{m}})$.

\subsection{Proof of Lemma \ref{lem:yao}} \label{sebsec:yao}  Let $P(A)$ be the probability of  choosing an algorithm $A$. Further, let $Q(X)$ be the probability of choosing meta input $X\in \mathcal{X}$. We start from the right hand side of Equation \eqref{eq:yao}: 

\begin{align*}
    \max_{A \in \mathcal{A}_d}~\E_{X\sim Q}[\E_{\widetilde X\sim D(X)}[C_A(\widetilde{X}))]] &= \sum_{A \in \mathcal{A}_d}P(A) \cdot 
    \max_{A \in \mathcal{A}_d}~\E_{X\sim Q}[\E_{\widetilde X\sim D(X)}[C_A(\widetilde{X}))]]
    \\&\geq  
    \sum_{A \in \mathcal{A}_d}P(A) \cdot \E_{X\sim Q}[\E_{\widetilde X\sim D(X)}[C_A(\widetilde{X}))]]
    \\&=
    \sum_{A \in \mathcal{A}_d}P(A) \cdot \sum_{X\in \mathcal X} Q(X)\cdot \E_{\widetilde X\sim D(X)}[C_A(\widetilde{X}))]
    \\&=
      \sum_{X \in \mathcal{X}}Q(X) \cdot \sum_{A\in \mathcal{A}_d} P(A)\cdot \E_{\widetilde X\sim D(X)}[C_A(\widetilde{X}))]
    \\&= \sum_{X \in \mathcal{X}}Q(X)\cdot \E_{A\sim P}[\E_{\widetilde X\sim D(X)}[C_A(\widetilde X)]] \\
    &~\ge~ \min_{X \in \mathcal{X}}~\E_{A\sim P}[\E_{\widetilde X\sim D(X)}[C_A(\widetilde X)]]\,,
\end{align*}
where the last inequality is the desired result.

\section{Proof of Theorem \ref{thm:uppersmallp}} \label{subsec:proofoflowp}
We define a family of arriving sequence $\mathbf{F}_1$ as follows: The input family $\mathbf{F}_1$ contains all $(h,\ell)$ such that $h \in [0,m/p]$ and $\ell=\frac{10000\cdot m}{p}$. For any $h \in [0,m/p]$, we then denote $I_h$ as an ordered random arrival sequence under which $n_2 \sim \text{Bin}(\ell,1-p)$ type $2$ agents arrive followed by $n_1 \sim \text{Bin}(h,1-p)$ type $1$ agents.
Then, we characterize an upper bound on the CR of any deterministic algorithm under the family $\mathbf{F}_1$. Similar to the proof of Theorem \ref{thm:upperbound}, given that the online arrival sequences $I_h$ are all ordered,  any deterministic algorithm has to decide about how many type $2$ agents they accept provided that they observe $s_1$ samples from type $1$ agents. Put differently, any deterministic algorithm can be represented by a mapping that maps $s_1$ to the number of type $2$ agents it accepts. 

We pick two sequences $I_0$ and $I_{\tilde{h}}$ inside the family $\mathbf{F}_1$,
where under $I_0$, no type $1$ agents arrive and under  $I_{\tilde{h}}$, the market size for type $1$ agents is $\tilde{h}$.
Here, $\tilde{h}$ is the smallest integer such that $\tilde{h}p+\sqrt{\tilde{h} \log(\tilde{h})} \geq 1$. It is easy to see that $\tilde{h}\le 1/p$ and $\tilde{h}=O(1/p)$. 

For the sequence $I_0$, $\Pr(s_1=0)=1$ because $h=0$. For the sequence $I_{\tilde{h}}$, $\Pr(s_1=0)\geq 1-\frac{1}{\tilde{h}^2}$ because by Hoeffding's inequality,
\begin{equation} \label{eq:concentration1}
\Pr\left(s_1 - \E[s_1] \geq \sqrt{\tilde{h} \log(\tilde{h})}\right) \leq \frac{1}{\tilde{h}^2}.
\end{equation}
As $\E[s_1]=\tilde{h}p$, Equation \eqref{eq:concentration1} is equivalent to 
\[
\Pr\left(s_1 - \tilde{h}p  \geq \sqrt{\tilde{h} \log(\tilde{h})}\right) ~=~\Pr(s_1  \geq \tilde{h}p+\sqrt{\tilde{h} \log(\tilde{h})}) ~\leq~  \Pr(s_1 \geq 1) \leq \frac{1}{\tilde{h}^2},
\]
which implies that 
\[
\Pr(s_1 =0) \geq 1-\frac{1}{\tilde{h}^2}\,.
\]

Now suppose that an algorithm accepts $m-x$ type $2$ agents if the decision-maker observes $s_1=0$. Then the CR under $I_0$ is $\frac{m-x}{m}=1-\frac{x}{m}$. The expected CR under $I_{\tilde{h}}$ (denoted by $\text{CR}(I_{\tilde{h}})$) is bounded by 
\begin{align*}
\text{CR}(I_{\tilde{h}}) &= \E[\text{CR}(I_{\tilde{h}}) \vert s_1=0]\Pr(s_1=0) + \E[\text{CR}(I_{\tilde{h}}) \vert s_1\geq 1]\Pr(s_1 \geq 1) \\&\leq \frac{(m-x)r_2+xr_1}{\min\{\tilde{h},m\}r_1+(m-\tilde{h})^{+}r_2}+\frac{1}{\tilde{h}^2}\,,
\end{align*}

where we use the fact that $\Pr(s_1=0) \leq 1$ and $\E[\text{CR}(I_{\tilde{h}}) \vert s_1\geq 1] \leq 1$. Therefore, the CR is bounded by 
\[
\min\{1-\frac{x}{m}, \frac{(m-x)r_2+xr_1}{\min\{\tilde{h},m\}r_1+(m-\tilde{h})^{+}r_2}+\frac{1}{\tilde{h}^2} \}.
\]
To find the optimal $x$, we can let the two terms above be equal, and we have
\[
x^{*}=m-\frac{m^2}{m+\min\{\tilde{h},m\}(1-\alpha)}-m\frac{1}{\tilde{h}^2},
\]
and the upper bound of CR can be calculated as 
\[
1-\frac{x^{*}}{m} \leq \frac{m}{m+\min\{\tilde{h},m\}(1-\alpha)}+\frac{1}{\tilde{h}^2}\,,
\]
which is the desired result.

The proof for randomized algorithms is omitted because we can apply Lemma \ref{lem:yao} and use the similar method in the proof of Theorem \ref{thm:3} to derive the upper bound among all randomized algorithms.

\section{Proof of Proposition \ref{prop:limitedm}}
\label{sec:prooflimitm}
The proof is  similar to the proof of Theorem \ref{thm:2}, where there is no capacity constraint during the test period. When there is a capacity constraint (i.e., when the number of resources during the test period is $m_t=\omega(\sqrt{m})$), the number of samples is not longer $s_1$ and $s_2$ under the case of $s_1+s_2> m_t$. This, in turn, influences the  estimation of $\widehat{r_i}$, $i \in \{1,2\}$, and the probability of the good event that $\widehat{r_1} > \widehat{r_2}$. 

We note that  in the proof of  Theorem \ref{thm:2}, the only place we needed a (lower) bound on the probability of the good event is in the analysis of the third region $\mathcal R_3 =\{(h,\ell):  h\ge \mzero , \ell \ge \mzero  \}$. The results for the other two regions are still valid. 
Hence, in the following, we focus on the third region.

Given that $h,\ell \geq \sqrt{m}$ in the third region and $s_1+s_2>m_t$, we split the analysis into two cases: 
\begin{enumerate}
    \item $s_1 > m_t/2$, $s_2 > m_t/2$. In this case, the decision-maker accepts $m_t/2$ each type agents. Then, $\widehat{r_i}=\frac{\text{Bin}(m_t/2,r_i)}{m_t/2}$, $i \in \{1,2\}$. Then, we replace $s_i$ by $m_t/2$ in the proof of Lemma \ref{lem:cr3largem} and as $m_t=\omega(\sqrt{m})$, we can get the same asymptotic result.
    
    \item $s_1 < m_t/2$, $s_2 > m_t/2$ or $s_1 > m_t/2$, $s_2 < m_t/2$. Without loss of generality,  we analyze the first case (i.e., $s_1 < m_t/2$, $s_2 > m_t/2$). In this case, we accept $s_1$ type $1$ agents and $m_t-s_1$ type $2$ agents. Then, $\widehat{r_1}$ remains the same and $\widehat{r_2}=\frac{\text{Bin}(m_t-s_1,r_2)}{m_t-s_1}$. Then, we replace $s_2$ by $m_t-s_1$ in the proof of Lemma \ref{lem:cr3largem} and as $m_t=\omega(\sqrt{m})$, we can also get the same asymptotic result.
\end{enumerate}

\section{Proof of Proposition \ref{prop:negativereal}} \label{sec:proof:prop:negativereal}
Consider the instance $I$ which contains $\ell = \frac{10000\cdot m}{p\cdot r_2}$ type $2$ agents and $0$ type $1$ agent. Then, in the online arrival sequence, with probability $1$, the number of  realized rewards of  $1$ is larger than $m$. Therefore, $\textsc{opt-re}(\mathcal R , I)=m$ almost surely. 

As the realization of each arriving agent is independent and follows $\text{Ber}(r_2)$, any algorithm cannot predict the realization of the arriving agent. If one accepts an agent, then the expectation of reward we get is $1 \cdot \Pr(\text{Ber}(r_2)=1)+0 \cdot \Pr(\text{Ber}(r_2)=0)=r_2$. The best algorithm can accept $m$ type $2$ agents due to the resource constraint.  Therefore, the realized  competitive ratio is at most $\frac{mr_2}{m}=r_2$.

\section{Proof of Proposition \ref{thm:heterop}} \label{app:proofheterop}
We first consider region $\mathcal{R}_1$. From Lemma \ref{lem:c1}, CR is $1-\Theta(1/\sqrt{m})$ which does not depend on $p$ in all cases except $\{n_2<m-n_1\}$ or $\{n_2 \geq m-n_1\} \cap \{n_1 \leq \frac{1-p_1}{p_1}s_1 \}$. Now, suppose that we apply the modified version of Algorithm \ref{alg:ptas}.This means  we give protection level $\min\{\frac{1-p_1}{p_1}s_1,m  \}$ to type $1$ agent when the algorithm believes type $1$ is the high-reward type  and give protection level $\min\{\frac{1-p_2}{p_2}s_2,m  \}$ to type $2$ agent when the algorithm believes type $2$ is the high-reward type. As the CR does not depend on $p$, in all cases except $\{n_2<m-n_1\}$ or $\{n_2 \geq m-n_1\} \cap \{n_1 \leq \frac{1-p_1}{p_1}s_1 \}$, the CR is still $1-\Theta(1/\sqrt{m})$. 

For the case of $n_2<m-n_1$,  we have
              \begin{align*}
            \frac{\rewGood}{\textsc{opt}(\mathbf{n})}&=\frac{\min\left\{n_2,\left(m-\frac{1-p_1}{p_1}s_1\right)^{+}\right\}\cdot r_2+n_1 r_1}{n_2r_2+n_1r_1}\\&\geq \frac{\min\left\{n_2,\left(m-\frac{1-p_1}{p_1}s_1\right)^{+}\right\}}{n_2} = \min\left\{1,\frac{\left(m-\frac{1-p_1}{p_1}s_1\right)^{+}}{n_2} \right\} \\& = \frac{\left(m-\frac{1-p_1}{p_1}s_1\right)^{+}}{n_2}
            > \frac{\left(m-\frac{1-p_1}{p_1}s_1\right)^{+}}{m-n_1} \geq \frac{\left(m-\frac{1-p_1}{p_1}h\right)^{+}}{m} \geq \left( 1-\frac{1-p_1}{p_1}\frac{1}{\mzero} \right)^{+}\,,
        \end{align*}
        where in the last  inequality holds because in region $\mathcal R_1$, $h< \sqrt{m}$.

For the case $\{n_2 \geq m-n_1\} \cap \{n_1 \leq \frac{1-p_1}{p_1}s_1 \}$, we have
\begin{align*}
             \frac{\rewGood}{\textsc{opt}(\mathbf{n})}&=\frac{\min\{n_2,(m-\frac{1-p_1}{p_1}s_1)^{+}\}r_2+n_1r_1}{(m-n_1)r_2+n_1r_1} \\&=\frac{(m-\frac{1-p_1}{p_1}s_1)^{+}r_2+n_1r_1}{(m-n_1)r_2+n_1r_1} \geq \frac{(m-\frac{1-p_1}{p_1}s_1)^{+}}{m-n_1} \geq \frac{(m-\frac{1-p_1}{p_1}h)^{+}}{m} \geq \left(1-\frac{1-p_1}{p_1}\frac{1}{\mzero}\right)^{+}.
\end{align*}

For any other cases, we still have the CR is $1-\Theta(1/\sqrt{m})$. Totally, we have the modified version of Algorithm \ref{alg:ptas} in $\mathcal{R}_1$ has a CR of $1-\Theta(1/(p_1\sqrt{m}))$. 

Next, from Lemma \ref{lem:c3}, in region $\mathcal{R}_2$, CR is $1-\Theta(1/\sqrt{m})$, which does not depend on $p$. Therefore, the modified version also achieves a CR of $1-\Theta(1/\sqrt{m})$.

Finally, from Lemma \ref{thm:c4}, in region $\mathcal{R}_3$, for $m \geq m_1$, $\overline{\text{CR}_3}=1-\Theta(1/(p^2m))$; and for $m<m_1$, $\underline{\text{CR}_3}=1-\Theta(1/\sqrt{m})$. Again, if we apply the modified version of Algorithm \ref{alg:ptas}, for $m<m_1$, the CR is still $1-\Theta(1/\sqrt{m})$. Otherwise, from the proof of Lemma \ref{lem:cr3largem}, we have
\[
\overline{\text{CR}_3} \geq (1-\frac{1}{(\mzero p_1-m^{1/4})^2})(1-\frac{1}{(\mzero p_2-m^{1/4})^2}) = 1-\Theta(\frac{1}{\min\{p_1,p_2\}^2m}).
\]

To summarize, we have the modified version of Algorithm \ref{alg:ptas} has a CR of
\[
1-\Theta(\max\{\frac{1}{p_1\sqrt{m}}, \frac{1}{\min\{p_1,p_2\}^2m} \}).
\]

\section{Algorithm for Sampling Probability $p=O(1/\sqrt{m})$}  \label{sec:algpsmall}

\begin{algorithm}[htbp] 
\caption{Algorithm 3 for $p=O(1/\sqrt{m})$}
\begin{itemize}
\item[] \textbf{Input:} The number of resources $m$.

\end{itemize}
\begin{enumerate}
  \item With probability $1/2$, we only accept  type $1$ agents. With probability $1/2$,  
 we only accept  type $2$ agents.
\end{enumerate}

\label{alg:naive}
\end{algorithm}

 \begin{theorem}[Competitive Ratio of Algorithm \ref{alg:naive}] \label{thm:naive1}
Consider the model presented in Section \ref{sec:model}, where the expected rewards of $r_i$, $i\in \{1, 2\}$ is unknown to the decision-maker. For $p=O(1/\sqrt{m})$ and for any $m > 0$, the CR of Algorithm \ref{alg:naive} is $1/2$.
\end{theorem}

\proof{Proof of Theorem \ref{thm:naive1}}
With probability $0.5$, we reject all type $2$ agents and accept all type $1$ agents, and with probability $0.5$, we reject all type $1$ agents and accept all type $2$ agents. Then, we have
\[
\frac{\E[\rew_{A}\left(I\right)]}{\textsc{opt}(I)} = \frac{\frac{1}{2}\min\{n_1,m\}r_1+\frac{1}{2}\min\{n_2,m\}r_2}{\min\{n_1,m\}r_1+\min\{n_2,(m-n_1)^{+}\}r_2} \geq \frac{1}{2},
\]
where the last inequality is due to $\min\{n_2,(m-n_1)^{+}\} \leq \min\{n_2,m\}$.

\Halmos
\endproof

\section{Partially Known Sampling Probability} \label{subsec:partialknown}
In our setting, we assume that the sampling probability $p$ is known to the decision-maker, and we show that under this assumption, Algorithm \ref{alg:ptas} obtains an asymptotically optimal CR.  However, as we show in    Section \ref{sec:negative}, obtaining the same result is not possible when $p$ is completely unknown, justifying our assumption. Nonetheless, here we would like to study a setting in which  the decision-maker (before the test period starts)  obtains an estimate of the sampling probability, denoted $\hat{p}$, where $\hat{p} \in [p(1-\delta),p(1+\delta)]$ for some $\delta \in (0,1)$. Here, $p$ is the true (unknown) sampling probability and $\delta$ measures how accurate the estimate is. Such an estimate can be obtained using previous allocation periods in which the same outreach program is used.
 
 For this setting, let us consider a variation of Algorithm \ref{alg:ptas} under which the algorithm pretends the estimated sampling probability to be the true probability. 
That is, when the algorithm protects type $1$ agents (i.e., $\widehat r_1 > \widehat r_2$), the algorithm   gives a protection level of  $ \min\{m,s_1\frac{1-\hat{p}}{\hat{p}} \} $  to type $1$ agent. Otherwise, the algorithms gives $\min\{m,s_2\frac{1-\hat{p}}{\hat{p}} \}$ protection level to type $2$ agent. The following example shows that  Algorithm \ref{alg:ptas} is quite robust to the estimation error in the sampling probability under adversarial arrivals.  Recall that in Section \ref{sec:simulation}, we show that Algorithm \ref{alg:ptas} with noisy $p$ performs extremely good under stochastic arrivals. 

\begin{example}
Let $\text{CR}(p;\hat{p})$ be the worst case CR of Algorithm \ref{alg:ptas} with input $\hat{p}$ and the true sampling probability $p$. Then, we measure the robustness of Algorithm \ref{alg:ptas} using the following measure that we refer to as \emph{robustness ratio}: 
\[
\sup_{\hat{p} \in [p(1-\delta),p(1+\delta)]} \frac{\text{CR}(p; p)-\text{CR}(p; \hat{p})}{\text{CR}(p; p)}.
\]
Observe that the smaller the robustness ratio, the better. The CR is calculated by taking the minimum value among $5000$ generated arrival instances.

Figures \ref{fig:figrobmsmall} and  \ref{fig:sigmalarge} present the robustness ratio versus $p$
with $\delta \in \{0.1, 0.2\}$, $r_1 =0.9$, and $r_2= 0.5$  for the case of $m =30$ and $m= 200$, respectively. The figures show the robustness ratio mainly depends on the estimation error $\delta$. When $\delta =0.1$, the robustness ratio is in the range of $1\%-6\%$. That is, the CR can drop by at most $6\%$ when we have 10\% estimation error in $p$. Similarly,  when with $\delta =0.2$, the robustness ratio is in the range of $1-12\%$. These results confirm the robustness of Algorithm \ref{alg:ptas} to the estimation error in $p$.

\begin{figure}[ht] 
  \centering
  \includegraphics[width=.5\linewidth]{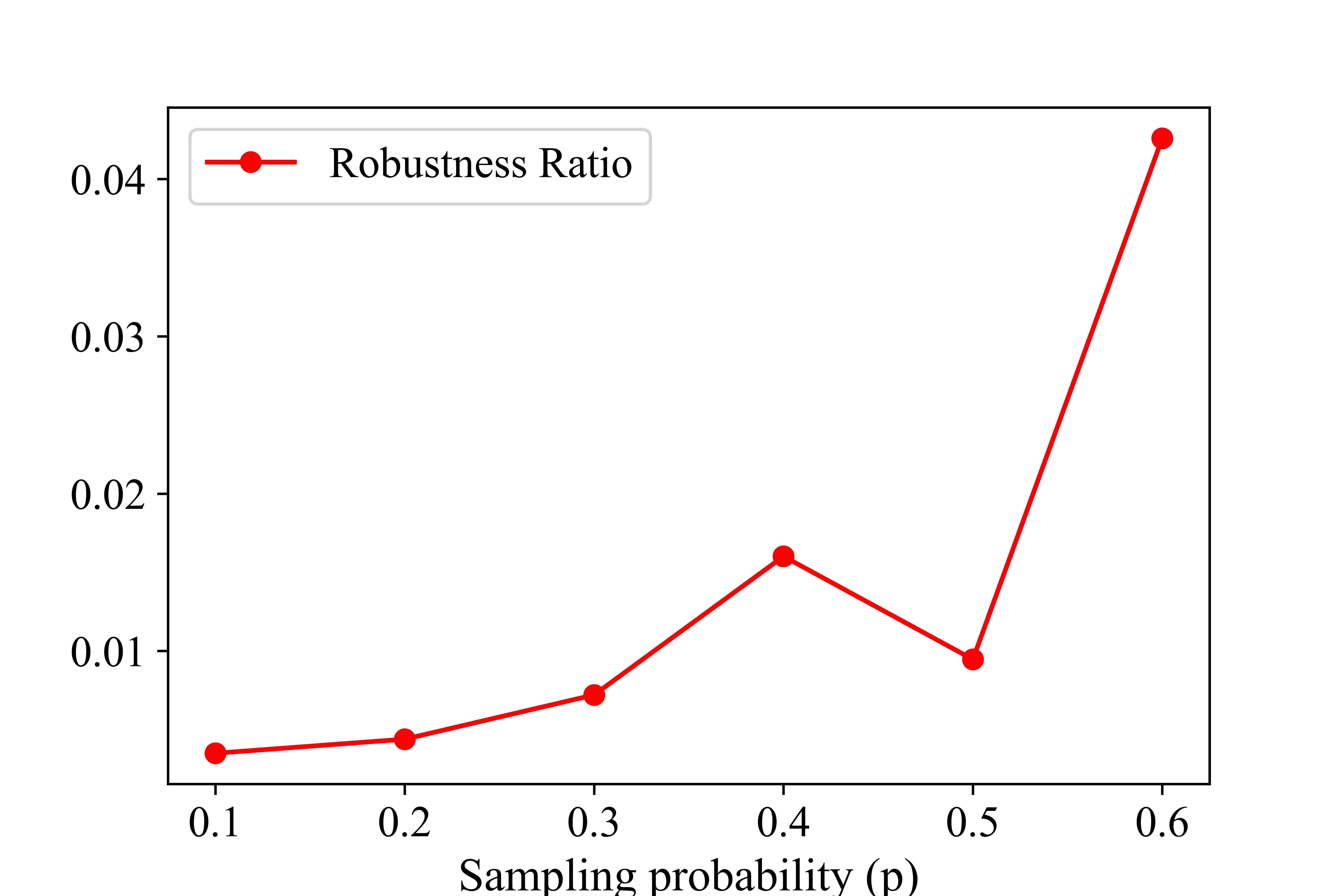}  
  \caption{ $\delta = 0.1$, $m=30$}
  \label{fig:sigmasmallmsmall}
  \centering
  \includegraphics[width=.5\linewidth]{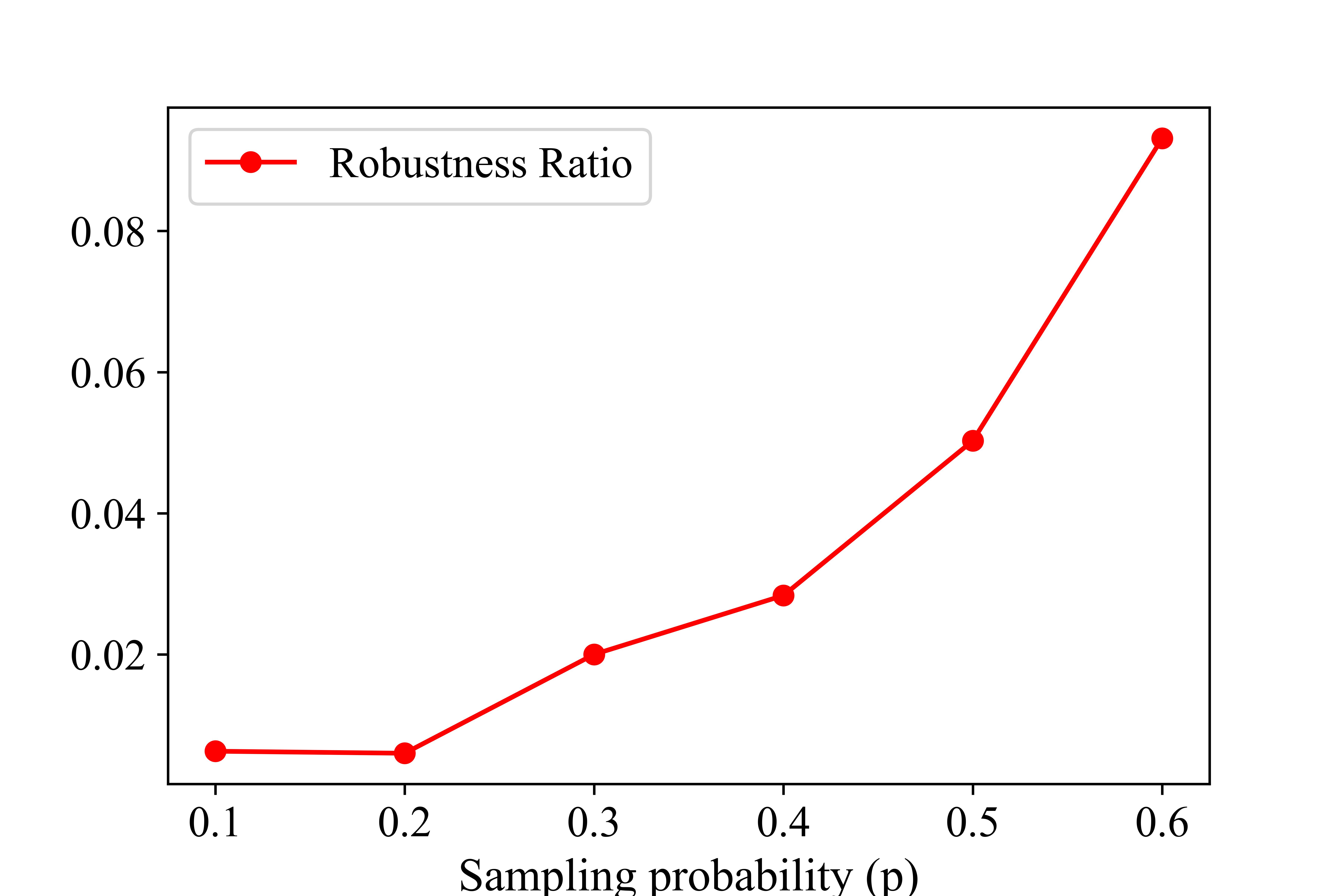}  
  \caption{ $\delta = 0.2$, $m=30$}
  \label{fig:sigmalargemsmall}
\caption{Robustness of Algorithm \ref{alg:ptas} to the estimation error in $p$ for $m=30$. Here, $r_1=0.9$, and  $r_2=0.5$.}
\label{fig:figrobmsmall}
\end{figure}

\begin{figure}[ht] 
  \centering
  \includegraphics[width=.5\linewidth]{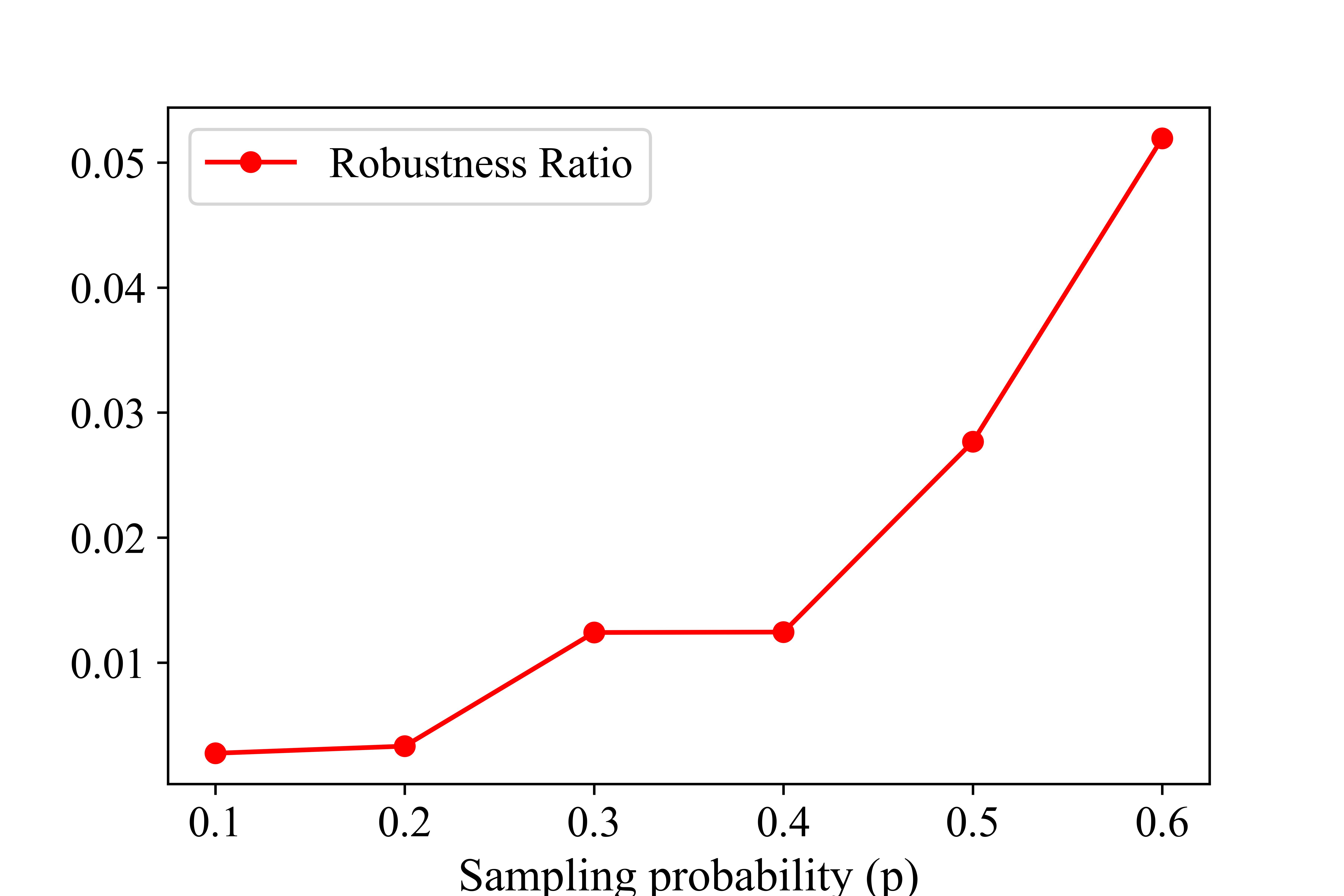}  
  \caption{ $\delta = 0.1$, $m=200$}
  \label{fig:sigmasmall}
  \centering
  \includegraphics[width=.5\linewidth]{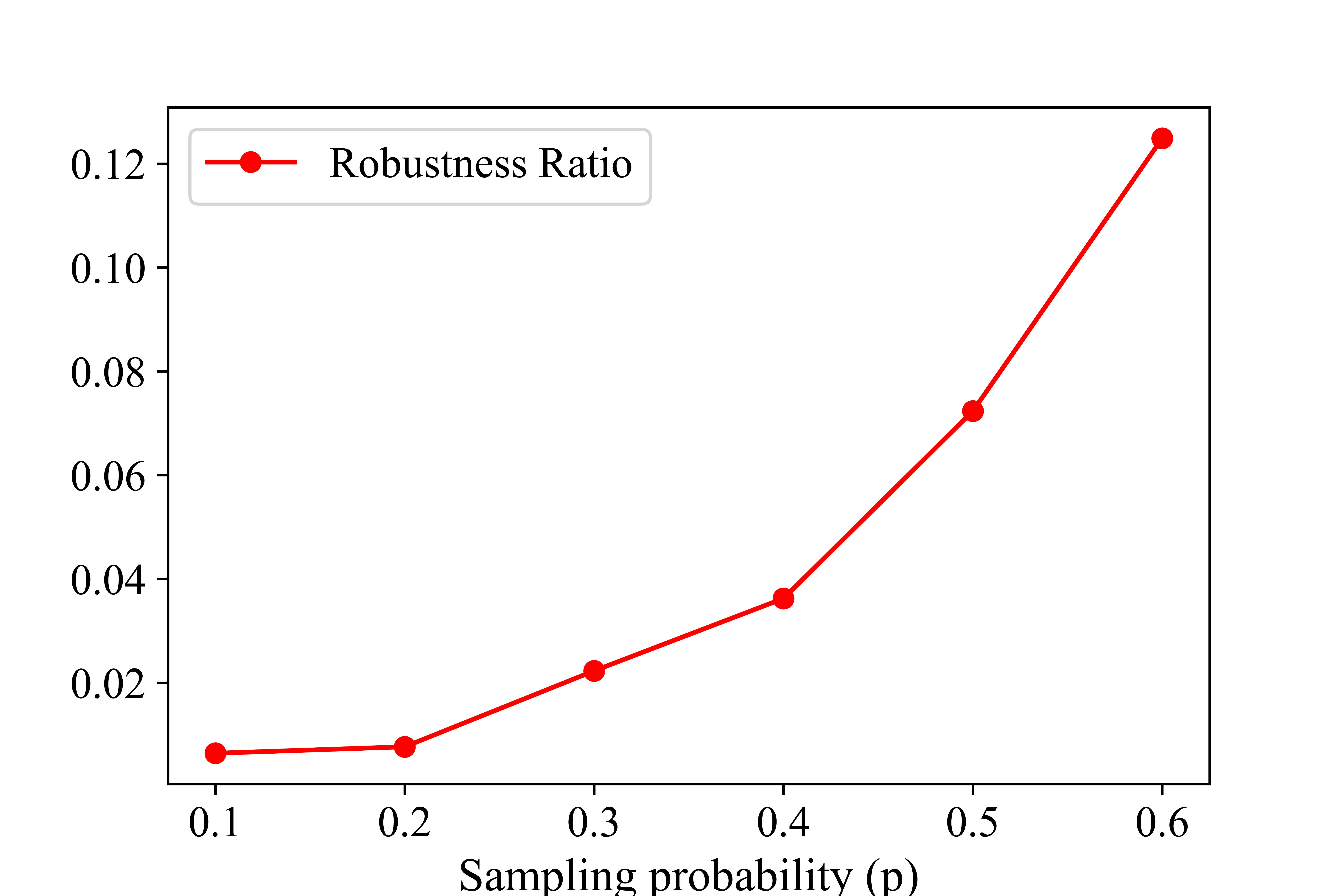}  
  \caption{ $\delta = 0.2$, $m=200$}
  \label{fig:sigmalarge}
\caption{Robustness ratio of Algorithm \ref{alg:ptas} in $p$ for $m=200$. Here, $r_1=0.9$, and  $r_2=0.5$. }
\label{fig:figrob}
\end{figure}

\end{example}

\end{APPENDICES}

\end{document}